%% file: LDM_for_MLECS_Elsarticle.tex
\documentclass[preprint,12pt]{elsarticle}

\usepackage[a4paper, total={6.6in, 9in}]{geometry}
\usepackage{lipsum}
\usepackage{amsfonts}
\usepackage{graphicx}
\usepackage{epstopdf}
\usepackage{tikz,graphicx}
\usetikzlibrary{3d,perspective} 
\usetikzlibrary{backgrounds,calc,shadows,positioning,fit}
\usepackage{amssymb}
\usepackage{amsmath}
\usepackage{enumerate} 
\usepackage{enumitem}
\usepackage[colorlinks=true,linkcolor=black,citecolor=blue,urlcolor=blue,]{hyperref}
\usepackage{algorithm}
\allowdisplaybreaks[4] 

\usepackage{amsthm}
\newtheorem{remark}{Remark}

\newtheorem{problem}{Problem}
\newtheorem{consequence}{Consequence}
\newtheorem{assumption}{Assumption}
\newtheorem{proposition}{Proposition}
\newtheorem{theorem}{Theorem}
\newtheorem{definition}{Definition}

\begin{document}

\begin{frontmatter}



\title{A layer decomposition method for multi-layer elastic contact systems with interlayer Tresca friction}


\author[Southeast]{Zhizhuo Zhang}
\ead{zhizhuo\_zhang@163.com}
\author[Southeast]{Xiaobing Nie}
\ead{xbnie@seu.edu.cn}
\author[UPVD]{Mikaël Barboteu}
\ead{barboteu@univ-perp.fr}
\author[Southeast]{Jinde Cao\corref{cao}}
\cortext[cao]{corresponding author: Jinde Cao (Email addresses: jdcao@seu.edu.cn)}

\affiliation[Southeast]{organization={School of Mathematics, Southeast University},
            city={Nanjing},
            postcode={211189}, 
            country={China}}
\affiliation[UPVD]{organization={Laboratoire de Mathématiques et Physique, Université de Perpignan Via Domitia},
            addressline={52 Avenue Paul Alduy}, 
            city={Perpignan},
            postcode={66860}, 
            country={France}}

\begin{abstract}
With the increasing demand for the accuracy of numerical simulation of pavement mechanics, the variational inequality model and its induced finite element method which can simulate the interlayer contact state becomes a potential solution. In this paper, a layer decomposition algorithm for solving variational inequality models of multi-layer elastic contact systems with interlayer Tresca friction conditions is studied. Continuous and discrete versions of the algorithm and their convergence theorems have been proposed and proved successively. Then, the algebraic form of the executable optimization algorithm and the numerical experimental results verify the practicability of the variational inequality model and its algorithm in the pavement mechanics modeling.
\end{abstract}


\end{frontmatter}


\input{Main}





\section*{Acknowledgements}
This work was supported by the National Key Research and Development Program of China under Grant No.2020YFA0714300, the National Natural Science Foundations of China under Grant No.61833005 and No.61673111.

\bibliographystyle{plain}
\bibliography{references.bib}


\end{document}


\maketitle

Due to space limitations, the detailed setting of the problem, derivation of variational inequality, explanation and verification of Assumption 4.2, and implementation of the algorithm are placed in this supplementary file.

\section{Setting of the problem}

\subsection{The physical model and PDEs}

In this section, the multilayer elastic system with interlayer Tresca friction will be illustated. The mechanical model of the multilayer elastic system considered in this study is shown in Fig 1. The multilayer elastic system consists of $n$ layers and the open bounded connected region occupied by each layer of deformable bodies made of elastic materials is denoted as $\Omega^{i}\subset\mathbb{R}^{d}$, where $d=3$ and $i=1,2,\ldots,n$. Therefore, the total domain occupied by the system is denoted as $\Omega=\cup_{i=1}^{n}\Omega^{i}$. Then, the boundary $\Gamma^{i}=\partial \Omega^{i}$ is assumed to be Lipschitz continuous and can be decomposed as follows: $\Gamma^{i}=\cup_{j=1}^{3} \bar{\Gamma}_{j}^{i}$, where 
$\Gamma_{1}^{i}$, $\Gamma_{2}^{i}$, $\Gamma_{3}^{i}$ are mutually disjoint relatively open sets. For the convenience of notation, let $\Gamma_{j}=\cup_{i=1}^{n}\bar{\Gamma}_{j}^{i}$, $j=1,2,3$.
Then, $\operatorname{meas}_{d}\left(\cdot\right)$ is used to represent the d-dimensional Lebesgue-measure, and suppose $\operatorname{meas}_{d-1}\left(\Gamma_{1}^{i}\right)>0$ for $i=1,2,\ldots,n$ and $\operatorname{meas}_{d-1}\left(\Gamma_{3}^{n}\right)>0$. 
For $x\in\Gamma_{j}^{i}$, the unit outward normal vector at $x$ is denoted by $\nu$. It is also assumed that the boundary $\Gamma_{2}^{i+1}$ and $\Gamma_{3}^{i}$ coincide when no external force is applied to the system.
However, in the mechanical model, these two boundaries are not considered to be completely equivalent, because they belong to different elastic bodies and their unit outward normal vectors are in opposite directions.
For the convenience of calculation, the unit outward normal vectors of the boundary $\Gamma_{2}^{i}$ and $\Gamma_{3}^{i}$ ($i=1,2,\ldots,n$) are denoted as $\boldsymbol{\alpha}^{i}$ and $\boldsymbol{\beta}^{i}$, respectively, and let the friction boundary $\Gamma_{c}^{i}=\Gamma_{3}^{i}\cap \Gamma_{2}^{i+1}$.

Then, based on the above definition of elastic system, the external forces and boundary conditions in the mechanical model will be given. First, the body force on the elastic body of the $i$-th layer is denoted as $\boldsymbol{f}_{0}^{i}$, $i=1,2,\ldots,n$. Then, the elastic body of the $i$-th layer is fixed on the boundary $\Gamma_{1}^{i}$ and $\Gamma_{3}^{n}$, that is, the zero displacements will be imposed on the boundary $\Gamma_{1}^{i}\cup \Gamma_{3}^{n}$. The surface traction force $\boldsymbol{f}_{2}$ acts on the boundary $\Gamma_{2}^{1}$ of the first layer of elastic body. On the $i$-th layer of elastic body ($i=2,3,\ldots,n$), its boundary $\Gamma_{2}^{i}$ will be subjected to surface forces caused by deformation and friction of the upper layer. Finally, on the boundary $\Gamma_{3}^{i}$ ($i=1,2,\ldots,n-1$), the body is in Tresca frictional contact with the underlying layer of system.

\begin{remark}\label{rem:3.1}
It is worth noting that the Saint-Venant principle ensures that the boundary condition can simulate the real road conditions well, and the quasi-static condition ensures that the small deformation can reach equilibrium in the system in a very short time without considering the effect of acceleration. Therefore, under these two assumptions, the model can well simulate the instantaneous mechanical response of pavement.
\end{remark}

Under the action of force, the elastic body of the $i$-th layer will deform, and its displacement field is denoted as $\boldsymbol{u}^{i}: \Omega^{i} \to \mathbb{R}^{d}$ ($i=1,2,\ldots,n$). At the same time, the stress tensor inside the object is denoted as $\boldsymbol{\sigma}^{i}: \Omega^{i} \to \mathbb{S}^{d}$, where $\mathbb{S}^{d}$ represents the space of second order symmetric tensors on $\mathbb{R}^{d}$. 
Then, "$\cdot$" and "$:$" represent the inner product in $\mathbb{R}^{d}$ and $\mathbb{S}^{d}$, respectively, and "$|\cdot|$" represents the Euclidean norm. And the space of displacement field and stress function is defined as:
$$
\begin{aligned}
V^{i} &=\left\{\boldsymbol{v}=\left(v_{k}\right) \in\left(H^{1}(\Omega^{i})\right)^{d} ~\Big|~ \boldsymbol{v}=\mathbf{0} \text { on } \Gamma^{i}_{1} \right\} , \\
Q^{i} &=\left\{\boldsymbol{\tau}=\left(\tau_{kl}\right) \in\left(L^{2}(\Omega^{i})\right)^{d \times d}~\Big|~ \tau_{lk}=\tau_{lk}, 1 \leqslant k, l \leqslant d\right\}, \\
Q^{i}_{1} &=\left\{\boldsymbol{\tau} \in Q^{i}~\Big|~ \operatorname{Div} \boldsymbol{\tau} \in\left(L^{2}(\Omega^{i})\right)^{d}\right\},
\end{aligned}
$$
where $i=1,2\ldots,n$ and $H^{1}(\cdot)=W^{1,2}(\cdot)$ is Sobolev space. Therefore, on the domain $\Omega$, let $\boldsymbol{u}=\left(\boldsymbol{u}^{1}, \boldsymbol{u}^{2}, \ldots , \boldsymbol{u}^{n} \right)$ and $\boldsymbol{\sigma}=\left(\boldsymbol{\sigma}^{1}, \boldsymbol{\sigma}^{2}, \ldots , \boldsymbol{\sigma}^{n} \right)$. Then $\boldsymbol{u}\in V$ and $\boldsymbol{\sigma}\in Q_{1}$, where
$$
V = V^{1} \times V^{2} \times \cdots \times V^{n} \text{ and } 
Q_{1}= Q_{1}^{1} \times Q_{1}^{2} \times \cdots \times Q_{1}^{n}.
$$
Based on the definition, it is easy to verify that the above spaces are Hilbert spaces, so the standard norm can be defined as follows:
$$
(\boldsymbol{u}^{i},\boldsymbol{v}^{i})_{H^{1}\left(\Omega^{i}\right)^{d}} = \sum_{k=1}^{d} (u^{i}_{k},v^{i}_{k})_{H^{1}\left(\Omega^{i}\right)},~~
(\boldsymbol{\sigma}^{i},\boldsymbol{\tau}^{i})_{L^{2}\left(\Omega^{i}\right)^{d}} = \sum_{k,l=1}^{d} (\sigma^{i}_{kl},\tau^{i}_{kl})_{L^{2}\left(\Omega^{i}\right)},
$$
where $\boldsymbol{u}^{i},\boldsymbol{v}^{i}\in V^{i}$ and $\boldsymbol{\sigma}^{i},\boldsymbol{\tau}^{i}\in Q^{i}$ ($i=1,2,\ldots,n$). From this, the standard inner product on spaces $V$ and $Q$ can be defined as:
$$
(\boldsymbol{u},\boldsymbol{v})_{H^{1}} = \sum_{i=1}^{n} \sum_{k=1}^{d} (u^{i}_{k},v^{i}_{k})_{H^{1}\left(\Omega^{i}\right)},~~
(\boldsymbol{\sigma},\boldsymbol{\tau})_{L^{2}} = \sum_{i=1}^{n} \sum_{k,l=1}^{d} (\sigma^{i}_{kl},\tau^{i}_{kl})_{L^{2}\left(\Omega^{i}\right)}.
$$
Based on the assumption of small deformation, the stress tensor $\boldsymbol{\varepsilon}\left(\boldsymbol{v}^{i}\right)\in Q^{i}$ can be defined as:
$$
\boldsymbol{\varepsilon}\left(\boldsymbol{v}^{i}\right)=\frac{1}{2}\left(\nabla \boldsymbol{v}^{i}+\left(\nabla \boldsymbol{v}^{i}\right)^{T}\right).
$$
Since $\operatorname{meas}_{d-1}\left(\Gamma_{1}^{i}\right)>0$, Korn's inequality holds \cite{3kikuchi1988contact}:
\begin{equation}\label{ieq:3:korn}
\|\boldsymbol{v}^{i}\|_{H^{1}(\Omega^{i})} \leqslant c_{k}^{i}\|\boldsymbol{\varepsilon}(\boldsymbol{v^{i}})\|_{L^{2}\left(\Omega^{i}\right)^{d}} \quad \forall \boldsymbol{v}^{i} \in V^{i},
\end{equation}
where $c_{k}^{i}>0$ is a constant depending only on $\Omega^{i}$ and $\Gamma_{1}^{i}$. Based on Korn's inequality, the inner product can be defined:
$$
(\boldsymbol{u}^{i}, \boldsymbol{v}^{i})_{V^{i}} =(\boldsymbol{\varepsilon}(\boldsymbol{u}^{i}), \boldsymbol{\varepsilon}(\boldsymbol{v}^{i}))_{L^{2}\left(\Omega^{i}\right)^{d}} \quad \forall \boldsymbol{u}^{i}, \boldsymbol{v}^{i} \in V^{i},
$$
and the norm $\|\cdot\|_{V^{i}}$ induced by this inner product is equivalent to the standard norm $\|\cdot\|_{H^{1}(\Omega^{i})^{d}}$ in $V^{i}$ space. Therefore, $\left(V^{i}, \|\cdot\|_{V^{i}}\right)$ is a real Hilbert space.
Similarly, the inner product can be defined in space $V$:
$$
(\boldsymbol{u}, \boldsymbol{v})_{V} = \sum_{i=1}^{n} (\boldsymbol{\varepsilon}(\boldsymbol{u}^{i}), \boldsymbol{\varepsilon}(\boldsymbol{v}^{i}))_{L^{2}\left(\Omega^{i}\right)^{d}} \quad \forall \boldsymbol{u}, \boldsymbol{v} \in V,
$$
It is easy to verify that $\left(V, \|\cdot\|_{V}\right)$ is also a real Hilbert space. Furthermore, according to the trace theorem on Sobolev spaces \cite{3adams2003sobolev}, there exists a constant $c_{t}^{i}$ that depends only on $\Omega^{i}$ and $\Gamma^{i}$, such that the following inequality holds:
\begin{equation}\label{ieq:3:trace}
\|\gamma^{i}_{j}\boldsymbol{v}^{i}\|_{L^{2}\left(\Gamma_{c}^{j}\right)^{d}} \leqslant c_{t}^{i}\|\boldsymbol{v}^{i}\|_{V^{i}},~ \forall \boldsymbol{v}^{i} \in V^{i}, j=i \text{ or }i-1,
\end{equation}
where $\gamma^{i}_{j}: V^{i}\to L^2\left(\Gamma_{c}^{j} \right)$ ($j=i-1$ or $i$) is the usual trace operator. Moreover, $c_{t}^{\min}$ and $c_{t}^{\max}$ can be defined as:
$$
c_{t}^{\min} = \min\{c_{t}^{1},\ldots,c_{t}^{n-1}\} \text{ and } c_{t}^{\max} = \max\{c_{t}^{1},\ldots,c_{t}^{n-1}\}.
$$

In order to facilitate the definition of contact boundary conditions, the displacement and stress on the boundaries $\Gamma^{i}_j$ ($j=2,3$) are decomposed as follows:
\begin{align*}
&v_{\beta}^{i}=\boldsymbol{v}^{i} \cdot \boldsymbol{\beta}^{i}, && \boldsymbol{v}_{\eta}^{i}=\boldsymbol{v}^{i}-v_{\beta}^{i}\cdot\boldsymbol{\beta}^{i}, && \boldsymbol{v}^{i}\in L^{2}\left(\Gamma^{i}_3\right);\\
&v_{\alpha}^{i}=\boldsymbol{v}^{i} \cdot \boldsymbol{\alpha}^{i}, && \boldsymbol{v}_{\tau}^{i}=\boldsymbol{v}^{i}-v_{\alpha}^{i}\cdot\boldsymbol{\alpha}^{i}, && \boldsymbol{v}^{i}\in L^{2}\left(\Gamma^{i}_2\right);\\
&\sigma_{\beta}^{i}=\boldsymbol{\beta}^{i}\cdot\boldsymbol{\sigma}^{i} \cdot \boldsymbol{\beta}^{i}, && \boldsymbol{\sigma}_{\eta}^{i}=\boldsymbol{\sigma}^{i} \cdot \boldsymbol{\beta}^{i}-\sigma_{\beta}^{i}\cdot\boldsymbol{\beta}^{i}, && \boldsymbol{\sigma}^{i}\in Q^{i};\\
&\sigma_{\alpha}^{i}=\boldsymbol{\alpha}^{i}\cdot\boldsymbol{\sigma}^{i} \cdot \boldsymbol{\alpha}^{i}, && \boldsymbol{\sigma}_{\tau}^{i}=\boldsymbol{\sigma}^{i} \cdot \boldsymbol{\alpha}^{i}-\sigma_{\alpha}^{i}\cdot\boldsymbol{\alpha}^{i}, && \boldsymbol{\sigma}^{i}\in Q^{i}.
\end{align*}
On the boundary $\Gamma_{3}^{i}$, $v_{\beta}^{i}$ and $\boldsymbol{v}_{\eta}^{i}$ represent the normal and tangential displacements, respectively, and $\sigma_{\beta}^{i}$ and $\boldsymbol{\sigma}_{\eta}^{i}$ represent the normal and tangential stresses, respectively. Similarly, on the boundary $\Gamma_{2}^{i}$, $v_{\alpha}^{i}$, $\boldsymbol{v}_{\tau}^{i}$, $\sigma_{\alpha}^{i}$ and $\boldsymbol{\sigma}_{\tau}^{i}$ represent the displacements and stresses on the corresponding components.
Moreover, let $v_{N}^{i}=v_{\beta}^{i}$, $\boldsymbol{v}_{T}^{i}=\boldsymbol{v}_{\eta}^{i}$, $\sigma_{N}^{i}=\sigma_{\beta}^{i}$ and $\boldsymbol{\sigma}_{T}^{i}=\boldsymbol{\sigma}_{\eta}^{i}$ on $\Gamma_{c}^{i}$.
In order to characterize the discontinuity of the displacement field on the friction boundary, the jump operator $[\cdot]$ is defined as:
\begin{align*}
\left[{v}_{N}^{i}\right]=\boldsymbol{v}^{i} \cdot \boldsymbol{\beta}^{i}+\boldsymbol{v}^{i+1} \cdot \boldsymbol{\alpha}^{i+1},~~
\left[\boldsymbol{v}_{T}^{i}\right]=\boldsymbol{v}_{\eta}^{i} - \boldsymbol{v}_{\tau}^{i+1},
\end{align*}
where $\boldsymbol{v}^{i}\in L^{2}\left(\Gamma^{i}_{3}\right)$ and $\boldsymbol{v}^{i+1}\in L^{2}\left(\Gamma^{i+1}_{2}\right)$ ($i=1,2,\ldots,n-1$).
Note that the contact between the two layers on the contact boundary $\Gamma_{c}^{i}$ ($i=1,2,\ldots n-1$) satisfies the non-penetration condition, that is $\left[{v}_{N}^{i}\right]\leqslant 0$.

Under the above conditions, the frictional contact problem of the multilayer elastic system can be described by the following partial differential equation:
\begin{problem}[$P_{0}$]\label{prb:3:p_0}
Find a displacement field $\boldsymbol{u}^{i}:\Omega^{i} \rightarrow \mathbb{R}^{d}$ and the stress field $\boldsymbol{\sigma}^{i}: \Omega^{i} \rightarrow \mathbb{S}^{d}$ ($i=1,2,\ldots,n$) such that:
\begin{align}
&\boldsymbol{\sigma}^{i}=\mathcal{A}^{i} \boldsymbol{\varepsilon}(\boldsymbol{u}^{i}) && \text { in } \Omega^{i},\label{eq:3:constitutive}\\
&\operatorname{Div} \boldsymbol{\sigma}^{i}+\boldsymbol{f}_{0}^{i}=\mathbf{0} && \text { in } \Omega^{i},\label{eq:3:balance}\\
&\boldsymbol{u}^{i}=\mathbf{0} && \text { on } \Gamma_{1}^{i}, \label{eq:3:boubdary_displacement}\\
&\boldsymbol{\sigma}^{1} \cdot \boldsymbol{\alpha}^{1}=\boldsymbol{f}_{2} && \text { on } \Gamma_{2}^{1},\label{eq:3:boundary_stress}\\
&\sigma^{i}_{\alpha} = -\sigma^{i-1}_{\beta},~
\boldsymbol{\sigma}_{\tau}^{i} = \boldsymbol{\sigma}_{\eta}^{i-1},
&&\text { on } \Gamma_{2}^{i}, ~i\ne 1, \label{eq:3:boundary_contact_1}\\
&\left.
\begin{aligned}
& \left[{u}_{N}^{i}\right] \leqslant 0,~ \sigma^{i}_{N}\cdot[u_{N}^{i}] = 0,~ |\boldsymbol{\sigma}_{T}^{i}|\leqslant g^{i}(\mathbf{x})\\
& |\boldsymbol{\sigma}_{T}^{i}| < g^{i}(\mathbf{x}) \Rightarrow
[\boldsymbol{u}_{T}^{i}] = 0\\
& |\boldsymbol{\sigma}_{T}^{i}| = g^{i}(\mathbf{x}) \Rightarrow
\boldsymbol{\sigma}_{T}^{i} = -\lambda[\boldsymbol{u}_{T}^{i}], ~\lambda\geqslant 0
\end{aligned}
\right\} &&  \text { on } \Gamma_{c}^{i},~i\ne n.\label{eq:3:boundary_contact_2}
\end{align}
\end{problem}
In the above partial differential equation, formula (\ref{eq:3:constitutive}) represents the constitutive relation of stress and strain, where $\mathcal{A}^{i}$ is a given nonlinear operator, which is called elastic operator. Equation (\ref{eq:3:balance}) represents the equilibrium equation between stress and body force; formulations (\ref{eq:3:boubdary_displacement}) and (\ref{eq:3:boundary_stress}) are displacement field constraints (Dirichlet boundary conditions) and stress boundary conditions, respectively.

Finally, formulations (\ref{eq:3:boundary_contact_1}) and (\ref{eq:3:boundary_contact_2}) represent the frictional contact conditions between the layers of elastic system. Here, Equation (\ref{eq:3:boundary_contact_1}) indicates that the force between the two elastic bodies is equal in magnitude and opposite in direction. Since both sides of the contact surface $\Gamma_{c}^{i}$ are elastic bodies that deform under the action of force, under the Tresca frictional law, the frictional function $g^{i}(x)$ is relative with coordination $x$.
Furthermore, it is easy to prove that 
formulation (\ref{eq:3:boundary_contact_2}) is equivalent to:
\begin{equation}\label{eq:4:boundary_contact_3}
\left\{\begin{aligned}
&\sigma_{N}^{i} \cdot\left[u_{N}^{i}\right]=0,~~\left|\boldsymbol{\sigma}_{T}^{i}\right| \leqslant g^{i}\left(x\right) \\
&g^{i}\left(x\right)\left|\left[\boldsymbol{u}_{T}^{i}\right]\right|+\boldsymbol{\sigma}_{T}^{i} \cdot \left[\boldsymbol{u}_{T}^{i}\right]=0
\end{aligned}\right. .
\end{equation}

\begin{remark}\label{3:rem:1}
It is worth noting that the Tresca frictional contact conditions in (\ref{eq:3:boundary_contact_2}) are usually only applicable to bilateral contact problems, because the friction law requires that the friction region is not changed during deformation. However, in the unilateral contact problem of this model, the friction contact conditions usually satisfy the Coulomb friction law, namely:
$$
\left.
\begin{aligned}
& \left[{u}_{N}^{i}\right] \leqslant 0,~ \sigma^{i}_{N}\cdot[u_{N}^{i}] = 0,~ |\boldsymbol{\sigma}_{T}^{i}|\leqslant g^{i}(\sigma^{i}_{N})\\
& |\boldsymbol{\sigma}_{T}^{i}| < g^{i}(\sigma^{i}_{N}) \Rightarrow
[\boldsymbol{u}_{T}^{i}] = 0\\
& |\boldsymbol{\sigma}_{T}^{i}| = g^{i}(\sigma^{i}_{N}) \Rightarrow
\boldsymbol{\sigma}_{T}^{i} = -\lambda[\boldsymbol{u}_{T}^{i}], ~\lambda\geqslant 0
\end{aligned}
\right\}   \text { on } \Gamma_{c}^{i},~i\ne n-1.
$$
However, because the theoretical analysis and algorithm design of the multilayer elastic system with Coulomb friction law are more difficult, how to approximate it by a simpler model has a strong practical significance. The existing studies show that by designing and solving a fixed point algorithm which computes a unilateral contact problem with Tresca friction law at each step, the unilateral contact problem with Coulomb friction law can be approximately solved \cite{3eck1998existence,3laborde2008fixed}. So, the unilateral contact problem with Tresca friction law has also become a new research direction in the contact problem. In addition, in the study of pavement mechanics, the interlayer contact model lacks theoretical standards, so it is necessary to try multiple contact conditions and carry out comparative studies. Therefore, it is of great significance to study the model of multilayer elastic system with Tresca friction law from both theoretical and application points of view.
\end{remark}

\subsection{Variational inequality}
In order to derive the variational inequality of problem $P_{0}$, the space $\mathcal{K}$ needs to be defined as follows:
$$
\mathcal{K}=\left\{\boldsymbol{v} \in V ~\Big|~ ~\left[v_{N}^{i}\right] \leqslant 0 \text {, a.e. on } \Gamma_{c}^{i},~ i=1,2,\ldots,n-1\right\} \subset V.
$$
Obviously, the the space $\mathcal{K}$ is the closed subspace of $V$.
Then, the elasticity operators $\mathcal{A}^{i}$ in problem $P_{0}$ are assumed to satisfy the following properties:
\begin{equation}\label{character:3:A}
\left\{\begin{aligned}
&\text { (a) } \mathcal{A}^{i}: \Omega^{i} \times \mathbb{S}^{d} \rightarrow \mathbb{S}^{d}; \\ 
&\text { (b) There exists } L_{A}^{i}>0 \text { such that } \\ 
&~~~~~~~~\left|\mathcal{A}^{i}\left(\boldsymbol{x}, \boldsymbol{\varepsilon}_{1}\right)-\mathcal{A}^{i}\left(\boldsymbol{x}, \boldsymbol{\varepsilon}_{2}\right)\right| \leqslant L_{A}^{i}\left|\boldsymbol{\varepsilon}_{1} - \boldsymbol{\varepsilon}_{2}\right|, ~\forall \boldsymbol{\varepsilon}_{1}, \boldsymbol{\varepsilon}_{2} \in \mathbb{S}^{d} \text{ a.e. } \boldsymbol{x} \in \Omega^{i};  \\ 
&\text { (c) There exists } M^{i}>0 \text { such that } \\ 
&~~~~~~~~\left(\mathcal{A}^{i}\left(\boldsymbol{x}, \boldsymbol{\varepsilon}_{1}\right)-\mathcal{A}^{i}\left(\boldsymbol{x}, \boldsymbol{\varepsilon}_{2}\right)\right) : \left(\boldsymbol{\varepsilon}_{1} - \boldsymbol{\varepsilon}_{2}\right) \geqslant M^{i}\left|\boldsymbol{\varepsilon}_{1}-\boldsymbol{\varepsilon}_{2}\right|^{2},  \\
&~~~~~~~~\forall \boldsymbol{\varepsilon}_{1}, \boldsymbol{\varepsilon}_{2} \in \mathbb{S}^{d} \text { a.e. } \boldsymbol{x} \in \Omega^{i}; \\ 
&\text { (d) For any } \varepsilon \in \mathbb{S}^{d}, \boldsymbol{x} \mapsto \mathcal{A}^{i}(\boldsymbol{x}, \boldsymbol{\varepsilon}) \text { is Lebesgue measurable on } \Omega^{i}; \\ 
&\text { (e) The mapping } \boldsymbol{x} \mapsto \mathcal{A}^{i}(\boldsymbol{x}, \mathbf{0}) \in Q^{i}.
\end{aligned}\right.
\end{equation}
And the contact functions $g^{i}$ ($i=1,2,\cdots,n$) satisfies the following constraints:
\begin{equation}\label{character:3:g}
\left\{\begin{aligned}
&\text { (a) } g^{i}: \Gamma_{3}^{i} \rightarrow \mathbb{R}_{+}; \\
&\text { (b) There exists an } L_{j}^{i}>0 \text { such that } \\
&~~~~~~~~\left|g^{i}\left(\boldsymbol{x}\right)-g^{i}\left(\boldsymbol{x}\right)\right| \leqslant L_{j}^{i}\left|u_{1}-u_{2}\right| \\ 
&~~~~~~~~\forall u_{1}, u_{2} \in \mathbb{R} \text { a.e. in } \Omega^{i}; \\
&\text { (c) } \boldsymbol{x} \mapsto g^{i}(\boldsymbol{x}) \in L^{2}\left(\Gamma_{3}^{i}\right) \text { is measurable}.
\end{aligned}\right.
\end{equation}
Note that the constraints on the contact functions here are rather trivial, the strictest of which is that $g^{i}$ is globally Lipschitz continuous. In fact, limited by its own material properties, the stress that an elastic system can withstand is often limited, so this condition can be relaxed to local Lipschitz continuity.

Furthermore, the body force $\boldsymbol{f}_{0}^{i}$ and surface traction force $\boldsymbol{f}_{2}$ of the layered elastic system are assumed to satisfy the following conditions:
\begin{equation}\label{character:3:f}
\boldsymbol{f}_{0}^{i} \in \left(L^{2}(\Omega^{i})\right)^{d},~ 
\boldsymbol{f}_{2} \in \left(L^{2}\left(\Gamma_{2}^{1}\right)\right)^{d}.
\end{equation}

Next, $L^{i}(\boldsymbol{v})$ and $L(\boldsymbol{v})$ is defined as:
\begin{equation}\label{def:3:L}
\begin{aligned}
& L^{i}(\boldsymbol{v}^{i}) = \int_{\Omega^{i}} \boldsymbol{f}_{0}^{i} \cdot \boldsymbol{v}^{i} d x + \int_{\Gamma_{2}^{i}} \boldsymbol{f}_{2} \cdot \boldsymbol{v}^{i} d l, ~ i=1,\ldots,n\\
& L(\boldsymbol{v})=\sum_{i=1}^{n} \int_{\Omega^{i}} \boldsymbol{f}_{0}^{i} \cdot \boldsymbol{v}^{i} d x + \int_{\Gamma_{2}^{1}} \boldsymbol{f}_{2} \cdot \boldsymbol{v}^{1} d l \triangleq (\boldsymbol{f}, \boldsymbol{v})_{V},
\end{aligned}
\end{equation}
$\forall \boldsymbol{v} \in V$. Obviously, $L(\boldsymbol{v})$ is a linear functional about $\boldsymbol{v}$, and because $V$ is a real Hilbert space, the Riesz representation theorem guarantees $\boldsymbol{f} \in V$. Let $j^{i}: \Omega \times V^{i} \times V^{i+1} \rightarrow \mathbb{R}$ and $j: \Omega \times V \rightarrow \mathbb{R}$ be the functional
\begin{equation}\label{def:3:j}
\begin{aligned}
& j^{i}(\mathbf{x}, \boldsymbol{w}^{i}, \boldsymbol{w}^{i+1})= \int_{\Gamma_{c}^{i}} g^{i}\left(\mathbf{x}\right)|\boldsymbol{w}_\eta^i-\boldsymbol{w}_\tau^{i+1}| d l,~ i=1,\ldots,n-1;\\
& j(\mathbf{x}, \boldsymbol{w})= \sum_{i=1}^{n-1} j^{i}(\mathbf{x}, \boldsymbol{w}^{i}, \boldsymbol{w}^{i+1})= \sum_{i=1}^{n-1} \int_{\Gamma_{c}^{i}} g^{i}\left(\mathbf{x}\right)|[\boldsymbol{w}_{T}^{i}]| d l,~~~ \forall \boldsymbol{w} \in V.
\end{aligned}
\end{equation}
Furthermore, the double-variable operator $a^{i}(\cdot,\cdot): V^{i}\times V^{i} \to \mathbb{R}$ ($i=1,2,\ldots,n$) and $a(\cdot,\cdot): V\times V \to \mathbb{R}$ are defined as:
\begin{align}
a^{i}(\boldsymbol{v}^{i}, \boldsymbol{w}^{i})&=\int_{\Omega^{i}} \mathcal{A}^{i} \boldsymbol{\varepsilon}(\boldsymbol{v}^{i}): \boldsymbol{\varepsilon}(\boldsymbol{w}^{i}) \mathrm{d} x 
= \left(\mathcal{A}^{i} \boldsymbol{\varepsilon}(\boldsymbol{v}^{i}), \boldsymbol{\varepsilon}(\boldsymbol{w}^{i})\right)_{L^{2}\left(\Omega^{i}\right)}
:= \left(A^{i}\boldsymbol{v}^{i}, \boldsymbol{w}^{i}\right)_{V}\label{def:3:ai}\\
a(\boldsymbol{v}, \boldsymbol{w})&=\sum_{i=1}^{n} a^{i}(\boldsymbol{v}^{i}, \boldsymbol{w}^{i}):= \left(A\boldsymbol{v}, \boldsymbol{w}\right)_{V}\label{def:3:a}
\end{align}
$\forall \boldsymbol{v}, \boldsymbol{w} \in V$. Note that $A^{i}:V^{i}\to V^{i}$ is a operator defined by $\mathcal{A}^{i}$. Thus operator $a(\cdot,\cdot)$ is nonlinear with respect to the first variable and linear with respect to the second variable.

Based on Green's formula in tensor form, the solution to the displacement field of problem $P_{0}$ can be transformed into a solution to a variational inequality problem of the form:
\begin{problem}[$P_{0}^{v}$]\label{prb:3.2:P_0_v}
Find a displacement $\boldsymbol{u}\in \mathcal{K}$ which satisfies:
\begin{equation}\label{ieq:prb:3.2:var}
a(\boldsymbol{u},\boldsymbol{v}-\boldsymbol{u}) + j(\mathbf{x},\boldsymbol{v}) - j(\mathbf{x},\boldsymbol{u}) \geqslant L(\boldsymbol{v}-\boldsymbol{u}), ~\forall \boldsymbol{v} \in \mathcal{K}.
\end{equation}
\end{problem}

\begin{proof}
Multiplying both sides of the equilibrium equation (\ref{eq:3:balance}) for the $i$-layer elastic body by $\boldsymbol{v}^{i}-\boldsymbol{u}^{i}$ and integrating in the region $\Omega^{i}$ yields:
\begin{equation}\label{proof:3.1.1}
\begin{aligned}
&\int_{\Gamma_{2}^{i}} \boldsymbol{\sigma}^{i} \cdot \alpha^{i} \cdot \left(\boldsymbol{v}^{i} - \boldsymbol{u}^{i}\right) d a + \int_{\Gamma_{3}^{i}} \boldsymbol{\sigma}^{i} \cdot \beta^{i} \cdot \left(\boldsymbol{v}^{i} - \boldsymbol{u}^{i}\right) d a \\
& + \int_{\Omega^{i}} \boldsymbol{f}^{i}_{0}\left(\boldsymbol{v}^{i} - \boldsymbol{u}^{i}\right) d x 
= \int_{\Omega^{i}} \boldsymbol{\sigma}^{i} :  \left(\boldsymbol{\varepsilon}(\boldsymbol{v}^{i}) - \boldsymbol{\varepsilon}(\boldsymbol{u}^{i})\right) d x
\end{aligned}
\end{equation}
In this process, Green's formula in tensor form holds:
$$
(\boldsymbol{\sigma}^{i}, \boldsymbol{\varepsilon}(\boldsymbol{v}^{i}))_{L^{2}\left(\Omega^{i}\right)^{d}} + (\operatorname{Div} \boldsymbol{\sigma}, \boldsymbol{v}^{i})_{L^{2}\left(\Omega^{i}\right)} = \int_{\Gamma^{i}} \boldsymbol{\sigma}^{i} \cdot \nu \cdot \boldsymbol{v}^{i} \mathrm{d} a ~~~ \forall \boldsymbol{v} \in V^{i},
$$
where $\nu$ represents unit outer normal on $\Gamma^{i}$.
Then, by accumulating the formula (\ref{proof:3.1.1}), it can be obtained:
\begin{align}\label{proof:3.1.2}
a(\boldsymbol{u},\boldsymbol{v}-\boldsymbol{u})
=& \sum_{i=1}^{n} \int_{\Gamma_{2}^{i}} \boldsymbol{\sigma}^{i} \cdot \alpha^{i} \cdot \left(\boldsymbol{v}^{i} - \boldsymbol{u}^{i}\right) d l + \sum_{i=1}^{n}\int_{\Gamma_{3}^{i}} \boldsymbol{\sigma}^{i} \cdot \beta^{i} \cdot \left(\boldsymbol{v}^{i} - \boldsymbol{u}^{i}\right) d l + \sum_{i=1}^{n}\int_{\Omega^{i}} \boldsymbol{f}^{i}_{0}\left(\boldsymbol{v}^{i} - \boldsymbol{u}^{i}\right) d x \nonumber \\
=& \sum_{i=1}^{n-1}\left( \int_{\Gamma_{2}^{i+1}} \boldsymbol{\sigma}^{i+1} \cdot \alpha^{i+1} \cdot \left(\boldsymbol{v}^{i+1} - \boldsymbol{u}^{i+1}\right) d l + \int_{\Gamma_{3}^{i}} \boldsymbol{\sigma}^{i} \cdot \beta^{i} \cdot \left(\boldsymbol{v}^{i} - \boldsymbol{u}^{i}\right) d l \right) + L(\boldsymbol{v}-\boldsymbol{u}).
\end{align}
Based on the decomposition of force and displacement in normal and tangential directions and boundary conditions (\ref{eq:3:boundary_contact_1}) and (\ref{eq:3:boundary_contact_2}), the first term of equation (\ref{proof:3.1.2}) satisfies:
\begin{equation}\label{proof:3.1.3}
\begin{aligned}
&\sum_{i=1}^{n-1}\left( \int_{\Gamma_{2}^{i+1}} \boldsymbol{\sigma}^{i+1} \cdot \alpha^{i+1} \cdot \left(\boldsymbol{v}^{i+1} - \boldsymbol{u}^{i+1}\right) d l + \int_{\Gamma_{3}^{i}} \boldsymbol{\sigma}^{i} \cdot \beta^{i} \cdot \left(\boldsymbol{v}^{i} - \boldsymbol{u}^{i}\right) d l \right)\\
=&\sum_{i=1}^{n-1}\Bigg(\int_{\Gamma_{2}^{i+1}} {\sigma}^{i+1}_{\alpha} \cdot \left({v}^{i+1}_{\alpha} - {u}^{i+1}_{\alpha}\right) d l + \int_{\Gamma_{2}^{i+1}} \boldsymbol{\sigma}^{i+1}_{\tau} \cdot \left(\boldsymbol{v}^{i+1}_{\tau} - \boldsymbol{u}^{i+1}_{\tau}\right) d l \\
&+ \int_{\Gamma_{3}^{i}} {\sigma}^{i}_{\beta} \cdot \left({v}^{i}_{\beta} - {u}^{i}_{\beta}\right) d l + \int_{\Gamma_{3}^{i}} \boldsymbol{\sigma}^{i}_{\eta} \cdot \left(\boldsymbol{v}^{i}_{\eta} - \boldsymbol{u}^{i}_{\eta}\right) d l \Bigg)\\
=& \sum_{i=1}^{n-1} \int_{\Gamma_{c}^{i}} \boldsymbol{\sigma}^{i}_{T} \cdot \left(\left[\boldsymbol{v}_{T}^{i}\right] - \left[\boldsymbol{u}_{T}^{i}\right]\right) d l
\end{aligned}
\end{equation}
According to the mutual non-penetration condition, ${\sigma}^{i}_{N}\cdot [u_{N}^{i}] = 0$ and when $\boldsymbol{v}\in\mathcal{K}$, ${\sigma}^{i}_{N}\cdot [v_{N}^{i}] \geqslant 0$.
Therefore, by substituting formulations (\ref{proof:3.1.3}) into equation (\ref{proof:3.1.2}) and combining with boundary conditions (\ref{eq:3:boundary_contact_1}) and (\ref{eq:3:boundary_contact_2}), the following inequality relation can be derived:
\begin{align*}
a(\boldsymbol{u},\boldsymbol{v}-\boldsymbol{u})
=&  L(\boldsymbol{v}-\boldsymbol{u}) + \sum_{i=1}^{n-1} \int_{\Gamma_{c}^{i}} \boldsymbol{\sigma}^{i}_{T} \cdot \left(\left[\boldsymbol{v}_{T}^{i}\right] - \left[\boldsymbol{u}_{T}^{i}\right]\right) d a\\
\geqslant & L(\boldsymbol{v}-\boldsymbol{u}) - \sum_{i=1}^{n-1} \int_{\Gamma_{c}^{i}} g^{i}(\mathbf{x})) \cdot \left(|[\boldsymbol{v}_{T}^{i}]| - |[\boldsymbol{u}_{T}^{i}]|\right) d a\\
\geqslant & L(\boldsymbol{v}-\boldsymbol{u}) - j(\mathbf{x},\boldsymbol{v}) + j(\mathbf{x},\boldsymbol{u}).
\end{align*}

Hence, problem $P_{0}^{v}$ can be derived from problem $P_{0}$, that is, the solution of problem $P_{0}$ must be the solution of problem $P_{0}^{v}$. If the solution of problem $P_{0}^{v}$ exists and is unique, then problem $P_{0}^{v}$ is equivalent to problem $P_{0}$.
\end{proof}

Under the framework of variational inequalities, the constraints of the original problem $P_{0}$ are greatly simplified. Our previous research has proved that when the above elastic operators $\mathcal{A}^i$ and contact functions $g^i$ satisfy the regularity conditions (\ref{character:3:A}), (\ref{character:3:g}), the solution of the variational inequality (\ref{ieq:prb:3.2:var}) of the multilayer elastic system is unique, and the numerical solution obtained by the finite element method also satisfies the convergence theorem. However, the above conclusions remain in the theoretical level, and the algorithm for solving the finite element numerical solution of the multilayer elastic system has not been proposed. 
Therefore, based on the domain decomposition algorithm\cite{3haslinger2014domain}, a layer decomposition algorithm suitable for this problem is proposed.
Thereafter, the elastic operator $\mathcal{A}^i$ is restricted to a linearly symmetric fourth-order tensor, and thus $a(\cdot,\cdot)$ is commutative bilinear form.

\section{Assumption 4.2}

\begin{remark}\label{rem:3.2}
In order to demonstrate the rationality of Assumption 4.2, the following example is given. The earliest example can be found in Ref \cite{3bjorstad1986iterative}.
Suppose that problems (2.5) are regular as follows:
\begin{enumerate}[label = (\roman*)]
\item $\exists \epsilon \in(0,1 / 2)$ such that for any $\lambda^{i} \in H^{1 / 2+\epsilon}\left(\Gamma_c^{i}\right)$ the solution $\chi^i_{j} \lambda^{i}$ of (2.5) belongs to $\left(H^{1+\epsilon}\left(\Omega^j\right)\right)^d \cap {V}^j$ and
$$
\exists c:=c(\epsilon)=\text{const.}>0:\left\|\chi^i_j \lambda^{i} \right\|_{1+\epsilon, \Omega^j} \leq c\|\lambda^{i}\|_{1 / 2+\epsilon, \Gamma_c^{i,i+1}},~ \forall \lambda^{i} \in H^{1 / 2+\epsilon}\left(\Gamma_c^{i}\right)
$$
(for the definition of the fractional Sobolev spaces and their norms we refer to [Sobolev spaces]\cite{3adams2003sobolev}).

Further we suppose that the following inverse inequality is satisfied for elements of ${W}_h$
\item there exists a constant $c>0$ independent of $h>0$ such that
$$
\left\|\lambda_h^{i}\right\|_{1 / 2+\epsilon, \Gamma_c^{i,i+1}} \leq c h^{-\epsilon}\left\|\lambda_h^{i} \right\|_{1 / 2, \Gamma_c^{i,i+1}} \quad \forall \lambda_h^{i} \in {W}_h^{i} .
$$
\end{enumerate}
Then under (i) and (ii), the assumption (4.5) easily follows. Indeed, let $\lambda_h^{i} \in {W}_h^{i}$ be arbitrary. Then
$$
\begin{aligned}
\| \lambda_h^{i} \|_{\Gamma_{c h}^{i,i+1}} & \stackrel{(4.3)}{=} \|\hat{\chi}_{i+1}^i \lambda_h^{i} \|_{1, \Omega^{i+1}} \leq \|\chi^{i}_{i+1} \lambda_{h}^{i} \|_{1, \Omega^{i+1}} +\|\chi_{i+1}^i \lambda_h^{i} - \hat{\chi}_{i+1}^i \lambda_h^i \|_{1, \Omega^{i+1}} \\
& \stackrel{(2.6)}{=} \|\lambda_{h}^{i}\|_{\Gamma_{c}^{i,i+1}} + c h^{ \epsilon } \|\chi^i_{i+1}\lambda_{h}^{i} \|_{1+\epsilon, \Omega^i} \stackrel{(i)+(i i)}{\leq} c \| \lambda_h^{i} \|_{\Gamma_c^{i,i+1}}
\end{aligned}
$$
using that $\hat{\chi}_{i+1}^i \lambda_h^{i}$ is the Galerkin approximation of $\chi_{i+1}^i \lambda_h^{i}$.
\end{remark}

\section{Implementable Algorithm}

In main text, the discrete form of the layer decomposition algorithm have been deduced. However, this algorithm can be used to calculate the mechanical response of multilayer elastic systems only after it is be transformed into an algebraic optimization model. Therefore, in this section, there are two topics need to be discussed, that is: discrete optimization model, algebraic optimization model.  

\subsection{Discrete optimization model}

In previous section, the the layer decomposition algorithm have been presented in discrete variational form. However, the algorithm in this form is convenient to be discussed in theory and difficult to be used in practical calculation. According to the theory of variational inequalities \cite{3han2002quasistatic}, the elliptic variational inequality of the second kind, which can be expressed as:
$$
\begin{aligned}
&\text{Find } u \in V \text{ such that:}\\
&~~~~~~ a(u,v-u) + j(v) - j(u) \geqslant L(v-u),~ \forall v \in V,
\end{aligned}
$$
is equivalent to the following optimization problem:
$$
\min _{v \in V} J(v) = \frac{1}{2} a(v, v) +j(v) -L(v).
$$
Therefore, it is not difficult to verify that the Problem Problem $P_{d h}^i\left(\boldsymbol{\lambda}_h\right)$ can be equivalently deformed to 

\begin{problem}[ $P_{oh}^i(\boldsymbol{\lambda}_{h})$]\label{prb:3.9:oh}
\begin{equation}\label{prb:3.9:oh1}
\left\{
\begin{aligned}
&\text{Find } \boldsymbol{u}^1_{h}:=\boldsymbol{u}^1_{h}(\boldsymbol{\lambda}_{h}) \in K^{1}_{h}(\lambda^{1}_{h} \cdot \boldsymbol{\alpha}^{2}) \text{ such that:}\\
& \boldsymbol{u}^1_{h} = \arg\min_{\boldsymbol{v}^1_{h} \in K^{1}_{h}(\lambda^{1}_{h} \cdot\boldsymbol{\alpha}^{2})} \frac{1}{2} a^1\left(\boldsymbol{v}^1_{h}, \boldsymbol{v}^1_{h} \right) + j^1\left(\mathbf{x}, \boldsymbol{v}^1_{h}, \lambda^1_{h} \right) - L^1\left(\boldsymbol{v}^1_{h} \right) ,
\end{aligned}
\right.  \tag{\text{$P^{1}_{oh}(\boldsymbol{\lambda})$}}
\end{equation}
\begin{small}
\begin{equation}\label{prb:3.9:ohi}
\left\{
\begin{aligned}
&\text{Find } \boldsymbol{u}^i_{h} :=\boldsymbol{u}^i_{h}(\boldsymbol{\lambda}_{h}) \in K^{i}_{h} (\lambda^{i}_{h} \cdot \boldsymbol{\alpha}^{i+1}) \cap V^{i}_{2h}(\lambda^{i-1}_{h} ), ~i=2,\ldots,n-1 \text{ such that:}\\
& \boldsymbol{u}^i_{h} = \arg\min_{\boldsymbol{v}^i_{h} \in K^{i}_{h}(\lambda^{i}_{h} \cdot \boldsymbol{\alpha}^{i+1}) \cap V^{i}_{2h}(\lambda^{i-1}_{h})} \frac{1}{2} a^i\left(\boldsymbol{v}^i_{h}, \boldsymbol{v}^i_{h} \right) + j^i\left(\mathbf{x}, \boldsymbol{v}^i_{h}, \lambda^i_{h} \right) - L^i\left(\boldsymbol{v}^i_{h} \right),
\end{aligned}
\right.  \tag{\text{$P^{i}_{oh}(\boldsymbol{\lambda})$}}
\end{equation}
\end{small}
\begin{equation}\label{prb:3.9:ohn}
\left\{
\begin{aligned}
&\text{Find } \boldsymbol{u}^n_{h} := \boldsymbol{u}^n_{h}(\boldsymbol{\lambda}_{h}) \in V_{2 h}^n\left(\lambda_h^{n-1}\right) \text{ such that:}\\
&\boldsymbol{u}^n_{h} = \arg\min_{\boldsymbol{v}^n_{h} \in  V^{n}_{2h}(\lambda^{n-1}_{h})} \frac{1}{2}a^n\left(\boldsymbol{v}^n_{h}, \boldsymbol{v}^n_{h} \right) - L^n\left(\boldsymbol{v}^n_{h} \right).
\end{aligned}
\right.  \tag{\text{$P^{n}_{oh}(\boldsymbol{\lambda})$}}
\end{equation}
\end{problem}

Hence, the layer decomposition algorithm was designed by solving the optimization problems and the equation problems.

\subsection{Algebra optimization model}

The definition of finite element space $V_h^i$ ($i=1,\ldots,n$) is indispensable prerequisite for translating the above optimization algorithm into corresponding algebraic optimization algorithm. Since the finite element spaces $V_h^i$ and bases in spaces have been defined, some basic parameters need to be defined. 
First, $n^{i}=\dim V_{h}^{i}=d\cdot N^{i}_{b}$ ($i=1,\ldots,n$) represents the dimension of finite element space $V_{h}^{i}$, where $N_{b}^{i}$ is used to represent the number of finite element nodes in $\Omega^{i}$, $m^{i}$ ($i=1,\ldots,n-1$) stands for the number of contact nodes of ${\Gamma}_{3 h}^j$ and $d=3$ is dimension of displacement space. Then, in finite element spaces, the operators $a^{i}(\cdot,\cdot)$, $L^{i}(\cdot)$ and $j(\cdot,\cdot)$ can be expressed  in the following algebraic forms:
\begin{small}
\begin{align}
& a^i\left(\boldsymbol{u}^i_{h}, \boldsymbol{v}^i_{h}\right) = \int_{\Omega^i} \mathcal{A}^i \boldsymbol{\varepsilon}\left(\boldsymbol{u}^i_{h}\right): \boldsymbol{\varepsilon}\left(\boldsymbol{v}^i_{h}\right) \mathrm{d} x = (\boldsymbol{v}^i_{v})^{\top} A_{h}^{i} \boldsymbol{u}^i_{v}, \label{3:eq:a.m}\\
& L^{i}(\boldsymbol{v}^i_{h}) = \int_{\Omega^i} f_0^i \cdot \boldsymbol{v}^i_{h} d x + \int_{\Gamma_2^i} f_2 \cdot \boldsymbol{v}^i_{h} d l= (\boldsymbol{v}^i_{v})^{\top} \boldsymbol{f}_{0h}^i + (\boldsymbol{v}^i_{v})^{\top} \boldsymbol{f}_{2h}, \label{3:eq:L.m}\\
&  j^i\left(\mathbf{x}, \boldsymbol{w}^i_{h}, \lambda^{i}_{h}\right)=\int_{\Gamma_c^i} g^i(\mathbf{x})\left|\boldsymbol{w}_{\eta h}^i-\lambda_{\tau h}^{i}\right| d l = \int_{\Gamma_c^i} g^i(\mathbf{x}) \left| N_{2}^{i} \cdot \Phi^{i} \cdot \left(T^{i}_{1} \boldsymbol{w}_{v}^i -T^{i}_{1} \lambda_{v}^{i} \right) \right| d l,\label{3:eq:j.m.1}
\end{align}
\end{small}
where $\boldsymbol{u}^i_{v}$, $\boldsymbol{v}^i_{v}$ and $\boldsymbol{w}^i_{v} \in \mathbb{R}^{n^{i}}$ are used to represent the displacement vector for each finite element node on region $\Omega^{i}$, and the displacement vector of the contact nodes on the boundary $\Gamma_{3h}^{i}$ is denoted by $\lambda^{i}_{v}\in\mathbb{R}^{n^{i}}$ . Furthermore,  $A_h^i\in \mathbb{R}^{n^{i}\times n^{i}}$ and $\boldsymbol{f}_{0h}^i$ ($\boldsymbol{f}_{2h}^i$) are the symmetric, positive definite stiffness matrix and load vector, respectively.  Before defining the matrix $T_{v}^{i} \in \mathbb{R}^{n^{i}\times n^{i}}$, let $N^{i}_{1}\in \mathbb{R}^{m^{i}\times n^{i}}$ and $I_{j}^{i}\in\mathbb{R}^{n^{i}\times n^{i}}$ ($j=1,2$), where each row of $N^{i}_{1}$ contains at most $d$ non-zero entries, which corresponds the unit outward normal vector of finite element node on ${\Gamma}_{3 h}^i$. And the matrices $I_{1}^{i}$ and $I_{2}^{i}$ are symmetric, where non-zero elements (equal only to $1$) corresponding to contact nodes on $\Gamma_{3}^{i}$ and $\Gamma_{2}^{i}$, respectively, exist only above the diagonal. Therefore, $N^{i}_{1}\cdot \boldsymbol{u}^i_{v} \in \mathbb{R}^{m^{i}}$ is used to represent the displacement length of every contact node on ${\Gamma}_{3 h}^i$ along direction of the unit outward normal vector and $I^{i}_{1}\cdot \boldsymbol{u}^i_{v}, I^{i}_{2}\cdot \boldsymbol{u}^i_{v} \in \mathbb{R}^{n^{i}}$ are used to obtain the displacement vector of every contact node on ${\Gamma}_{3 h}^i$ and ${\Gamma}_{2 h}^i$, respectively. So, it is not hard to verify that $I^{i}_{1}\cdot \boldsymbol{u}^i_{v} - (N^{i}_{1})^{\top}\cdot N^{i}_{1} \cdot \boldsymbol{u}^i_{v} \in \mathbb{R}^{n^{i}}$ is the tangent displacement vector of contact node along the contact surface $\Gamma_{c}^{i}$. Then, let $T_{1}^{i} = I^{i}_{1} - (N^{i}_{1})^{\top}\cdot N^{i}_{1}$. It should be noticed that matrices $N^{i}_{1}$ and $T_{1}^{i}$ enable us to write the algebraic form of the non-penetration and the friction condition, respectively. Then, $\Phi^{i}=diag(\phi^{i}_1, \ldots, \phi^{i}_{n^{i}})\in\mathbb{R}^{n^{i}\times n^{i}}$ is basis function matrix, that is, each diagonal element of $\Phi^{i}$ corresponds to a basis function of the finite element space $V_{h}^{i}$, and the function of matrix $N_2^i\in\mathbb{R}^{3\times n^{i}}$ is to add the coaxially based functions. For convenience of calculation, let $I_{3}^{i}, I_{4}^{i}\in \mathbb{R}^{3m^{i}\times n^{i}}$ be the matrices consisting of all the non-zero rows of the matrix $I_{1}^{i}$ and $I_{2}^{i}$, respectively. If we let $|\boldsymbol{w}^i_{v}|=\sum_{j=1}^{n^{i}}|{w}^i_{v,j}|$, $|\boldsymbol{w}^i_{v}|^{v} = \left(\left|w_{v,1}^i\right|,\left|w_{v,2}^i\right|, \ldots,\left|w_{v,n^i}^i\right|\right)^{\top}$ and $|\boldsymbol{w}^i_{v}|^{v}_{j} = |{w}^i_{v,j}|$ for $\boldsymbol{w}^i\in \mathbb{R}^{n^{i}}$, then under the assumption that the displacement vectors of adjacent contact nodes have the same sign, the formulation (\ref{3:eq:j.m.1}) can be rewritten as 
\begin{equation}\label{3:eq:j.m.2}
j^i\left(\mathrm{x}, \boldsymbol{w}_h^i, \lambda_h^i\right)=\int_{\Gamma_c^i} g^i(\mathrm{x})\left|N_2^i \cdot \Phi^i \cdot\left(T_1^i \boldsymbol{w}_v^i-T_1^i \lambda_v^i\right)\right| d l = \boldsymbol{g}^{i}\cdot \left|T_1^i \boldsymbol{w}_v^i-T_1^i \lambda_v^i\right|^{v},
\end{equation}
where $\boldsymbol{g}^{i} = \left(g^{i}_{1},g^{i}_{2},\ldots,g^{i}_{n^{i}}\right)^{\top}$ and $g^{i}_{j}= \int_{\Gamma_{c}^{i}}g^i(\mathrm{x}) \phi^{i}_{j} dl$.

Finally, the finite element nodes on boundary $\Gamma_{1}^{i}$ need be dealt with. The number of boundary nodes on Dirichlet boundary $\Gamma_{1}^{i}$ is denoted by $m_{0}^{i}$, then the matrix $T_{0}^{i}\in \mathbb{R}^{3m^{i}_{0}\times n^{i}}$ is used to take out the coordinates of boundary nodes. Therefore, the displacement boundary condition in algorithm can be illustrated as $T_{0}^{i}\cdot \boldsymbol{v}_{v}^{i} = \boldsymbol{0}$. Of course, based on this fact, the finite element space $V_{h}^{i}$ has been pre-processed, that is, rows and columns of every matrix or vector corresponding to the finite element nodes on zero displacement boundary $\Gamma_{1}^{i}$ have been deleted, which does note affect the symmetry of matrix. Therefore, all matrix and vectors used in algorithm have been pre-processed by default. 

The algebraic counterpart of Discrete Layer Decomposition Algorithm 4.1  read as follow:

\begin{algorithm}[!ht]
\caption{Algebraic Form of Discrete Layer Decomposition Algorithm}
\label{algorithm:3.3}
Let $\lambda_{0,v}^{i}\in \mathbb{R}^{n^{i}}$, $i=1,\ldots,n-1$ and $\theta>0$ be given.

\textbf{Output:} $\boldsymbol{u}^{i}_{k,v}$, $\boldsymbol{p}^{i}_{k,v}$, $\boldsymbol{q}^{i}_{kv,}$, $\lambda_{k,v}^{i}$

\textbf{For} $k\leq N$ \textbf{do}

~~~~$
\left\{
\begin{aligned}
& \boldsymbol{u}^{1}_{k,v} := \arg\min \frac{1}{2} \left(\boldsymbol{v}_v^1\right)^{\top} A_h^1 \boldsymbol{v}_v^1 + \boldsymbol{g}^1 \cdot\left|T_1^1 \boldsymbol{v}_v^1-T_1^1 \lambda_{k-1,v}^1\right|^v - \left(\left(\boldsymbol{v}_v^1\right)^{\top} 
\boldsymbol{f}_{0h}^1 + \left(\boldsymbol{v}_v^1\right)^{\top} \boldsymbol{f}_{2h} \right) \\ 
& ~~~~~ \text{ s.t. } N_{1}^{1} \cdot\left( \boldsymbol{v}_v^1 - \lambda_{k-1,v}^{1} \right) \leqslant 0, \\
& \boldsymbol{u}^{i}_{k,v} := \arg\min \frac{1}{2} \left(\boldsymbol{v}_v^i\right)^{\top} A_h^i \boldsymbol{v}_v^i + \boldsymbol{g}^i \cdot\left|T_1^i \boldsymbol{v}_v^i-T_1^i \lambda_{k-1,v}^i\right|^v - \left( \boldsymbol{v}_v^i \right)^{\top} 
\boldsymbol{f}_{0h}^i  \\ 
& ~~~~~ \text{ s.t. } N_{1}^{i} \cdot\left( \boldsymbol{v}_v^i - \lambda_{k-1,v}^{i} \right) \leqslant 0 \text{ and } I_{4}^{i}\cdot \boldsymbol{v}_v^i = I_{3}^{i-1}\lambda_{k-1,v}^{i-1},\\
& \boldsymbol{u}^{n}_{k,v} := \arg\min \frac{1}{2} \left(\boldsymbol{v}_v^n\right)^{\top} A_h^n \boldsymbol{v}_v^n - \left( \boldsymbol{v}_v^n \right)^{\top}  \boldsymbol{f}_{0h}^n  \\ 
& ~~~~~ \text{ s.t. } I_{4}^{n}\cdot \boldsymbol{v}_v^n = I_{3}^{n-1}\lambda_{k-1,v}^{n-1}
\end{aligned}
\right.
$

~~~~$
\left\{
\begin{aligned}
& A_{h}^{i+1}\boldsymbol{p}_{k,v}^{i+1}= \frac{1}{2} (I_{4}^{i+1})^{\top} \left( I^{i+1}_{4}\left( A_{h}^{i+1}\boldsymbol{u}_{k,v}^{i+1} -\boldsymbol{f}_{0h}^{i+1} \right) + I^{i}_{3}\left( A_{h}^{i} \boldsymbol{u}_{k,v}^{i} -\boldsymbol{f}_{0h}^{i} \right) \right)\\
& A_{h}^{i}\boldsymbol{q}_{k,v}^{i}= \frac{1}{2} (I_{3}^{i})^{\top} \left( I^{i+1}_{4}\left( A_{h}^{i+1}\boldsymbol{u}_{k,v}^{i+1} -\boldsymbol{f}_{0h}^{i+1} \right) + I^{i}_{3}\left( A_{h}^{i} \boldsymbol{u}_{k,v}^{i} -\boldsymbol{f}_{0h}^{i} \right) \right)
\end{aligned}
\right.
$

~~~~$I_{3}^{i}\lambda_{k,v}^i=I_{3}^{i}\lambda_{k-1,h}^{i} -\theta\left( I_{4}^{i+1}\boldsymbol{p}^{i+1}_{kh} + I_{3}^{i}\boldsymbol{q}^{i}_{kh}\right)$.

\textbf{end}
\end{algorithm}

Since there are non-differentiable terms in above optimization problems, it is challenging to solve the algorithm directly. But, it can be noticed that the dual formulation of the above optimization problems fairly simplify the implementation of this Discrete  Layer Decomposition Algorithm by removing the non-differentiable terms and reducing the inequality constrains. 

First, the Lagrange functions for all optimization problems in Step $2$ of the Algorithm \ref{algorithm:3.3} can be expressed as follow:
\begin{equation}\label{3:eq:Lagrange}
\begin{aligned}
\mathcal{L}^{i}\left(\boldsymbol{v}_{v}^{i}, \boldsymbol{\omega}^{i}\right) = &\frac{1}{2}\left(\boldsymbol{v}_v^i\right)^{\top} A_h^i \boldsymbol{v}_v^i - \left(\boldsymbol{v}_v^i\right)^{\top} \boldsymbol{b}^i + \boldsymbol{g}^i \cdot\left|T_1^i \boldsymbol{v}_v^i-T_1^i \lambda_{k-1, v}^i\right|^v \\ 
&+ \left(\boldsymbol{\omega}^{i}_{1N}\right)^{\top} N_1^i \cdot\left(\boldsymbol{v}_v^i-\lambda_{k-1, v}^i\right) + \left(\boldsymbol{ \omega}^{i}_{2} \right)^{\top} \left(I_4^n \cdot \boldsymbol{v}_v^n - I_3^{n-1} \lambda_{k-1, v}^{n-1}\right)
\end{aligned}
\end{equation}
where $\boldsymbol{\omega}_{1 N}^i\geqslant\boldsymbol{0}\in\mathbb{R}^{m^{i}}$, $\left|\boldsymbol{\omega}_{1 T}^i\right|^{v}\leqslant I_{3}^{i} \boldsymbol{g}^{i} \in\mathbb{R}^{3m^{i}}$, $\boldsymbol{\omega}_{2}^i\in\mathbb{R}^{3m^{i-1}}$ and $\boldsymbol{b}^{i}=\boldsymbol{f}^{i}_{0h}$ or $\boldsymbol{f}_{0h}^{1}+\boldsymbol{f}_{2h}$.  
Then, the domain $\Upsilon^{i}$ is defined as:
$$
\Upsilon^{i} := \left\{\boldsymbol{\omega}^{i} = \left( \begin{array}{c}
\boldsymbol{\omega}^{i}_{1N}\\
\boldsymbol{\omega}^{i}_{1T}\\
\boldsymbol{\omega}^{i}_{2}
\end{array} \right) \in \mathbb{R}^{4m^{i}+3m^{i-1}} :\boldsymbol{\omega}_{1 N}^i\geqslant\boldsymbol{0}, \left|\boldsymbol{\omega}_{1 T}^i\right|^{v}\leqslant I_{3}^{i} \boldsymbol{g}^{i}, \boldsymbol{\omega}^{1}_{2}= \boldsymbol{0}, \boldsymbol{\omega}^{n}_{1N}= \boldsymbol{0}, \boldsymbol{\omega}^{n}_{1T}= \boldsymbol{0} \right\}
$$
Therefore, this optimization problems can be equivalently transformed into the following saddle-point problems:
\begin{equation}\label{3:prb:saddle_problem}
\left.\begin{array}{l}
\text { Find }\left(\boldsymbol{u}_{k,v}^{i}, \boldsymbol{\varpi}_{k}^{i}\right) \in \mathbb{R}^{n^{i}} \times \Upsilon^{i} \text { such that } \\
\mathcal{L}^{i}\left(\boldsymbol{u}_{k,v}^{i}, \boldsymbol{\omega}^{i}\right) \leq \mathcal{L}^{i} \left(\boldsymbol{u}_{k,v}^{i}, \boldsymbol{\varpi}_{k}^{i}\right) \leq \mathcal{L}^{i} \left(\boldsymbol{v}_{v}^{i}, \boldsymbol{\varpi}_{k}^{i}\right) ~~~~ \forall \boldsymbol{v}_{v}^{i} \in \mathbb{R}^{n^{i}},~ \forall \boldsymbol{\omega}^{i} \in \Upsilon^{i}
\end{array}\right\}.
\end{equation}
Since the saddle-point problems (\ref{3:prb:saddle_problem}) is unconstrained respect to the displacement filed variables $\boldsymbol{v}_{v}^{i}$, it is easy to derive the following dual problems: 
\begin{equation}\label{3:prb:dual}
\boldsymbol{\varpi}_k^i = \arg\min_{\boldsymbol{\omega}^{i}\in\Upsilon^{i}} \frac{1}{2}\left( \boldsymbol{\omega}^{i} \right)^{\top} \mathbf{C}^{i} \boldsymbol{\omega}^{i} + \left( \boldsymbol{\omega}^{i} \right)^{\top} \mathbf{d}^{i}
\end{equation}
where matrix $\mathbf{C}^{i}$ and vector $\mathbf{d}^{i}$ defined by:
$$
\begin{aligned}
\mathbf{C}^{i}&=\left(\begin{array}{ccc}
N_{1}^{i} \left(A_{h}^{i}\right)^{-1} \left(N_{1}^{i}\right)^{\top} & N_{1}^{i} \left(A_{h}^{i}\right)^{-1} \left(T_{1}^{i}\right)^{\top} \left(I_{3}^{i}\right)^{\top} & N_{1}^{i} \left(A_{h}^{i}\right)^{-1} \left(I_{4}^{i}\right)^{\top} \\
I_{3}^{i} T_{1}^{i} \left(A_{h}^{i}\right)^{-1}  \left(N_{1}^{i}\right)^{\top} & I_{3}^{i} T_{1}^{i} \left(A_{h}^{i}\right)^{-1} \left(T_{1}^{i}\right)^{\top} \left(I_{3}^{i}\right)^{\top} & I_{3}^{i} T_{1}^{i} \left(A_{h}^{i}\right)^{-1}  \left(I_{4}^{i}\right)^{\top} \\
I_{4}^{i} \left(A_{h}^{i}\right)^{-1} \left(N_{1}^{i}\right)^{\top} & I_{4}^{i} \left(A_{h}^{i}\right)^{-1} \left(T_{1}^{i}\right)^{\top} \left(I_{3}^{i}\right)^{\top} & I_{4}^{i} \left(A_{h}^{i}\right)^{-1} \left(I_{4}^{i}\right)^{\top}
\end{array}\right),\\
&= \mathbf{B}^{i} \left(A_{h}^{i}\right)^{-1} \left(\mathbf{B}^{i}\right)^{\top} \\
\mathbf{d}^{i}&=\left(\begin{array}{c}
N_{1}^{i}\left(\lambda_{k-1,v}^{i} - \left(A_{h}^{i}\right)^{-1}\boldsymbol{b}^{i}\right)\\
I_{3}^{i}T_{1}^{i}\left(\lambda_{k-1,v}^{i} - \left(A_{h}^{i}\right)^{-1}\boldsymbol{b}^{i}\right)\\
I_{3}^{i-1} \lambda_{k-1,v}^{i-1} - I_{4}^{i}\left(A_{h}^{i}\right)^{-1}\boldsymbol{b}^{i}
\end{array}\right), ~~\text{where}~~ \mathbf{B}^{i} =  \left( \begin{array}{c}
N_1^i \\
I_3^i T_1^i \\
I_4^i
\end{array} \right).
\end{aligned}
$$
By solving the dual problems (\ref{3:prb:dual}), the optimal dual solutions $\boldsymbol{\varpi}_{k}^{i}$ can be obtained, then the unknown displacement field $\boldsymbol{u}_{k,v}^{i}$ can be deduced by 
\begin{equation}\label{3:eq:u_kv}
\boldsymbol{u}_{k,v}^{i} = \left(A_{h}^{i}\right)^{-1} \left( \boldsymbol{b}^{i} - \left( \mathbf{B}^{i} \right)^{\top} \cdot \boldsymbol{\varpi}_{k}^{i} \right),~i=1,2,\ldots,n.
\end{equation}
By bringing above formula (\ref{3:eq:u_kv}) into step 3 of Algorithm \ref{algorithm:3.3}, it can be obtained that:
\begin{equation}
\begin{aligned}
A_h^{i+1} \boldsymbol{p}_{k, v}^{i+1}&= -\frac{1}{2}\left(I_4^{i+1}\right)^{\top} \left(I_4^{i+1}\left( \mathbf{B}^{i+1} \right)^{\top} \boldsymbol{\varpi}_{k}^{i+1} +I_3^i\left( \mathbf{B}^{i} \right)^{\top} \boldsymbol{\varpi}_{k}^{i} \right),\\
A_h^i \boldsymbol{q}_{k, v}^i &= -\frac{1}{2}\left(I_3^i\right)^{\top} \left(I_4^{i+1}\left( \mathbf{B}^{i+1} \right)^{\top} \boldsymbol{\varpi}_{k}^{i+1} +I_3^i\left( \mathbf{B}^{i} \right)^{\top} \boldsymbol{\varpi}_{k}^{i} \right).
\end{aligned}
\end{equation}
Therefore, the sequence $\left\{\boldsymbol{u}_{k,v}^{i}\right\}$ in Algorithm \ref{algorithm:3.3} can be replaced by $\left\{\boldsymbol{\varpi}_{k}^{i+1}\right\}$ and the displacement field $\boldsymbol{u}_{\bar{k},v}^{i}$ is computed from $\boldsymbol{\varpi}_{\bar{k}}^{i+1}$ by (\ref{3:eq:u_kv}) for an index $\bar{k}$ determined by an appropriate stopping criterion.

Finally, the dual algebra form of layer decomposition algorithm can be illustrated as follow:

\begin{algorithm}[!ht]
\caption{Dual Algebraic Form of Discrete Layer Decomposition Algorithm}
\label{algorithm:3.4}
Let $\lambda_{0,v}^{i}\in \mathbb{R}^{n^{i}}$, $i=1,\ldots,n-1$, $\theta>0$, $\epsilon_0:=1$ and $tol$ be given.

\textbf{Output:} $\boldsymbol{u}^{i}_{k,v}$, $\boldsymbol{p}^{i}_{k,v}$, $\boldsymbol{q}^{i}_{kv,}$, $\lambda_{k,v}^{i}$

\textbf{While} $k\leq N$ and $\epsilon_k>tol$ \textbf{do}

~~~~$
\boldsymbol{\varpi}_k^i = \arg\min_{\boldsymbol{\omega}^{i}\in\Upsilon^{i}} \frac{1}{2}\left( \boldsymbol{\omega}^{i} \right)^{\top} \mathbf{C}^{i} \boldsymbol{\omega}^{i} + \left( \boldsymbol{\omega}^{i} \right)^{\top} \mathbf{d}^{i}, ~ i = 1,\ldots,n;
$

~~~~$
\left\{
\begin{aligned}
A_h^{i+1} \boldsymbol{p}_{k, v}^{i+1}&= -\frac{1}{2}\left(I_4^{i+1}\right)^{\top} \left(I_4^{i+1}\left( \mathbf{B}^{i+1} \right)^{\top} \boldsymbol{\varpi}_{k}^{i+1} +I_3^i\left( \mathbf{B}^{i} \right)^{\top} \boldsymbol{\varpi}_{k}^{i} \right);\\
A_h^i \boldsymbol{q}_{k, v}^i &= -\frac{1}{2}\left(I_3^i\right)^{\top} \left(I_4^{i+1}\left( \mathbf{B}^{i+1} \right)^{\top} \boldsymbol{\varpi}_{k}^{i+1} +I_3^i\left( \mathbf{B}^{i} \right)^{\top} \boldsymbol{\varpi}_{k}^{i} \right);
\end{aligned}
\right.
$

~~~~$I_{3}^{i}\lambda_{k,v}^i=I_{3}^{i}\lambda_{k-1,h}^{i} -\theta\left( I_{4}^{i+1}\boldsymbol{p}^{i+1}_{kh} + I_{3}^{i}\boldsymbol{q}^{i}_{kh}\right)$;

~~~~$\epsilon_k:=\frac{\sum_{i=1}^{n}\left\|\lambda_{k,v}^{i}-\lambda_{k-1,v}^{i}\right\|_{\mathbb{R}^{3 m^{i}}}}{\sum_{i=1}^{n}\left\|\lambda_{k,v}^{i}\right\|_{\mathbb{R}^{3 m^{i}}}}$ and $k = k + 1$.

\textbf{end}
\end{algorithm}
where the prescribed tolerance $tol>0$. The parameter $tol$ ensures that the algorithm will stop when the target accuracy is reached.

\bibliographystyle{siamplain}
\bibliography{references}

%% file: Main.tex
\section{Introduction}
The contact problem is a research hotspot in solid mechanics and has a wide application prospect. Especially with the introduction of variational inequality method, the theory of contact problem has experienced a rapid development \cite{3kikuchi1988contact,3han2019numerical}. 
However, when the variational inequality model and its finite element approximation are used for numerical simulation of practical problems, it must be supported by corresponding numerical algorithms. 
The difficulty of designing this kind of algorithm is to deal with the non-penetration condition and the friction condition in the contact zone \cite{3haslinger2014domain}. At present, the proposed and studied algorithms include Lagrange multiplier method \cite{3hild2010stabilized,3FRANCESCHINI2022114632}, penalty function method \cite{3zang2011contact}, augmented Lagrange method \cite{3simo1992augmented}, semi-smooth Newton method \cite{3ito2003semi}, Neumann-Dirichlet algorithm \cite{3bayada2008convergence}, interior point method \cite{3kuvcera2013interior}, nonsmooth multiscale method \cite{3krause2009nonsmooth}, primal-dual active set algorithm \cite{3hueber2008primal, 3abide2022unified}, Nitsche’s methods \cite{3BEAUDE2023116124} and so on. But, most of the problems solved by these algorithms can only involve a single elastic body. 
With the complexity of application scenarios and the development of variational inequality theory, the contact problem between elastic bodies has been paid more and more attention. For example, domain decomposition method \cite{3haslinger2014domain}, semi-smooth Newton method and non-smooth multi-scale method are gradually developed to deal with the contact problem of multi-elastic domains.
However, due to the limitations of the actual physical model, the above methods only consider the design and research of the algorithm for the contact problem of two elastic bodies. Therefore, according to the specific physics research background, it is necessary to construct the corresponding variational inequality model and design the efficient finite element algorithm to solve the practical problems.

The mechanical model of pavement has always been the core issue in the field of transportation \cite{3nawaz2013soil,3brown1996soil,3huang2004pavement}. In theory, it is the intersection of materials science, mechanics and mathematics, so it has been widely concerned by researchers \cite{3liu2011coupled}. In terms of application, the mechanical model of pavement is closely related to the paving, maintenance and service life prediction of pavement, so it is directly related to the development strategy of cities and even countries \cite{3brownjohn2007structural}. 
Therefore, for the study of such an important problem, it is necessary to construct the theoretical framework of mathematical model to understand and guide the research of pavement mechanics.
With the development of pavement laying technology and computer computing power, the theory of mathematical model is also constantly improving. For example, the elastic half-space model \cite{3hung2001elastic}, potential function \cite{3burmister1945general}, rheology theory \cite{3yusoff2011modelling} have been concerned and studied. 
In particular, partial differential equation model also provides an important way to simulate pavement mechanics from the perspective of continuum mechanics \cite{3zhang2016weak}. With the deepening of numerical theory research such as finite element, this model has become a very important reference standard in the study of pavement mechanics. And with the improvement of calculation accuracy, the contact condition between pavement layers has gradually become an important factor in the pavement mechanics model, which also provides a theoretical basis and development space for the application of variational inequality model \cite{3raposeiras2012influence, 3zokaei2014finite}. 
The existing pavement experiments and field observations have proved that the material difference of pavement layer and the discontinuity of contact conditions make the numerical simulation of pavement as a whole inconsistent with the actual pavement conditions, thus affecting the fitting accuracy of numerical results \cite{3ma2021analytical, 3kim2011numerical}. 
The variational inequality model developed on the basis of partial differential equation model has great advantages in dealing with contact problems, so it has great potential research value in improving the precision of numerical simulation.

In the previous work, a class of multi-layer elastic contact systems with interlayer contact conditions has been constructed for the characteristics of pavement physical models \cite{3zhang2022variational, 3zhangvariational}. Then, a partial differential equation model and variational inequality model based on this system are proposed, and the existence and uniqueness of their solutions and the convergence of numerical solutions under the finite element method are proved successively. 
However, only the theoretical research is not enough to make the model be used in the actual mechanical model of road surface, we must design an effective algorithm to verify that the relevant variational inequality theory can play a practical function.
Inspired by the Neumann-Neumann domain decomposition algorithm \cite{3haslinger2014domain}, a layer decomposition algorithm suitable for multi-layer elastic contact systems with interlayer Tresca friction conditions is proposed and studied. The proposed algorithm not only solves the finite element numerical solution problem of variational inequality model corresponding to multi-layer elastic contact system, but also proves the feasibility of applying variational inequality model to pavement mechanics modeling through practical examples, which lays a foundation for the subsequent application in pavement mechanics research.

The paper is organized as follows. 
In the Section 2, the pavement physical model is introduced, and then the corresponding partial differential equation model and variational inequality model are presented successively. 
In the Section 3, a layer decomposition algorithm is presented, and based on the fixed point theorem, the convergence condition of the algorithm is proved. 
In order to verify that the algorithm is also convergent in finite element space, in Section 4, the discrete form of the algorithm in piecewise linear polynomial space is presented and its convergence theorem is proved. 
Based on the discrete form of the algorithm, in Section 5, two numerical experiments simulating real pavement stress models are performed, and the nephograms of the numerical solution of displacement field of the pavement are presented and analyzed.
Finally, in Section 6, this study is summarized and some perpesctives are stated.

\section{Setting of the problem}

\subsection{PDE model and variational inequality for multi-layer elastic contact systems}

\begin{figure}[!t]
\centering
\begin{minipage}{0.45\linewidth}
\centering
\begin{tikzpicture}[3d view={-45}{22}, scale=0.5]
\draw[->] (0,0,0) -- (0.5,0,0) node[pos=1,right]{x};
\draw[->] (0,0,0) -- (0,0.5,0) node[pos=1,left]{y};
\draw[->] (0,0,0) -- (0,0,0.5) node[pos=1,above]{z};
\foreach \n in {0,1} 
{
\draw (0,0,0+\n) -- (4,0,0+\n) -- (4,8,0+\n) -- (0,8,0+\n) -- (0,0,0+\n);
\draw (0,0,1+\n) -- (4,0,1+\n) -- (4,8,1+\n) -- (0,8,1+\n) -- (0,0,1+\n);
\draw (0,0,0+\n) -- (0,0,1+\n);
\draw (4,0,0+\n) -- (4,0,1+\n);
\draw (4,8,0+\n) -- (4,8,1+\n);
\draw (0,8,0+\n) -- (0,8,1+\n);
}
\foreach \n in {0,1} 
{
\draw (0,0,3+\n) -- (4,0,3+\n) -- (4,8,3+\n) -- (0,8,3+\n) -- (0,0,3+\n);
\draw (0,0,4+\n) -- (4,0,4+\n) -- (4,8,4+\n) -- (0,8,4+\n) -- (0,0,4+\n);
\draw (0,0,3+\n) -- (0,0,4+\n);
\draw (4,0,3+\n) -- (4,0,4+\n);
\draw (4,8,3+\n) -- (4,8,4+\n);
\draw (0,8,3+\n) -- (0,8,4+\n);
}
\foreach \n in {0.1,0.2,...,0.9}
{ \fill (0,0,2+\n) circle (0.5pt);
\fill (4,0,2+\n) circle (0.5pt);
\fill (4,8,2+\n) circle (0.5pt);
\fill (0,8,2+\n) circle (0.5pt);
}
\begin{scope}[canvas is xy plane at z=5]
\node at (2,4) [color = blue, transform shape,scale=1.5] {$\Gamma_2^1$};
\end{scope}
\filldraw[fill opacity=0.1,fill=blue] (0,0,5) -- (4,0,5) -- (4,8,5) -- (0,8,5) -- cycle ;
\draw [color=blue] (0,0,5) -- (4,0,5) -- (4,8,5) -- (0,8,5) -- (0,0,5);
\node at (4,0,4.5) [color = blue,scale=0.8,right] {$\Omega^{1}$};

\begin{scope}[canvas is xy plane at z=4]
\node at (2,4) [color = purple, transform shape,scale=1.5] {$\Gamma_2^2/\Gamma_3^1$};
\end{scope}
\filldraw[fill opacity=0.1,fill=purple] (0,0,4) -- (4,0,4) -- (4,8,4) -- (0,8,4) -- cycle ;
\draw [color=purple] (0,0,4) -- (4,0,4) -- (4,8,4) -- (0,8,4) -- (0,0,4);
\node at (4,0,3.5) [color = purple, scale=0.8,right] {$\Omega^{2}$};

\begin{scope}[canvas is xy plane at z=3]
\node at (2,4) [color=orange, transform shape,scale=1.5] {$\Gamma_2^3/\Gamma_3^2$};
\end{scope}
\filldraw[fill opacity=0.1,fill=orange] (0,0,3) -- (4,0,3) -- (4,8,3) -- (0,8,3) -- cycle ;
\draw [color=orange] (0,0,3) -- (4,0,3) -- (4,8,3) -- (0,8,3) -- (0,0,3);

\begin{scope}[canvas is xy plane at z=2]
\node at (2,4) [color=cyan, transform shape,scale=1.5] {$\Gamma_2^{n-1}/\Gamma_3^{n-2}$};
\end{scope}
\filldraw[fill opacity=0.1,fill=cyan] (0,0,2) -- (4,0,2) -- (4,8,2) -- (0,8,2) -- cycle ;
\draw [color=cyan] (0,0,2) -- (4,0,2) -- (4,8,2) -- (0,8,2) -- (0,0,2);
\node at (4,0,1.5) [color = cyan, scale=0.8,right] {$\Omega^{n-1}$};

\begin{scope}[canvas is xy plane at z=1]
\node at (2,4) [color = gray, transform shape,scale=1.5] {$\Gamma_2^{n}/\Gamma_3^{n-1}$};
\end{scope}
\filldraw[fill opacity=0.1,fill=gray] (0,0,1) -- (4,0,1) -- (4,8,1) -- (0,8,1) -- cycle ;
\draw [color=gray] (0,0,1) -- (4,0,1) -- (4,8,1) -- (0,8,1) -- (0,0,1);
\node at (4,0,0.5) [color = gray, scale=0.8,right] {$\Omega^{n}$};

\begin{scope}[canvas is xy plane at z=0]
\node at (2,4) [color = teal, transform shape,scale=1.5] {$\Gamma_1^{n}(\Gamma_3^{n})$};
\end{scope}
\filldraw[fill opacity=0.1,fill=teal] (0,0,0) -- (4,0,0) -- (4,8,0) -- (0,8,0) -- cycle ;
\draw [color=teal] (0,0,0) -- (4,0,0) -- (4,8,0) -- (0,8,0) -- (0,0,0);

\begin{scope}[canvas is yz plane at x=0]
\node at (1,0.5) [rotate=180,yscale=-1,transform shape,scale=1.2] {$\Gamma_1^{n}$};
\node at (1,1.5) [rotate=180,yscale=-1,transform shape,scale=1.2] {$\Gamma_1^{n-1}$};
\node at (1,3.5) [rotate=180,yscale=-1,transform shape,scale=1.2] {$\Gamma_1^{2}$};
\node at (1,4.5) [rotate=180,yscale=-1,transform shape,scale=1.2] {$\Gamma_1^{1}$};
\end{scope}
\begin{scope}[canvas is xz plane at y=0]
\node at (3,0.5) [xscale=1,transform shape,scale=1.2] {$\Gamma_1^{n}$};
\node at (3,1.5) [xscale=1,transform shape,scale=1.2] {$\Gamma_1^{n-1}$};
\node at (3,3.5) [xscale=1,transform shape,scale=1.2] {$\Gamma_1^{2}$};
\node at (3,4.5) [xscale=1,transform shape,scale=1.2] {$\Gamma_1^{1}$};
\end{scope}
\begin{scope}[canvas is yz plane at x=4]
\node at (7,0.5) [rotate=180,yscale=-1,transform shape,scale=1.2] {$\Gamma_1^{n}$};
\node at (7,1.5) [rotate=180,yscale=-1,transform shape,scale=1.2] {$\Gamma_1^{n-1}$};
\node at (7,3.5) [rotate=180,yscale=-1,transform shape,scale=1.2] {$\Gamma_1^{2}$};
\node at (7,4.5) [rotate=180,yscale=-1,transform shape,scale=1.2] {$\Gamma_1^{1}$};
\end{scope}
\begin{scope}[canvas is xz plane at y=8]
\node at (1,0.5) [xscale=1,transform shape,scale=1.2] {$\Gamma_1^{n}$};
\node at (1,1.5) [xscale=1,transform shape,scale=1.2] {$\Gamma_1^{n-1}$};
\node at (1,3.5) [xscale=1,transform shape,scale=1.2] {$\Gamma_1^{2}$};
\node at (1,4.5) [xscale=1,transform shape,scale=1.2] {$\Gamma_1^{1}$};
\end{scope}
\foreach \x in {1.5,1.75,2.0,2.25,2.5} \foreach \y in {3.5,3.75,4.0,4.25,4.5} {
\draw[->,color=red, line width=0.1pt] (\x,0.15+\y,5.4)--(\x,\y,5);}
\node at (2,4,6) [color=red,scale=0.8] {$f_2$};
\end{tikzpicture}
\caption{The physical model of multi-layer contact problem}
\label{3:fig1}
\end{minipage}
\begin{minipage}{0.45\linewidth}
\centering
\includegraphics[width=1\linewidth]{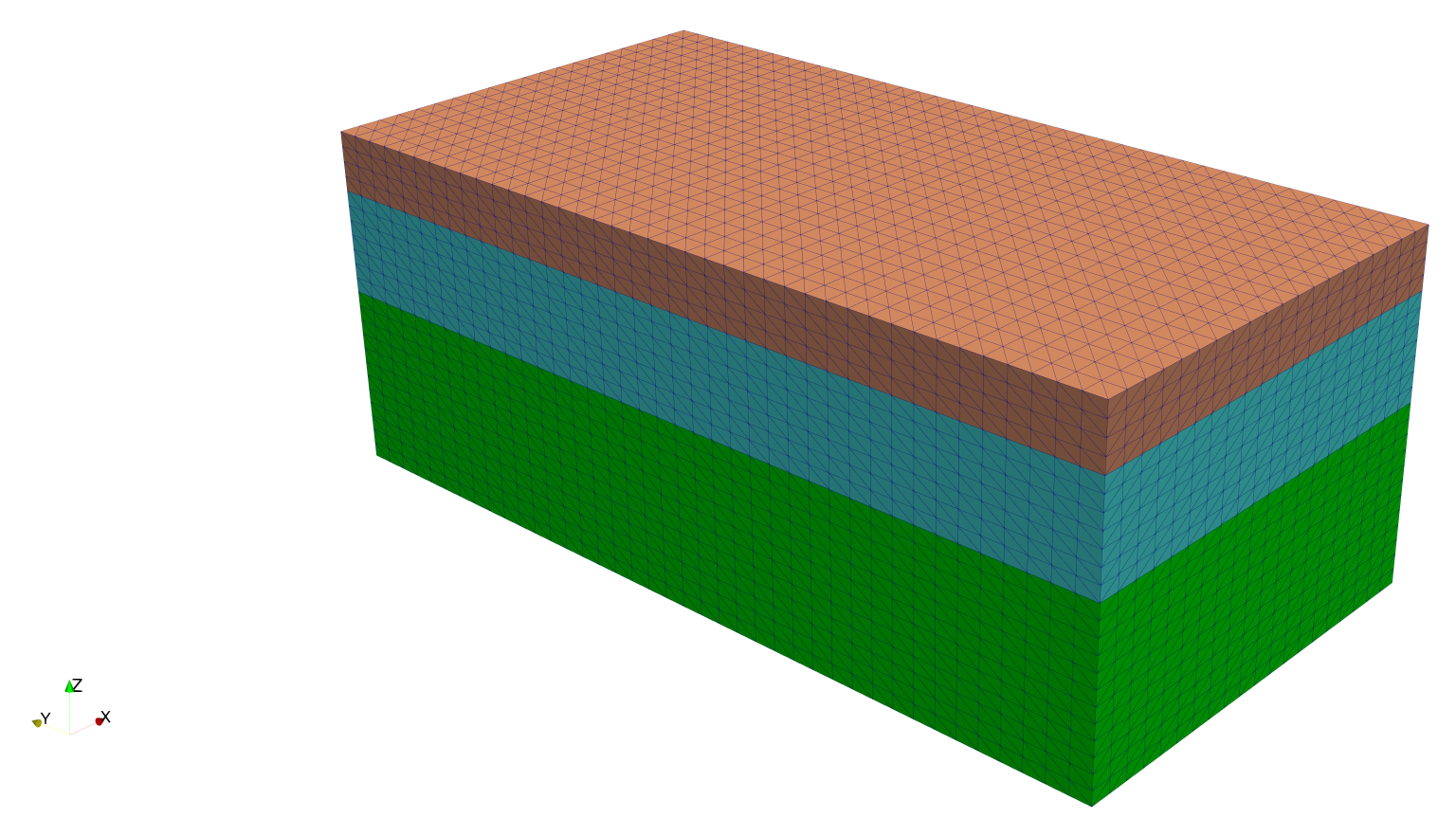}
\caption{The physical model and meshing of three-layer elastic system.}
\label{3:fig2}
\end{minipage}
\end{figure}

In previous works \cite{3zhang2022variational,3zhangvariational}, the variational inequalities for multi-layer elastic contact and visoelastic systems have been discussed. The physical model discussed this time continues previous research work, as shown in Fig \ref{3:fig1}. In the physical model, the open region occupied by the elastic system is denoted as $\Omega=\cup_{i=1}^{n}\Omega^{i}\subset\mathbb{R}^{d}$ ($d=3$), and the open Lipschitz boundary is denoted as $\Gamma^{i}=\cup_{j=1}^{3} \bar{\Gamma}_{j}^{i}$, where the displacement field on $\Gamma_{1}^{i}$ is $\boldsymbol{0}$. Let $\Gamma_{j}=\cup_{i=1}^{n}\bar{\Gamma}_{j}^{i}$ and unit outward normal vector at $x\in\Gamma_{j}^{i}$ be $\boldsymbol{\nu}^i$, where $\boldsymbol{\nu}^i=\boldsymbol{\alpha}^{i}$ and $\boldsymbol{\beta}^{i}$ when $j=2$ and $3$, respectively. Then, the frictional contact boundary is defined as $\Gamma_{c}^{i}=\Gamma_{3}^{i}\cap \Gamma_{2}^{i+1}$. The surface force on the boundary $\Gamma_{2}^{1}$ is denoted by $\boldsymbol{f}_{2}$, and the body force acting on the elastic body $\Omega^i$ is denoted by $\boldsymbol{f}_{0}^{i}$. 
For illustrating the characteristics of contact boundaries, the following decompositions are defined:
\begin{small}
\begin{align*}
&v_{\beta}^{i}=\boldsymbol{v}^{i} \cdot \boldsymbol{\beta}^{i},~  \boldsymbol{v}_{\eta}^{i}=\boldsymbol{v}^{i}-v_{\beta}^{i}\cdot\boldsymbol{\beta}^{i}, ~ v_{\alpha}^{i}=\boldsymbol{v}^{i} \cdot \boldsymbol{\alpha}^{i},~ \boldsymbol{v}_{\tau}^{i}=\boldsymbol{v}^{i}-v_{\alpha}^{i}\cdot\boldsymbol{\alpha}^{i}, \\
&\sigma_{\beta}^{i}=\boldsymbol{\beta}^{i}\cdot\boldsymbol{\sigma}^{i} \cdot \boldsymbol{\beta}^{i}, ~\boldsymbol{\sigma}_{\eta}^{i}=\boldsymbol{\sigma}^{i} \cdot \boldsymbol{\beta}^{i}-\sigma_{\beta}^{i}\cdot\boldsymbol{\beta}^{i}, ~
\sigma_{\alpha}^{i}=\boldsymbol{\alpha}^{i}\cdot\boldsymbol{\sigma}^{i} \cdot \boldsymbol{\alpha}^{i},~ \boldsymbol{\sigma}_{\tau}^{i}=\boldsymbol{\sigma}^{i} \cdot \boldsymbol{\alpha}^{i}-\sigma_{\alpha}^{i}\cdot\boldsymbol{\alpha}^{i}.
\end{align*}
\end{small}
The above decomposition realizes that the displacement and stress are projected on the normal vector and tangent plane of the boundary respectively. Then, let $v_{N}^{i}=v_{\beta}^{i}$, $\boldsymbol{v}_{T}^{i}=\boldsymbol{v}_{\eta}^{i}$, $\sigma_{N}^{i}=\sigma_{\beta}^{i}$ and $\boldsymbol{\sigma}_{T}^{i}=\boldsymbol{\sigma}_{\eta}^{i}$ on $\Gamma_{c}^{i}$. The jump operator $[\cdot]$ is defined as:
\begin{small}
\begin{align*}
\left[{v}_{N}^{i}\right]=\boldsymbol{v}^{i} \cdot \boldsymbol{\beta}^{i}+\boldsymbol{v}^{i+1} \cdot \boldsymbol{\alpha}^{i+1},~~
\left[\boldsymbol{v}_{T}^{i}\right]=\boldsymbol{v}_{\eta}^{i} - \boldsymbol{v}_{\tau}^{i+1}.
\end{align*}
\end{small}
Note that the non-penetration condition $\left[{v}_{N}^{i}\right]\leqslant 0$ is satisfied on contact zone $\Gamma_{c}^{i}$ ($i=1,2,\ldots n-1$).
Therefore, the partial differential equation model of the multi-layer elastic contact system is as follows:

\begin{problem}[$P_{0}$]\label{prb:3:p_0}
Find a displacement field $\boldsymbol{u}^{i}:\Omega^{i} \rightarrow \mathbb{R}^{d}$ and the stress field $\boldsymbol{\sigma}^{i}: \Omega^{i} \rightarrow \mathbb{S}^{d}$ ($i=1,2,\ldots,n$) such that:
\begin{small}
\begin{align}
&\boldsymbol{\sigma}^{i}=\mathcal{A}^{i} \boldsymbol{\varepsilon}(\boldsymbol{u}^{i}), ~\boldsymbol{\varepsilon}(\boldsymbol{u}^{i})=\frac{1}{2}(\nabla \boldsymbol{u}^{i}+(\nabla \boldsymbol{u}^{i})^{T}) && \text { in } \Omega^{i}, \label{3:pde.1}\\
&\operatorname{Div} \boldsymbol{\sigma}^{i}+\boldsymbol{f}_{0}^{i}=\mathbf{0} && \text { in } \Omega^{i},\label{3:pde.2}\\
&\boldsymbol{u}^{i}=\mathbf{0} && \text { on } \Gamma_{1}^{i},\Gamma_{3}^{n} \label{3:pde.3}\\
&\boldsymbol{\sigma}^{1} \cdot \boldsymbol{\alpha}^{1}=\boldsymbol{f}_{2} && \text { on } \Gamma_{2}^{1}, \label{3:pde.4}\\
&\sigma^{i}_{\alpha} = -\sigma^{i-1}_{\beta},~
\boldsymbol{\sigma}_{\tau}^{i} = \boldsymbol{\sigma}_{\eta}^{i-1},
&&\text { on } \Gamma_{2}^{i}, ~i\ne 1, \label{3:pde.5}\\
&\left\{
\begin{aligned}
& \left[{u}_{N}^{i}\right] \leqslant 0,~ \sigma^{i}_{N}\leqslant 0,~ \sigma^{i}_{N}\cdot[u_{N}^{i}] = 0,~ |\boldsymbol{\sigma}_{T}^{i}|\leqslant g^{i}(x)\\
& |\boldsymbol{\sigma}_{T}^{i}| < g^{i}(x) \Rightarrow
[\boldsymbol{u}_{T}^{i}] = 0\\
& |\boldsymbol{\sigma}_{T}^{i}| = g^{i}(x) \Rightarrow
\boldsymbol{\sigma}_{T}^{i} = -\lambda[\boldsymbol{u}_{T}^{i}], ~\lambda\geqslant 0
\end{aligned}
\right. &&  \text { on } \Gamma_{c}^{i},~i\ne n,\label{3:pde.6}
\end{align}
where $g^{i}(x) \geqslant 0$ ($x\in \Gamma^i_c$) is the bounded interlayer contact function.
\end{small}
\end{problem}
In the above PDEs, (\ref{3:pde.1}) and (\ref{3:pde.2}) are the constitutive equation and the equilibrium equation respectively, where $\mathcal{A}^{i}$ is the elastic operator. (\ref{3:pde.3}) and (\ref{3:pde.4}) represent the displacement boundary conditions and stress boundary conditions of the system. This boundary condition truly simulates the stress on the pavement when the vehicle is driving. Finally, (\ref{3:pde.5}) and (\ref{3:pde.6}) represent the Tresca frictional contact conditions on $\Gamma^{i}_{c}$. It can be found that the stresses of the elastic bodies on both sides of the contact surface are equal in magnitude and opposite in direction, and the contact function $g^{i}(x)$ represents the upper bound of the tangential stress.

And the variational inequality problem equivalent to this PDE model is as follows:
\begin{problem}[$P_{0}^{v}$]\label{prb:3.2:P_0_v}
Find a displacement $\boldsymbol{u} \in \mathcal{K}$ which satisfies:
\begin{small}
\begin{equation}\label{ieq:prb:3.2:var}
a(\boldsymbol{u}, \boldsymbol{v}-\boldsymbol{u})+j( \boldsymbol{v}) -j( \boldsymbol{u}) \geqslant L(\boldsymbol{v}-\boldsymbol{u}),~ \forall \boldsymbol{v} \in \mathcal{K},
\end{equation}
\end{small}
where
\begin{small}
\begin{align*}
& a(\boldsymbol{v}, \boldsymbol{w})=\sum_{i=1}^n a^i\left(\boldsymbol{v}^i, \boldsymbol{w}^i\right) =\sum_{i=1}^n\int_{\Omega^{i}} \mathcal{A}^{i} \boldsymbol{\varepsilon}(\boldsymbol{v}^{i}): \boldsymbol{\varepsilon}(\boldsymbol{w}^{i}) \mathrm{d} x,\\
& L(\boldsymbol{v})=\sum_{i=1}^n L^{i}(\boldsymbol{v}^{i}) = \sum_{i=1}^n \int_{\Omega^i} \boldsymbol{f}_0^i \cdot \boldsymbol{v}^i d x+\int_{\Gamma_2^1} \boldsymbol{f}_2 \cdot \boldsymbol{v}^1 d s,\\
& j( \boldsymbol{w})=\sum_{i=1}^{n-1} j^i\left( \boldsymbol{w}^i, \boldsymbol{w}^{i+1}\right)=\sum_{i=1}^{n-1} \int_{\Gamma_c^i} g^i(x)\left|\left[\boldsymbol{w}_T^i\right]\right| d s, \quad \forall \boldsymbol{w} \in V,\\
& \mathcal{K}=\left\{\boldsymbol{v} \in V \mid\left[v_N^i\right] \leqslant 0 \text {, a.e. on } \Gamma_c^i, i=1,2, \ldots, n-1\right\} \subset V.
\end{align*}
\end{small}
\end{problem}
In the above problem, vector space and second-order tensor space are defined as follows:
\begin{small}
$$
\begin{aligned}
V^{i} &=\left\{\boldsymbol{v}^{i}=\left(v^{i}_{k}\right) \in\left(H^{1}(\Omega^{i})\right)^{d} ~\Big|~ \boldsymbol{v}^{i}=\mathbf{0} \text { on } \Gamma^{i}_{1} \right\} , ~V = V^{1} \times \cdots \times V^{n},\\
Q^{i} &=\left\{\boldsymbol{\tau}^{i}=\left(\tau^{i}_{kl}\right) \in\left(L^{2}(\Omega^{i})\right)^{d \times d}~\Big|~ \tau^{i}_{lk}=\tau^{i}_{lk}, 1 \leqslant k, l \leqslant d\right\}, \\
Q^{i}_{1} &=\left\{\boldsymbol{\tau}^{i} \in Q^{i}~\Big|~ \operatorname{Div} \boldsymbol{\tau}^{i} \in\left(L^{2}(\Omega^{i})\right)^{d}\right\}, ~Q_{1}= Q_{1}^{1} \times \cdots \times Q_{1}^{n}.
\end{aligned}
$$
\end{small}
It is easy to verify that the above spaces are all Hilbert spaces, so the norms defined by standard inner product are recorded as $\|\cdot\|_{H^{1}\left(\Omega^{i}\right)^{d}}$, $\|\cdot\|_{H^{1}}$, $\|\cdot\|_{L^{2}\left(\Omega^{i}\right)^{d}}$, $\|\cdot\|_{L^{2}}$. 
Then, based on $\operatorname{meas}_{d-1}\left(\Gamma_{1}^{i}\right)>0$ and Korn's inequality \cite{3kikuchi1988contact}, the following inner products on $V^i$ and $V$ can be introduced:
\begin{small}
$$
(\boldsymbol{u}^{i}, \boldsymbol{v}^{i})_{V^{i}} =(\boldsymbol{\varepsilon}(\boldsymbol{u}^{i}), \boldsymbol{\varepsilon}(\boldsymbol{v}^{i}))_{L^{2}\left(\Omega^{i}\right)^{d}},~(\boldsymbol{u}, \boldsymbol{v})_{V} = \sum_{i=1}^{n} (\boldsymbol{\varepsilon}(\boldsymbol{u}^{i}), \boldsymbol{\varepsilon}(\boldsymbol{v}^{i}))_{L^{2}\left(\Omega^{i}\right)^{d}}, \forall \boldsymbol{u}^{i}, \boldsymbol{v}^{i} \in V^{i},
$$
\end{small}
and the norm $\|\cdot\|_{V^{i}}$ induced by this inner product is equivalent to $\|\cdot\|_{H^{1}(\Omega^{i})^{d}}$ in $V^{i}$.

Furthermore, according to the trace theorem on Sobolev spaces \cite{3adams2003sobolev}, there exists a constant $c_{t}^{i}$ that depends only on $\Omega^{i}$ and $\Gamma^{i}$, such that the following inequality holds:
\begin{small}
\begin{equation}\label{ieq:3:trace}
\|\gamma^{i}_{j}\boldsymbol{v}^{i}\|_{L^{2}\left(\Gamma_{c}^{j}\right)^{d}} \leqslant c_{t}^{i}\|\boldsymbol{v}^{i}\|_{V^{i}},~ \forall \boldsymbol{v}^{i} \in V^{i}, j=i \text{ or }i-1,
\end{equation}
\end{small}
where $\gamma^{i}_{j}: V^{i}\to L^2\left(\Gamma_{c}^{j} \right)$ ($j=i-1$ or $i$) is the usual trace operator. Moreover, $c_{t}^{\min}$ and $c_{t}^{\max}$ can be defined as 
$c_{t}^{\min} = \min\{c_{t}^{1},\ldots,c_{t}^{n-1}\}$ and $c_{t}^{\max} = \max\{c_{t}^{1},\ldots,c_{t}^{n-1}\}$.

\begin{remark}\label{3:rem:1}
It is worth noting that the Tresca frictional contact conditions in Problem \ref{prb:3:p_0} are usually only applicable to bilateral contact problems, because the friction law requires that the friction region is not changed during deformation. In the unilateral contact problem of this model, the friction contact conditions usually satisfy the Coulomb friction law \cite{3zhang2022variational}, where contact function $g^{i}(x)$ replaced by $g^{i}(\sigma^{i}_{N})$.
However, because the theoretical analysis and algorithm design of the multi-layer elastic contact system with Coulomb friction law are more difficult, how to approximate it by a simpler model has a strong practical significance. The existing studies show that by designing and solving a fixed point algorithm which computes a unilateral contact problem with Tresca friction law at each step, the unilateral contact problem with Coulomb friction law can be approximately solved \cite{3eck1998existence,3laborde2008fixed}. So, the unilateral contact problem with Tresca friction law has also become a new research direction in the contact problem. In addition, in the study of pavement mechanics, the interlayer contact model lacks theoretical standards, so it is necessary to try multiple contact conditions and carry out comparative studies. Therefore, it is of great significance to study the model of multi-layer elastic contact system with Tresca friction law from both theoretical and application points of view.
\end{remark}

Under the framework of variational inequalities, the constraints of the original problem $P_{0}$ are greatly simplified. Our previous research has proved that when the above elastic operators $\mathcal{A}^i$ and contact functions $g^i$ satisfy the some regularity conditions, the solution of the variational inequality (\ref{ieq:prb:3.2:var}) of the multi-layer elastic contact system is unique, and the numerical solution obtained by the finite element method also satisfies the convergence theorem. However, the above conclusions remain in the theoretical level, and the algorithm for solving the finite element numerical solution of the multi-layer elastic contact system has not been proposed. 
Therefore, based on the Neumann-Neumann domain decomposition algorithm \cite{3haslinger2014domain}, a layer decomposition algorithm suitable for this problem is proposed.
Thereafter, the elastic operator $\mathcal{A}^i$ is restricted to a linearly symmetric fourth-order tensor, and thus $a(\cdot,\cdot)$ is commutative bilinear form.

\begin{remark}
In Ref.\cite{3zhang2022variational}, a multi-layer elastic system with a generalized version of Coulomb friction is studied, and the existence and uniqueness of the solution of the corresponding variational inequality and the convergence of the finite element numerical solution are proved. It is not difficult to judge that the multi-layer elastic contact system with Tresca friction law in this study is a simplification of the contact problem with Coulomb friction. Therefore, through the same method, it can be verified that the existence and uniqueness theorem of the exact solution and the convergence theorem of the numerical solution in the Ref.\cite{3zhang2022variational} are also established in problem $P^v_0$.
\end{remark}

\subsection{Equivalent formulation of (\texorpdfstring{$P_{0}^{v}$}{})}

First, the space of functions on the contact boundary $\Gamma_{c}^{i}(i=1,\ldots n-1)$ needs to be defined. Next we shall suppose that $\Gamma_1^i \cap \Gamma_c^{i}=\Gamma_1^{i+1} \cap \Gamma_c^{i}$. Let
\begin{small}
\begin{equation}\label{def:3:H}
\left\{\begin{aligned}
&H^{1/2}\left(\Gamma_c^{i}\right)=\left\{\varphi \in \left(L^2\left(\Gamma_c^{i} \right) \right)^{d} ~\Big|~ \exists \boldsymbol{v}^{j} \in V^{j},~ \boldsymbol{v}^{j}=\varphi \text { on } \Gamma_c^{i},~j=i \text{ or } i+1\right\} \\
&{H}^{1/2}\left(\Gamma_c\right) = H^{1/2}\left(\Gamma_c^{1}\right) \times \cdots \times H^{1 / 2}\left(\Gamma_c^{n-1}\right)
\end{aligned}\right.
\end{equation}
\end{small}
Based on the assumptions of the collective position of $\Gamma_1^i$ and $\Gamma_c^{i}$, the definition of $H^{1/2}\left(\Gamma_c^{i}\right)$ is independent of the choice of $j \in\{i,i+1\}$. The standard norm of $H^{1/2} \left( \Gamma_c^{i} \right)$ is denoted by $\|\varphi^{i}\|_{1/2, \Gamma_c^{i}}$,
$\varphi^i=\left(\varphi_1^{i},\ldots,\varphi_d^{i}\right)\in H^{1/2}\left(\Gamma_c^i\right)$. Let $\chi^{i}_{j}:{H}^{1/2}\left(\Gamma_c^{i}\right) \mapsto {V}^{j},~j=i \text{ or } i+1$ be the extension mapping defined by
\begin{small}
\begin{equation}\label{def:3:chi}
\left.\begin{array}{rl}
\chi^{i}_{j} \varphi^{i} \in {V}^j: & a^j\left(\chi^i_j \varphi^{i}, \boldsymbol{v}^j\right)=0 \quad \forall \boldsymbol{v}^j \in {V}_{i-j+3}^j(0), \\
& \chi^{i}_{j} \varphi^{i} =\varphi^{i} \text { on } \Gamma_c^{i}, ~ j =i \text{ or } i+1
\end{array}\right\},
\end{equation}
\end{small}
where the space $V^{i}_{s}(0)$ is defined as: 
\begin{small}
$$
V^{i}_{s}(0)=\left\{\boldsymbol{v}^{i}\in \left(H^{1}(\Omega^{i})\right)^{d} ~\Big|~ \boldsymbol{v}^{i}=0,~ \text{ on } \Gamma^{i}_{1} \text{ and } \Gamma^{i+s-3}_{c}\right\},~s=2 \text{ or } 3.
$$
\end{small}
Then, the equivalent norms in ${H}^{1/2}\left(\Gamma_c^i\right)$ can be defined by
\begin{small}
\begin{equation}\label{def:3:norm_H_c_2}
\| \varphi^{i} \|_{\Gamma_c^{i,j}}:=\inf_{\boldsymbol{v}^j \in V^j} \left\{\left\|\boldsymbol{v}^j\right\|_{1, \Omega^j} \Big| \boldsymbol{v}^j=\varphi^{i} \text { on } \Gamma_c^{i}\right\}=\| \chi^i_{j} \varphi^{i} \|_{1, \Omega^j},~ j=i \text{ or } i+1,
\end{equation}
\end{small}
where $\|\cdot\|_{1,\Omega^j}$ stands for the energy norm in ${V}^j$ : $\left\| \boldsymbol{v}^j \right\|_{1, \Omega^j}^2:=a^j\left(\boldsymbol{v}^j, \boldsymbol{v}^j\right)$, $\forall\boldsymbol{v}^j \in {V}^j$.
The norms $\|\cdot\|_{1/2, \Gamma_c^{i}}$ and $\|\cdot\|_{\Gamma_c^{i,j}}, j=i,i+1$ are equivalent in ${H}^{1/2}\left(\Gamma_c^{i}\right)$, which satisfies the following requirements:
\begin{small}
$$
c^{i}_{1} \|\varphi^i\|_{\Gamma_c^{i,i+1}} \leqslant \|\varphi^i\|_{\Gamma_c^{i,i}} \leqslant c^{i}_{2} \|\varphi^i\|_{\Gamma_c^{i,i+1}},~ \forall \varphi^i \in H^{1/2}\left(\Gamma_c^i\right), ~ i=1\ldots,n-1.
$$
\end{small}

Let $\phi|_{\Gamma_c^{i}} \in L^2\left(\Gamma_c^{i}\right)$ be given and denote
\begin{small}
$$
{K}^{i}(\phi)=\left\{\boldsymbol{v}^i \in {V}^i ~\Big|~ \boldsymbol{v}^i \cdot \boldsymbol{\beta}^i \leq-\phi \text { on } \Gamma_c^{i}\right\}, ~~i=1,\ldots,n-1
$$
\end{small}
the convex, closed subset of ${V}^i$. Furthermore, the space $V^{i}_{s}(\boldsymbol{w})$ is defined as: 
\begin{small}
$$
V^{i}_{s}(\boldsymbol{w}^{j})=\left\{\boldsymbol{v}^{i}\in V^{i} ~\Big|~ \left\|\gamma^{i}_{k}\boldsymbol{v}^{i}- \gamma^{j}_{k}\boldsymbol{w}^{j}\right\|_{\Gamma^{k,k+1}_{c}}=0,~k=i+s-3 \right\},
$$
\end{small}
in which $s=2$ or $3$, $j=i+s-3$ or $i+s-2$ and if $i=1$ and $s=2$ or $i=n$ and $s=3$, we define that $V^{1}_{2}(\boldsymbol{w}) = V^{1}$ and $V^{n}_{3}(\boldsymbol{w}) = V^{n}$ for any $\boldsymbol{w}$, respectively.

In following proposition, it can be verified that problem $P_{0}^v$ is equivalent to the form of $n$ coupled problems, where each of them is formulated in its own region $\Omega^{i}$.

\begin{proposition}\label{prp:3.1}
A function $\boldsymbol{u}=(\boldsymbol{u}^1,\ldots,\boldsymbol{u}^{n}) \in \mathcal{K}$ is a solution of $P_{0}^v$ if and only if $\boldsymbol{u}^{i}\in K^{i}(\boldsymbol{u}^{i+1}\cdot\boldsymbol{\alpha}^{i+1})$, $i=1,\ldots,n-1$ and $\boldsymbol{u}^{n}\in V^{n}$ solve the following inequalities
\begin{small}
\begin{equation}\label{ieq:prp:3.1:1}
\left\{
\begin{aligned}
& a^{1}(\boldsymbol{u}^{1}, \boldsymbol{v}^{1}-\boldsymbol{u}^{1}) + j^{1}( \boldsymbol{v}^{1}, \boldsymbol{u}^{2}) -j^{1}( \boldsymbol{u}^{1}, \boldsymbol{u}^{2}) \geqslant L^{1}(\boldsymbol{v}^{1} -\boldsymbol{u}^{1}),\\
& a^{i}\left(\boldsymbol{u}^{i}, \boldsymbol{v}^{i}-\boldsymbol{u}^{i}\right) + j^{i}\left( \boldsymbol{v}^{i}, \boldsymbol{u}^{i+1}\right) - j^{i}\left( \boldsymbol{u}^{i}, \boldsymbol{u}^{i+1}\right) \geqslant  L^{i}\left(\boldsymbol{v}^{i} -\boldsymbol{u}^{i}\right)\\ 
& +L^{i-1}\left(\chi^{i-1}_{i-1}\gamma^{i}_{i-1} \left(\boldsymbol{v}^{i} -\boldsymbol{u}^{i}\right)\right) - a^{i-1}\left(\boldsymbol{u}^{i-1}, \chi^{i-1}_{i-1}\gamma^{i}_{i-1}\left(\boldsymbol{v}^{i}-\boldsymbol{u}^{i}\right)\right), \\
& a^{n}\left(\boldsymbol{u}^{n}, \boldsymbol{v}^{n}\right) = L^{n}\left(\boldsymbol{v}^{n}\right) + L^{n-1}\left(\chi^{n-1}_{n-1}\gamma^{n}_{n-1} \boldsymbol{v}^{n}\right) - a^{n-1}\left(\boldsymbol{u}^{n-1}, \chi^{n-1}_{n-1}\gamma^{n}_{n-1}\boldsymbol{v}^{n}\right), 
\end{aligned}
\right. 
\end{equation}
\end{small}
$\forall \boldsymbol{v}^{i} \in K^{i}(\boldsymbol{u}^{i+1}\cdot\boldsymbol{\alpha}^{i+1}),~i=1,\ldots,n-1$ and $\forall \boldsymbol{v}^{n} \in V^{n}$.
\end{proposition}

\begin{proof}
\textbf{($\Rightarrow$):}

Let $\boldsymbol{u} = (\boldsymbol{u}^1, \boldsymbol{u}^2, \ldots, \boldsymbol{u}^{n})$ solve Problem $P_{0}^{v}$.
From variational inequality (\ref{ieq:prb:3.2:var}), we can obtain:
\begin{small}
\begin{equation}\label{ieq:prb:3.2:var_1}
\sum_{i=1}^{n}a^{i}(\boldsymbol{u}^{i},\boldsymbol{v}^{i}-\boldsymbol{u}^{i}) + \sum_{i=1}^{n-1} \left(j^{i}(\boldsymbol{v}^{i},\boldsymbol{v}^{i+1}) - j^{i}(\boldsymbol{u}^{i},\boldsymbol{u}^{i+1})\right) \geqslant \sum_{i=1}^{n} L^{i}(\boldsymbol{v}^{i}-\boldsymbol{u}^{i}), ~\forall \boldsymbol{v} \in \mathcal{K}.
\end{equation}
\end{small}
According to the boundary conditions of Problem $P_{0}$, it can be known that 
\begin{small}
$$
\boldsymbol{u}^{i}\in K^{i}(\boldsymbol{u}^{i+1}\cdot \boldsymbol{\alpha}^{i+1}),~i=1,\ldots,n-1.
$$
\end{small}

(I) Let $\boldsymbol{v}\in \mathcal{K}$. According to the boundary conditions of $P_{0}$, we have:
\begin{small}
\begin{equation*}
a^{1}(\boldsymbol{u}^{1},\boldsymbol{v}^{1}-\boldsymbol{u}^{1}) + j^{1}(\boldsymbol{v}^{1},\boldsymbol{u}^{2}) - j^{1}(\boldsymbol{u}^{1},\boldsymbol{u}^{2}) \geqslant  L^{1}(\boldsymbol{v}^{1}-\boldsymbol{u}^{1}), ~\forall \boldsymbol{v}^{1} \in K^{1}(\boldsymbol{u}^{2}\cdot\boldsymbol{\alpha}^{2}).
\end{equation*}
\end{small}
Therefore, the first inequality of (\ref{ieq:prp:3.1:1}) is valied.

(II) Choose $\boldsymbol{v} = (\ldots,\boldsymbol{u}^{i-2}, \boldsymbol{u}^{i-1} + \chi^{i-1}_{i-1}\gamma^{i}_{i-1} \boldsymbol{w}^{i}, \boldsymbol{u}^{i} + \boldsymbol{w}^{i}, \boldsymbol{u}^{i+1}, \ldots) \in \mathcal{K}$, where $\boldsymbol{w}^{i} \in K^{i} \left(\left(\gamma^{i+1}_{i}\boldsymbol{u}^{i+1}-\gamma^{i}_{i}\boldsymbol{u}^{i}\right) \cdot \boldsymbol{\alpha}^{i+1}\right)$.
Then, we have
\begin{small}
\begin{equation}\label{proof:3.2.2}
\begin{aligned}
&a^{i-1}\left(\boldsymbol{u}^{i-1}, \chi^{i-1}_{i-1}\gamma^{i}_{i-1} \boldsymbol{w}^{i} \right) + a^{i}\left(\boldsymbol{u}^{i}, \boldsymbol{w}^{i}\right) + j^{i}\left( \boldsymbol{u}^{i} + \boldsymbol{w}^{i}, \boldsymbol{u}^{i+1}\right) - j^{i}\left( \boldsymbol{u}^{i}, \boldsymbol{u}^{i+1}\right) \\ 
& + j^{i-1}\left( \boldsymbol{u}^{i-1} + \chi^{i-1}_{i-1}\gamma^{i}_{i-1} \boldsymbol{w}^{i}, \boldsymbol{u}^{i} + \boldsymbol{w}^{i} \right) - j^{i-1}\left( \boldsymbol{u}^{i-1}, \boldsymbol{u}^{i}\right) \\
&+ j^{i-2}\left( \boldsymbol{u}^{i-2}, \boldsymbol{u}^{i-1} + \chi^{i-1}_{i-1}\gamma^{i}_{i-1} \boldsymbol{w}^{i} \right) - j^{i-2}\left( \boldsymbol{u}^{i-2}, \boldsymbol{u}^{i-1}\right)\\
\geqslant& L^{i-1}\left(\chi^{i-1}_{i-1}\gamma^{i}_{i-1} \boldsymbol{w}^{i}\right) + L^{i}\left(\boldsymbol{w}^{i}\right)
\end{aligned}.
\end{equation}
\end{small}
It can be noticed that
\begin{small}
\begin{equation*}
\begin{aligned}
&j^{i-1}\left( \boldsymbol{u}^{i-1} + \chi^{i}_{i-1}\gamma^{i}_{i-1} \boldsymbol{w}^{i}, \boldsymbol{u}^{i} + \boldsymbol{w}^{i} \right) =  j^{i-1}\left( \boldsymbol{u}^{i-1}, \boldsymbol{u}^i\right),\\
&j^{i-2}\left( \boldsymbol{u}^{i-2}, \boldsymbol{u}^{i-1} + \chi^{i-1}_{i-1}\gamma^{i}_{i-1} \boldsymbol{w}^{i} \right) = j^{i-2}\left( \boldsymbol{u}^{i-2}, \boldsymbol{u}^{i-1}\right),
\end{aligned}
\end{equation*}
\end{small}
where the second equality holds because $\left(\chi^{i-1}_{i-1}\gamma^{i}_{i-1} \boldsymbol{w}^{i}\right)_{\tau}=0$ on boundary $\Gamma_c^{i-2}$. Therefore, by letting $\boldsymbol{w}^i = \boldsymbol{v}^i - \boldsymbol{u}^i$ in (\ref{proof:3.2.2}), the second inequality of (\ref{ieq:prp:3.1:1}) holds $\forall \boldsymbol{v}^i \in K^i\left(\boldsymbol{u}^{i+1} \cdot \boldsymbol{\alpha}^{i+1}\right), i=2, \ldots, n-1$.

(III) Choose $\boldsymbol{v} = (\ldots,\boldsymbol{u}^{n-2}, \boldsymbol{u}^{n-1} \pm \chi^{n-1}_{n-1}\gamma^{n}_{n-1} \boldsymbol{w}^{n}, \boldsymbol{u}^{n} \pm \boldsymbol{w}^{n}) \in \mathcal{K}$. Then, we have
\begin{small}
\begin{equation}\label{proof:3.2.3}
\begin{aligned}
&a^{n-1}\left(\boldsymbol{u}^{n-1}, \pm\chi^{n-1}_{n-1}\gamma^{n}_{n-1} \boldsymbol{w}^{n} \right) + a^{n}\left(\boldsymbol{u}^{n}, \pm\boldsymbol{w}^{n}\right) \\
&+ j^{n-1}\left( \boldsymbol{u}^{n-1} \pm \chi^{n-1}_{n-1}\gamma^{n}_{n-1} \boldsymbol{w}^{n}, \boldsymbol{u}^{n} \pm \boldsymbol{w}^{n} \right) - j^{n-1}\left( \boldsymbol{u}^{n-1}, \boldsymbol{u}^{n}\right) \\
&+ j^{n-2}\left( \boldsymbol{u}^{n-2}, \boldsymbol{u}^{n-1} \pm \chi^{n-1}_{n-1}\gamma^{n}_{n-1} \boldsymbol{w}^{i} \right) - j^{n-2}\left( \boldsymbol{u}^{n-2}, \boldsymbol{u}^{n-1}\right)\\
\geqslant& L^{n-1}\left(\pm\chi^{n-1}_{n-1}\gamma^{n}_{n-1} \boldsymbol{w}^{n}\right) + L^{n}\left(\pm\boldsymbol{w}^{n}\right)
\end{aligned}
\end{equation}
\end{small}
Similar to the derivation process of (II), the last equation of (\ref{ieq:prp:3.1:1}) can be proved form the above (\ref{proof:3.2.3}).

\textbf{($\Leftarrow$):}

Let $\boldsymbol{u}^{i}\in K^{i}(\boldsymbol{u}^{i+1}\cdot\boldsymbol{\alpha}^{i+1})$, $i=1,\ldots,n-1$ and $\boldsymbol{u}^{n}\in V^{n}$ are the solution of (\ref{ieq:prp:3.1:1}) and $\boldsymbol{v}^i \in K^i(\boldsymbol{u}^{i+1} \cdot \boldsymbol{\alpha}^{i+1}), i=1, \ldots, n-1$ and $\boldsymbol{v}^n \in V^n$. Then, we have
\begin{small}
$$
\boldsymbol{w}^i\cdot\boldsymbol{\beta}^i + \boldsymbol{w}^{i+1} \cdot\boldsymbol{\alpha}^{i+1} \leqslant 0 \Rightarrow \boldsymbol{w}^i\cdot\boldsymbol{\beta}^i + \boldsymbol{w}^{i+1} \cdot\boldsymbol{\alpha}^{i+1} - \boldsymbol{u}^{i+1} \cdot\boldsymbol{\alpha}^{i+1} \leqslant - \boldsymbol{u}^{i+1} \cdot\boldsymbol{\alpha}^{i+1},
$$
\end{small}
where $\boldsymbol{w}^i\in K^{i}(\boldsymbol{w}^{i+1} \cdot\boldsymbol{\alpha}^{i+1})$, $i=1,\ldots,n-1$. By choosing $\boldsymbol{v}^i = \boldsymbol{w}^i + (\chi^{i}_{i}\gamma^{i+1}_{i}\boldsymbol{u}^{i+1} - \boldsymbol{w}^{i+1}) \in K^{i}(\boldsymbol{u}^{i+1} \cdot\boldsymbol{\alpha}^{i+1})$, $i=1,\ldots,n-1$ in (\ref{ieq:prp:3.1:1}), it can be obtained that
\begin{footnotesize}
\begin{equation}\label{proof:3.2.4}
\left\{
\begin{aligned}
&a^1\left(\boldsymbol{u}^1, \boldsymbol{w}^1 + \chi^{1}_{1}\gamma^{2}_{1}\left(\boldsymbol{u}^{2} - \boldsymbol{w}^{2}\right) -\boldsymbol{u}^1\right) +j^1\left( \boldsymbol{w}^1 + \chi^{1}_{1}\gamma^{2}_{1}\left(\boldsymbol{u}^{2} - \boldsymbol{w}^{2}\right), \boldsymbol{u}^2\right) -j^{1}\left( \boldsymbol{u}^1, \boldsymbol{u}^2\right)\\ 
&\geqslant L^1\left(\boldsymbol{w}^i + \chi^{1}_{1}\gamma^{2}_{1}\left(\boldsymbol{u}^{2} - \boldsymbol{w}^{2}\right)-\boldsymbol{u}^1\right),\\
& a^i\left(\boldsymbol{u}^i, \boldsymbol{w}^i + \chi^{i}_{i}\gamma^{i+1}_{i}\left(\boldsymbol{u}^{i+1} - \boldsymbol{w}^{i+1}\right) -\boldsymbol{u}^i\right) +j^i\left( \boldsymbol{w}^i + \chi^{i}_{i}\gamma^{i+1}_{i}\left(\boldsymbol{u}^{i+1} - \boldsymbol{w}^{i+1}\right), \boldsymbol{u}^{i+1}\right) \\
& -j^i\left( \boldsymbol{u}^i, \boldsymbol{u}^{i+1}\right) \geqslant L^i\left(\boldsymbol{w}^i + \chi^{i}_{i}\gamma^{i+1}_{i}\left(\boldsymbol{u}^{i+1} - \boldsymbol{w}^{i+1}\right) -\boldsymbol{u}^i\right) +L^{i-1}\Big(\chi_{i-1}^{i-1} \gamma^i_{i-1}\Big(\boldsymbol{w}^i + \\
& \chi^{i}_{i}\gamma^{i+1}_{i}\left(\boldsymbol{u}^{i+1} - \boldsymbol{w}^{i+1}\right) -\boldsymbol{u}^i\Big)\Big) -a^{i-1}\left(\boldsymbol{u}^{i-1}, \chi^{i-1}_{i-1} \gamma_{i-1}^i\left(\boldsymbol{w}^i + \chi^{i}_{i}\gamma^{i+1}_{i}\left(\boldsymbol{u}^{i+1} - \boldsymbol{w}^{i+1}\right) -\boldsymbol{u}^i\right)\right).
\end{aligned}\right.
\end{equation}
\end{footnotesize}
It is worth noting that:
\begin{footnotesize}
\begin{equation*}
\begin{aligned}
&j^i\left( \boldsymbol{w}^i +\chi_i^{i} \gamma_i^{i+1}\left(\boldsymbol{u}^{i+1}-\boldsymbol{w}^{i+1}\right), \boldsymbol{u}^{i+1}\right) = j^i\left( \boldsymbol{w}^i, \boldsymbol{w}^{i+1}\right),\\
& L^{i-1}\left(\chi_{i-1}^{i-1} \gamma_{i-1}^i\left(\boldsymbol{w}^i+\chi_i^{i} \gamma_i^{i+1}\left(\boldsymbol{u}^{i+1}-\boldsymbol{w}^{i+1}\right)-\boldsymbol{u}^i\right)\right)
= L^{i-1}\left(\chi_{i-1}^{i-1} \gamma_{i-1}^i\left(\boldsymbol{w}^i -\boldsymbol{u}^i\right)\right);\\
&a^{i-1}\left(\boldsymbol{u}^{i-1}, \chi_{i-1}^{i-1} \gamma_{i-1}^i\left(\boldsymbol{w}^i + \chi^{i}_{i}\gamma^{i+1}_{i}\left(\boldsymbol{u}^{i+1} - \boldsymbol{w}^{i+1}\right) -\boldsymbol{u}^i\right)\right)
= a^{i-1}\left(\boldsymbol{u}^{i-1}, \chi_{i-1}^{i-1} \gamma_{i-1}^i\left(\boldsymbol{w}^i  -\boldsymbol{u}^i\right)\right).
\end{aligned}
\end{equation*}
\end{footnotesize}
where $\chi_{i-1}^{i-1} \gamma_{i-1}^i \chi^{i}_{i}\gamma^{i+1}_{i}\left(\boldsymbol{u}^{i+1} - \boldsymbol{w}^{i+1}\right) = 0$ on $\Omega^{i-1}$. Therefore, (\ref{proof:3.2.4}) can be rewritten as:
\begin{small}
\begin{equation}\label{proof:3.2.5}
\left\{
\begin{aligned}
&a^1\left(\boldsymbol{u}^1, \boldsymbol{w}^1 -\boldsymbol{u}^1\right) + a^1\left(\boldsymbol{u}^1, \chi^{1}_{1}\gamma^{2}_{1}\left(\boldsymbol{u}^{2} - \boldsymbol{w}^{2}\right)\right)
+j^1\left( \boldsymbol{w}^1, \boldsymbol{w}^2\right) -j^{1}\left( \boldsymbol{u}^1, \boldsymbol{u}^2\right)\\ 
&\geqslant L^1\left(\boldsymbol{w}^i-\boldsymbol{u}^1\right) + L^1\left( \chi^{1}_{1}\gamma^{2}_{1}\left(\boldsymbol{u}^{2} - \boldsymbol{w}^{2}\right)\right),\\
& a^i\left(\boldsymbol{u}^i, \boldsymbol{w}^i -\boldsymbol{u}^i\right) 
+ a^i\left(\boldsymbol{u}^i, \chi^{i}_{i}\gamma^{i+1}_{i}\left(\boldsymbol{u}^{i+1} - \boldsymbol{w}^{i+1}\right) \right)
+j^i\left( \boldsymbol{w}^i, \boldsymbol{w}^{i+1}\right) 
-j^i\left( \boldsymbol{u}^i, \boldsymbol{u}^{i+1}\right) \\
& \geqslant L^i\left(\boldsymbol{w}^i -\boldsymbol{u}^i\right) 
+ L^i\left(\chi^{i}_{i}\gamma^{i+1}_{i}\left(\boldsymbol{u}^{i+1} - \boldsymbol{w}^{i+1}\right) \right)
+L^{i-1}\left(\chi_{i-1}^{i-1} \gamma_{i-1}^i\left(\boldsymbol{w}^i  -\boldsymbol{u}^i\right)\right) \\
&-a^{i-1}\left(\boldsymbol{u}^{i-1}, \chi_{i-1}^{i-1} \gamma_{i-1}^i\left(\boldsymbol{w}^i  -\boldsymbol{u}^i\right)\right),~~i=2,\ldots,n-1.
\end{aligned}\right.
\end{equation}
\end{small}
Furthermore, by choosing $\boldsymbol{v}^n = \boldsymbol{w}^n - \boldsymbol{u}^n$, we have:
\begin{small}
\begin{equation}\label{proof:3.2.6}
\begin{aligned}
a^n\left(\boldsymbol{u}^n, \boldsymbol{w}^n - \boldsymbol{u}^n\right) =& L^n\left(\boldsymbol{w}^n - \boldsymbol{u}^n\right) +L^{n-1}\left(\chi_{n-1}^{n-1} \gamma_{n-1}^n \left(\boldsymbol{w}^n - \boldsymbol{u}^n \right)\right)\\ 
&-a^{n-1}\left(\boldsymbol{u}^{n-1}, \chi_{n-1}^{n-1} \gamma_{n-1}^n \left(\boldsymbol{w}^n - \boldsymbol{u}^n \right)\right)
\end{aligned}
\end{equation}
\end{small}
Add formula (\ref{proof:3.2.5}) and formula (\ref{proof:3.2.6}) to get:
\begin{small}
\begin{equation*}
\sum_{i=1}^{n} a^i\left(\boldsymbol{u}^i, \boldsymbol{w}^i-\boldsymbol{u}^i\right)
+ \sum_{i=1}^{n-1}\left(j^i\left( \boldsymbol{w}^i, \boldsymbol{w}^{i+1}\right)-j^i\left( \boldsymbol{u}^i, \boldsymbol{u}^{i+1}\right)\right)  \geqslant \sum_{i=1}^{n} L^i\left(\boldsymbol{w}^i-\boldsymbol{u}^i\right),
\end{equation*}
\end{small}
where $\boldsymbol{w}=(\boldsymbol{w}^1,\ldots,\boldsymbol{w}^{n}) \in \mathcal{K}$, so it is equivalent to variational inequality (\ref{ieq:prb:3.2:var}).
\end{proof}

Proposition \ref{prp:3.1} proves that the original variational inequality (\ref{ieq:prb:3.2:var}) can be decomposed into $n$ coupled variational inequalities. In the next section, an algorithm to solve this coupling problem is presented.

\section{Layer decomposition algorithm for \texorpdfstring{$P_0^v$}{}}

Based on Proposition \ref{prp:3.1}, the decomposition algorithm for solving  Problem $P_0^v$ will be proposed and analyzed. 
The displacement parameter $\boldsymbol{\lambda}=\left(\lambda^{1},\ldots,\lambda^{n-1}\right)$ on the friction boundary $\Gamma^{i}_{c}$ is specified as: $\lambda^{i}\in H^{1/2}\left(\Gamma_c^i\right)$, $i=1,\ldots,n-1$. 
To simplify the formula, let's specify $\lambda^{0}=0$ and $\lambda^{n}=0$.
Let $\boldsymbol{u}=\left(\boldsymbol{u}^1, \boldsymbol{u}^2, \ldots, \boldsymbol{u}^n\right)$ be the solutions of the following decoupled problems: 
\begin{problem}[$P_d^i(\lambda)$]\label{prb:3.4:pd}
\begin{small}
\begin{align}
&\left\{
\begin{aligned}
&\text{Find } \boldsymbol{u}^1:=\boldsymbol{u}^1(\boldsymbol{\lambda}) \in K^{1}(\lambda^{1}\cdot\boldsymbol{\alpha}^{2}) \text{ such that:}\\
&a^1\left(\boldsymbol{u}^1, \boldsymbol{v}^1-\boldsymbol{u}^1\right)+j^1\left( \boldsymbol{v}^1, \lambda^1\right)-j^1\left( \boldsymbol{u}^1, \lambda^1\right) \geqslant L^1\left(\boldsymbol{v}^1-\boldsymbol{u}^1\right)\\
& \forall \boldsymbol{v}^1 \in K^{1}(\lambda^{1}\cdot\boldsymbol{\alpha}^{2})
\end{aligned}
\right.  \tag{\text{$P^{1}_{d}(\boldsymbol{\lambda})$}} \label{prb:3.4:d1}\\
&\left\{
\begin{aligned}
&\text{Find } \boldsymbol{u}^i:=\boldsymbol{u}^i(\boldsymbol{\lambda}) \in K^{i}(\lambda^{i}\cdot\boldsymbol{\alpha}^{i+1}), ~i=2,\ldots,n-1 \text{ such that:}\\
&a^i\left(\boldsymbol{u}^i, \boldsymbol{v}^i-\boldsymbol{u}^i\right)+j^i\left( \boldsymbol{v}^i, \lambda^i\right)-j^i\left( \boldsymbol{u}^i, \lambda^i\right) \geqslant L^i\left(\boldsymbol{v}^i-\boldsymbol{u}^i\right)\\
& \forall \boldsymbol{v}^i \in K^{i}(\lambda^{i}\cdot\boldsymbol{\alpha}^{i+1}) \cap V^{i}_{2}(\boldsymbol{u}^i), ~\boldsymbol{u}^i = \lambda^{i-1} \text{ on } \Gamma^{i-1}_{c}
\end{aligned}
\right.  \tag{\text{$P^{i}_{d}(\boldsymbol{\lambda})$}} \label{prb:3.4:di}\\
&\left\{
\begin{aligned}
&\text{Find } \boldsymbol{u}^n:=\boldsymbol{u}^n(\boldsymbol{\lambda}) \in V^{n} \text{ such that:}\\
&a^n\left(\boldsymbol{u}^n, \boldsymbol{v}^n\right)=L^n\left(\boldsymbol{v}^n\right)\\
& \forall \boldsymbol{v}^n \in V^{n}_{2}(0), ~ \boldsymbol{u}^n = \lambda^{n-1} \text{ on } \Gamma^{n-1}_{c}
\end{aligned}
\right.  \tag{\text{$P^{n}_{d}(\boldsymbol{\lambda})$}}\label{prb:3.4:dn}
\end{align}
\end{small}
\end{problem}
Therefore, if $\boldsymbol{\lambda}$ can be chosen in such way that $\sigma_\alpha^{i+1}(\boldsymbol{u}^{i+1})= -\sigma_\beta^{i}(\boldsymbol{u}^{i})$ and $\boldsymbol{\sigma}_\tau^{i+1}( \boldsymbol{u}^{i+1} )$ $= \boldsymbol{\sigma}_\eta^{i}( \boldsymbol{u}^{i} )$ on $\Gamma_{c}^{i}, ~i=1,\ldots,n-1$.
Then the $\boldsymbol{u}=\left(\boldsymbol{u}^1, \boldsymbol{u}^2, \ldots, \boldsymbol{u}^n\right)$ solving the above problems would be a solution of Problem $P_0^v$.
In order to find such $\boldsymbol{\lambda}$ which ensures continuity of the normal and tangential contact stress across $\Gamma_c^{i}$, the following auxiliary Neumann problems defined in $\Omega^i$ should be defined:
\begin{problem}[$P_p^{i+1}(\boldsymbol{\lambda})$ and $P_q^{i}(\boldsymbol{\lambda})$]\label{prb:3.5:P_pq}
Find $\boldsymbol{p}^{i+1}\in V^{i+1}_{3}(0)$, $\boldsymbol{q}^{i} \in V^{i}_{2}(0) (i=1,\ldots,n-1)$  such that:
\begin{small}
\begin{align}
a^{i+1}\left(\boldsymbol{p}^{i+1}, \boldsymbol{v}^{i+1}\right) =& \frac{1}{2} \Big(a^{i+1}\left(\boldsymbol{u}^{i+1}, \boldsymbol{v}^{i+1}\right)- L^{i+1}\left(\boldsymbol{v}^{i+1}\right)+ a^{i}\left(\boldsymbol{u}^{i}, \chi_{i}^{i} \gamma^{i+1}_{i}\boldsymbol{v}^{i+1}\right) \nonumber\\ 
&- L^{i}\left(\chi_{i}^{i} \gamma^{i+1}_{i}\boldsymbol{v}^{i+1}\right) \Big),~ \forall \boldsymbol{v}^{i+1} \in V^{i+1}_{3}(0) 
\tag{\text{$P^{i+1}_{p}(\boldsymbol{\lambda})$}} 
\label{prb:3.5:pi} \\
a^{i}\left(\boldsymbol{q}^{i}, \boldsymbol{v}^{i}\right) =& \frac{1}{2} \Big(a^{i+1}\left(\boldsymbol{u}^{i+1}, \chi^{i}_{i+1} \gamma^{i}_{i}\boldsymbol{v}^{i}\right)- L^{i+1}\left(\chi^{i}_{i+1} \gamma^{i}_{i} \boldsymbol{v}^{i}\right) \nonumber\\
&+ a^{i}\left(\boldsymbol{u}^{i},  \boldsymbol{v}^{i}\right) - L^{i}\left(\boldsymbol{v}^{i}\right) \Big),~ \forall \boldsymbol{v}^{i} \in V^{i}_{2}(0) 
\tag{\text{$P^{i}_{q}(\boldsymbol{\lambda})$}} 
\label{prb:3.5:qi}
\end{align}
\end{small}
where $\boldsymbol{u}^{i}:= \boldsymbol{u}^{i}(\boldsymbol{\lambda})$ are the solutions of $P_{d}^{i}(\boldsymbol{\lambda})$.
\end{problem}

\begin{algorithm}[!ht]
\caption{Layer Decomposition Algorithm}
\label{algorithm:3.1}
Let $\boldsymbol{\lambda}_{0} = (\lambda^{1}_{0},\ldots,\lambda^{n-1}_{0}) \in H^{1/2}(\Gamma^{i}_{c}) \times \cdots \times H^{1/2}(\Gamma^{n-1}_{c})$ and $\theta>0$ be given.

\textbf{Output:} $\boldsymbol{u}^{i}_{k}$, $\boldsymbol{p}^{i}_{k}$, $\boldsymbol{q}^{i}_{k}$, $\boldsymbol{\lambda}_{k}$

\textbf{For} $k\leq N$ \textbf{do}

~~~~$
\left\{
\begin{aligned}
& \boldsymbol{u}^{1}_{k}(\boldsymbol{\lambda}_{k-1}) \in K^1\left(\lambda^1_{k-1} \cdot \boldsymbol{\alpha}^2\right) \text{ solves } P^{1}_{d}(\boldsymbol{\lambda}_{k-1})\\
& \boldsymbol{u}^{i}_{k}(\boldsymbol{\lambda}_{k-1}) \in K^i\left(\lambda^i_{k-1} \cdot \boldsymbol{\alpha}^{i+1} \right) \text{ solves } P^{i}_{d}(\boldsymbol{\lambda}_{k-1}), ~i=2,\ldots,n-1\\
& \boldsymbol{u}^{n}_{k}(\boldsymbol{\lambda}_{k-1}) \in V^{n} \text{ solves } P^{n}_{d}(\boldsymbol{\lambda}_{k-1})
\end{aligned}
\right.
$

~~~~$
\left\{
\begin{aligned}
& \boldsymbol{p}^{i+1}_{k} \in V^{i+1} \text{ solves } P^{i+1}_{p}(\boldsymbol{\lambda}_{k-1})\\
& \boldsymbol{q}^{i}_{k} \in V^{i} \text{ solves } P^{i}_{q}(\boldsymbol{\lambda}_{k-1})
\end{aligned}
\right.
$

~~~~$\lambda_k^i=\lambda_{k-1}^{i} -\theta\left(\gamma^{i+1}_{i}\boldsymbol{p}^{i+1}_{k} + \gamma^{i}_{i}\boldsymbol{q}^{i}_{k}\right)$ on $\Gamma_c^i$.

\textbf{end}
\end{algorithm}

Therefore, the algorithm can be illustrated as following Algorithm \ref{algorithm:3.1}.
The convergence verification of the algorithm \ref{algorithm:3.1} is the core of this section.
\begin{theorem}\label{thm:3.1}
There exist: $0<\theta_*<4/(3c_{t}^{\min})^2$ and functions $\boldsymbol{\lambda}_{*} = (\lambda^{1}_{*}, \ldots,$ $\lambda^{n-1}_{*}) \in H^{1/2}(\Gamma^{i}_{c}) \times \cdots \times H^{1/2}(\Gamma^{n-1}_{c})$, $\boldsymbol{u}^{i}_*$, $\boldsymbol{p}^{i}_*$, $\boldsymbol{q}^{i}_*$, $i=1,\ldots,n$ such that $\forall\theta \in\left(0, \theta_*\right)$,
\begin{small}
\begin{align*}
&~\lambda_k^i \rightarrow \lambda_*^i \text { in } H^{1 / 2}\left(\Gamma_c^i\right), && i=1, \ldots, n-1\\
&\left.\begin{array}{l}
\boldsymbol{p}_k^{i+1} \rightarrow \boldsymbol{p}_*^{i+1} \text { in }\left(H^{1}\left(\Omega^{i+1}\right)\right)^d\\
\boldsymbol{q}_k^{i} \rightarrow \boldsymbol{p}_*^i \text { in }\left(H^{1}\left(\Omega^{i}\right)\right)^d
\end{array}\right\}, && i=1, \ldots, n-1 \\
&~\boldsymbol{u}_k^i \rightarrow \boldsymbol{u}_*^i \text { in } \left(H^{1}\left(\Omega^i\right)\right)^d, && i=1 \ldots, n
\end{align*}
\end{small}
as $k\to \infty$, where the sequence $\left\{\left(\boldsymbol{u}_k^i, \boldsymbol{p}_k^i, \boldsymbol{q}_k^i, \lambda_k^i\right)\right\}$ is generated by Algorithm \ref{algorithm:3.1}. In addition, the $\boldsymbol{u}_{*}=\left(\boldsymbol{u}^1_{*}, \ldots, \boldsymbol{u}^n_{*} \right) \in \mathcal{K}$ solves $P_1$.
\end{theorem}

In order to prove the Theorem \ref{thm:3.1}, some Propositions about the Problem \ref{prb:3.4:pd}, Problem \ref{prb:3.5:P_pq} need to be verified.

\begin{proposition}\label{prp:3.2}
Let $\boldsymbol{u}^{i}$ and $\tilde{\boldsymbol{u}}^{i}$ be the solutions of $P^i_d(\boldsymbol{\lambda}), P^i_d(\tilde{\boldsymbol{\lambda}})$ with $\boldsymbol{\lambda}, \tilde{\boldsymbol{\lambda}} \in H^{1 / 2}\left(\Gamma_c^i\right) \times \cdots \times H^{1 / 2}\left(\Gamma_c^{n-1}\right)$, $i=1,\ldots,n$, respectively. Then
\begin{small}
\begin{align}
& a^1\left(\boldsymbol{u}^1-\tilde{\boldsymbol{u}}^1, \boldsymbol{u}^1-\tilde{\boldsymbol{u}}^1\right) \leqslant a^1\left(\boldsymbol{u}^1-\tilde{\boldsymbol{u}}^1, \chi_1^1\left(\lambda^1-\tilde{\lambda}^1\right)\right) \label{ieq:prp:3.2:1} \\
& a^i\left(\boldsymbol{u}^i-\tilde{\boldsymbol{u}}^i, \boldsymbol{u}^i-\tilde{\boldsymbol{u}}^i\right) \leqslant a^i\left(\boldsymbol{u}^i-\tilde{\boldsymbol{u}}^i, \chi_i^i\left(\lambda^i-\tilde{\lambda}^i\right) + \chi_i^{i-1}\left(\lambda^{i-1}-\tilde{\lambda}^{i-1}\right)\right) \label{ieq:prp:3.2:2}\\
& \left\|\boldsymbol{u}^n-\tilde{\boldsymbol{u}}^n\right\|_{1, \Omega^n} = \left\| \lambda^{n-1}-\tilde{\lambda}^{n-1} \right\|_{\Gamma_c^{n-1,n}} \label{eq:prp:3.2:3}
\end{align}
\end{small}
\end{proposition}

\begin{proof}
(1) Let $\tilde{\boldsymbol{v}}^{1} = \boldsymbol{u}^1 - \chi_1^1 \left(\lambda^1-\tilde{\lambda}^1\right)$, then $\tilde{\boldsymbol{v}}^1 \in K^1\left(\tilde{\lambda}^1 \cdot \boldsymbol{\alpha}^{2}\right)$. Indeed, on $\Gamma^1_c$
\begin{small}
$$
\tilde{\boldsymbol{v}}^1 \cdot \boldsymbol{\beta}^1=\boldsymbol{u}^1 \cdot \boldsymbol{\beta}^1 -\lambda^1 \cdot \boldsymbol{\beta}^1 +\tilde{\lambda}^1 \cdot \boldsymbol{\beta}^1 = \boldsymbol{u}^1 \cdot \boldsymbol{\beta}^1 +\lambda^1 \cdot \boldsymbol{\alpha}^2 -\tilde{\lambda}^1 \cdot \boldsymbol{\alpha}^2 \leqslant -\tilde{\lambda}^1 \cdot \boldsymbol{\alpha}^2
$$
\end{small}
where $\boldsymbol{u}^1 \in K^1\left(\lambda^1 \cdot \boldsymbol{\alpha}^2 \right)$. Hence
\begin{small}
\begin{align*}
& a^1\left(\tilde{\boldsymbol{u}}^1, \boldsymbol{u}^1 - \chi_1^1 \left(\lambda^1-\tilde{\lambda}^1\right) -\tilde{\boldsymbol{u}}^1\right) +j^1\left( \boldsymbol{u}^1 - \chi_1^1 \left(\lambda^1-\tilde{\lambda}^1\right), \tilde{\lambda}^1\right) \\
&-j^1\left( \tilde{\boldsymbol{u}}^1, \tilde{\lambda}^1\right) \geqslant  L^1\left(\boldsymbol{u}^1 - \chi_1^1 \left(\lambda^1-\tilde{\lambda}^1\right) - \tilde{\boldsymbol{u}}^1\right).
\end{align*}
\end{small}
Similarly, let ${\boldsymbol{v}}^{1} = \tilde{\boldsymbol{u}}^1 - \chi_1^1 \left(\tilde{\lambda}^1-{\lambda}^1\right)$, then ${\boldsymbol{v}}^1 \in K^1\left({\lambda}^1 \cdot \boldsymbol{\alpha}^{2}\right)$. Hence
\begin{small}
\begin{align*}
& a^1\left({\boldsymbol{u}}^1, \tilde{\boldsymbol{u}}^1 - \chi_1^1 \left(\tilde{\lambda}^1-{\lambda}^1\right) -{\boldsymbol{u}}^1\right) +j^1\left( \tilde{\boldsymbol{u}}^1 - \chi_1^1 \left(\tilde{\lambda}^1- {\lambda}^1\right), {\lambda}^1\right) \\
&-j^1\left( {\boldsymbol{u}}^1, {\lambda}^1\right) \geqslant L^1\left(\tilde{\boldsymbol{u}}^1 - \chi_1^1 \left(\tilde{\lambda}^1-{\lambda}^1\right) - {\boldsymbol{u}}^1\right).
\end{align*}
\end{small}
It should be noted that
\begin{small}
\begin{align*}
j^1\left( \boldsymbol{u}^1-\chi_1^1\left(\lambda^1-\tilde{\lambda}^1\right), \tilde{\lambda}^1\right)
& = j^{1} \left( \boldsymbol{u}^1, {\lambda}^1\right)\\
j^1\left( \tilde{\boldsymbol{u}}^1 -\chi_1^1 \left(\tilde{\lambda}^1-\lambda^1\right), \lambda^1\right) & = j^1\left( \tilde{\boldsymbol{u}}^1, \tilde{\lambda}^1\right)
\end{align*}
\end{small}
Therefore, by summing above two inequalities, we have:
\begin{small}
$$
a^1\left(\boldsymbol{u}^1-\tilde{\boldsymbol{u}}^1, \boldsymbol{u}^1-\tilde{\boldsymbol{u}}^1\right) \leqslant a^1\left(\boldsymbol{u}^1-\tilde{\boldsymbol{u}}^1, \chi_1^1\left(\lambda^1-\tilde{\lambda}^1\right)\right)
$$
\end{small}

(2) Let $\tilde{\boldsymbol{v}}^{i} = \boldsymbol{u}^i - \chi_i^i \left(\lambda^i-\tilde{\lambda}^i\right) - \chi^{i-1}_{i} \left(\lambda^{i-1}-\tilde{\lambda}^{i-1}\right)$, then $\tilde{\boldsymbol{v}}^{i} \in K^i\left(\tilde{\lambda}^i \cdot \boldsymbol{\alpha}^{i+1}\right) \cap V_2^i\left(\tilde{u}^i\right)$. Indeed, on $\Gamma_{c}^{i}$
\begin{small}
$$
\tilde{\boldsymbol{v}}^i \cdot \boldsymbol{\beta}^i=\boldsymbol{u}^i \cdot \boldsymbol{\beta}^{i} -\lambda^i \cdot \boldsymbol{\beta}^i +\tilde{\lambda}^i \cdot \boldsymbol{\beta}^i = \boldsymbol{u}^i \cdot \boldsymbol{\beta}^i +\lambda^i \cdot \boldsymbol{\alpha}^{i+1} -\tilde{\lambda}^i \cdot \boldsymbol{\alpha}^{i+1} \leqslant -\tilde{\lambda}^i \cdot \boldsymbol{\alpha}^{i+1}
$$
\end{small}
and on $\Gamma_c^{i-1}$
\begin{small}
$$
\tilde{\boldsymbol{v}}^i=\boldsymbol{u}^i-\chi_i^{i-1}\left(\lambda^{i-1}-\tilde{\lambda}^{i-1}\right) = \chi_i^{i-1} \tilde{\lambda}^{i-1} = \tilde{\boldsymbol{u}}^{i-1}.
$$
\end{small}
Hence
\begin{small}
\begin{align*}
& a^i\left(\tilde{\boldsymbol{u}}^i, \boldsymbol{u}^i - \chi_i^i \left(\lambda^i-\tilde{\lambda}^i\right) - \chi^{i-1}_{i} \left(\lambda^{i-1}-\tilde{\lambda}^{i-1}\right) - \tilde{\boldsymbol{u}}^i\right) \\
&+j^i\left( \boldsymbol{u}^i - \chi_i^i \left(\lambda^i-\tilde{\lambda}^i\right) - \chi^{i-1}_{i} \left(\lambda^{i-1}-\tilde{\lambda}^{i-1}\right), \tilde{\lambda}^i\right) -j^i\left( \tilde{\boldsymbol{u}}^i, \tilde{\lambda}^i\right)\\
\geqslant& L^i\left(\boldsymbol{u}^i - \chi_i^i \left(\lambda^i-\tilde{\lambda}^i\right) - \chi^{i-1}_{i} \left(\lambda^{i-1}-\tilde{\lambda}^{i-1}\right) - \tilde{\boldsymbol{u}}^i\right).
\end{align*}
\end{small}
Similarly, let ${\boldsymbol{v}}^{i} = \tilde{\boldsymbol{u}}^i - \chi_i^i \left(\tilde{\lambda}^i-{\lambda}^i\right) - \chi^{i-1}_{i} \left(\tilde{\lambda}^{i-1}-{\lambda}^{i-1}\right)$, then ${\boldsymbol{v}}^{i} \in K^i\left({\lambda}^i \cdot \boldsymbol{\alpha}^{i+1}\right) \cap V_2^i\left({u}^i\right)$.
Hence:
\begin{small}
\begin{align*}
& a^i\left({\boldsymbol{u}}^i, \tilde{\boldsymbol{u}}^i - \chi_i^i \left(\tilde{\lambda}^i-{\lambda}^i\right) - \chi^{i-1}_{i} \left(\tilde{\lambda}^{i-1}-{\lambda}^{i-1}\right) - {\boldsymbol{u}}^i\right)\\ 
&+j^i\left( \tilde{\boldsymbol{u}}^i - \chi_i^i \left(\tilde{\lambda}^i- {\lambda}^i\right) - \chi^{i-1}_{i} \left(\tilde{\lambda}^{i-1} -{\lambda}^{i-1}\right), {\lambda}^i\right)
-j^i\left( {\boldsymbol{u}}^i, {\lambda}^i\right)\\
\geqslant& L^i\left(\tilde{\boldsymbol{u}}^i - \chi_i^i \left(\tilde{\lambda}^i-{\lambda}^i\right) - \chi^{i-1}_{i} \left(\tilde{\lambda}^{i-1}-{\lambda}^{i-1}\right) -{\boldsymbol{u}}^i\right).
\end{align*}
\end{small}
It can be noted that
\begin{small}
\begin{align*}
j^i\left( \boldsymbol{u}^i - \chi_i^i \left(\lambda^i-\tilde{\lambda}^i\right) - \chi^{i-1}_{i} \left(\lambda^{i-1}-\tilde{\lambda}^{i-1}\right), \tilde{\lambda}^i\right)= j^{i} \left( \boldsymbol{u}^i, {\lambda}^i\right),
\end{align*}
\end{small}
where $\chi_i^{i-1}\left(\lambda^{i-1}-\tilde{\lambda}^{i-1}\right)=0$ on $\Gamma_c^i$.
Then, it can be verified that
\begin{small}
$$
j^i\left( \tilde{\boldsymbol{u}}^i - \chi_i^i \left(\tilde{\lambda}^i-{\lambda}^i\right) - \chi^{i-1}_{i} \left(\tilde{\lambda}^{i-1}-{\lambda}^{i-1}\right), \lambda^i\right)  = j^i\left( \tilde{\boldsymbol{u}}^i, \tilde{\lambda}^i\right).
$$
\end{small}
Therefore, by summing above two inequalities, we have:
\begin{small}
$$
a^i\left(\boldsymbol{u}^i-\tilde{\boldsymbol{u}}^i, \boldsymbol{u}^i-\tilde{\boldsymbol{u}}^i\right) \leqslant a^i\left(\boldsymbol{u}^i-\tilde{\boldsymbol{u}}^i, \chi_i^i\left(\lambda^i-\tilde{\lambda}^i\right)+\chi_i^{i-1}\left(\lambda^{i-1}-\tilde{\lambda}^{i-1}\right)\right).
$$
\end{small}

(3) The function $\boldsymbol{u}^n-\tilde{\boldsymbol{u}}^n$ satisfies
\begin{small}
$$
\left.\begin{array}{l}
a^n\left(\boldsymbol{u}^n, \boldsymbol{v}^n\right)=L\left(\boldsymbol{v}^n\right) \\
a^n\left(\tilde{\boldsymbol{u}}^n, \boldsymbol{v}^n\right)=L\left(\boldsymbol{v}^n\right)
\end{array}\right\rangle a^n\left(\boldsymbol{u}^n-\tilde{\boldsymbol{u}}^n, \boldsymbol{v}^n\right)=0,~~ \forall  \boldsymbol{v}^n \in V_2^n(0)
$$
\end{small}
and
\begin{small}
$$
\boldsymbol{u}^n-\tilde{\boldsymbol{u}}^n=\lambda^{n-1}-\tilde{\lambda}^{n-1} \text { on } \Gamma_c^{n-1}
$$
\end{small}
Therefore, based on the definition of $\|\cdot\|_{\Gamma_c^{n-1, n}}$, the equality is verified.
\end{proof}

\begin{consequence}\label{con:3.1}
There exists a constant $c^{i}_{2}>0$ which is independent of $\lambda^{i}$ and $\tilde{\lambda}^i \in H^{1/2}\left(\Gamma_c^i\right)$ such that:
\begin{small}
\begin{align}
\left\|\boldsymbol{u}^1-\tilde{\boldsymbol{u}}^1\right\|_{1, \Omega^1} &\leqslant \left\|\lambda^{1}-\tilde{\lambda}^{1}\right\|_{\Gamma_c^{1,1}} \leqslant c^{1}_{2} \left\|\lambda^{1}-\tilde{\lambda}^{1}\right\|_{\Gamma_c^{1,2}} \label{ieq:con:3.1:1}\\
\left\|\boldsymbol{u}^i-\tilde{\boldsymbol{u}}^i\right\|_{1, \Omega^i} &\leqslant c^{i}_{2}\left\|\lambda^{i}-\tilde{\lambda}^{i}\right\|_{\Gamma_c^{i,i+1}} + \left\|\lambda^{i-1}-\tilde{\lambda}^{i-1}\right\|_{\Gamma_c^{i-1,i}} \label{ieq:con:3.1:2}
\end{align}
\end{small}
where $i=1,\ldots,n-1$.
\end{consequence}
\begin{proof}
From proposition \ref{prp:3.2}, it follows that 
\begin{small}
$$
\left\|\boldsymbol{u}^1-\tilde{\boldsymbol{u}}^1\right\|_{1, \Omega^1}^{2} \leqslant \left\|\boldsymbol{u}^1-\tilde{\boldsymbol{u}}^1\right\|_{1, \Omega^1} \cdot \left\|\chi_1^1\left(\lambda^1-\tilde{\lambda}^1\right) \right\| = \left\|\boldsymbol{u}^1-\tilde{\boldsymbol{u}}^1\right\|_{1, \Omega^1} \cdot \left\|\lambda^1-\tilde{\lambda}^1\right\|_{\Gamma_c^{1,1}}
$$
\end{small}
The last inequality of (\ref{ieq:con:3.1:1}) follows from the fact that the norms $\left\|\cdot\right\|_{\Gamma_c^{1,1}}$ and $\left\|\cdot\right\|_{\Gamma_c^{1,2}}$ are equivalent in $H^{1/2}(\Gamma^{i}_c)$. By same way, (\ref{ieq:con:3.1:2}) can be verified.
\end{proof}

\begin{proposition}\label{prp:3.3}
Let $\boldsymbol{p}^{i+1}, \tilde{\boldsymbol{p}}^{i+1} \in {V}^{i+1}_{3}(0)$ be the solutions of $({P}_p^{i+1}(\boldsymbol{\lambda}))$ and $({P}_p^{i+1}(\tilde{\boldsymbol{\lambda}}))$ with $\boldsymbol{\lambda}$, $\tilde{\boldsymbol{\lambda}} \in H^{1/2}\left( \Gamma^{1}_c\right) \times \cdots \times H^{1/2} \left(\Gamma^{n-1}_c\right)$, respectively. Then, $\exists {c}_{p}^{i+1}>1/2$ which is independent of $\boldsymbol{\lambda}$ and $\tilde{\boldsymbol{\lambda}}$ such that for $i=1$, we have:
\begin{small}
\begin{equation}\label{ieq:prp:3.3:1}
\left\|\boldsymbol{p}^{2} - \tilde{\boldsymbol{p}}^{2} \right\|_{1, \Omega^2} \leqslant {c}_p^{2} \left(  \left\|\lambda^1-\tilde{\lambda}^1\right\|_{\Gamma_c^{1, 2}} + \left\|\lambda^{2} -\tilde{\lambda}^{2} \right\|_{\Gamma_c^{2, 3}} \right),
\end{equation}
\end{small}
for $i=2,\ldots,n-2$, we have:
\begin{footnotesize}
\begin{equation}\label{ieq:prp:3.3:2}
\left\|\boldsymbol{p}^{i+1} - \tilde{\boldsymbol{p}}^{i+1} \right\|_{1, \Omega^{i+1}} \leqslant {c}_p^{i+1} \left( \left\|\lambda^{i-1} -\tilde{\lambda}^{i-1} \right\|_{\Gamma_c^{i-1, i}} + \left\|\lambda^i-\tilde{\lambda}^i\right\|_{\Gamma_c^{i, i+1}} + \left\|\lambda^{i+1} -\tilde{\lambda}^{i+1} \right\|_{\Gamma_c^{i+1, i+2}} \right),
\end{equation}
\end{footnotesize}
and for $i=n-1$, we have:
\begin{small}
\begin{equation}\label{ieq:prp:3.3:3}
\left\|\boldsymbol{p}^n-\tilde{\boldsymbol{p}}^n\right\|_{1, \Omega^n} \leqslant {c}_p^{n} \left( \left\|\lambda^{n-2} -\tilde{\lambda}^{n-2} \right\|_{\Gamma_c^{n-2, n-1}} + \left\|\lambda^{n-1}-\tilde{\lambda}^{n-1}\right\|_{\Gamma_c^{n-1, n}} \right).
\end{equation}
\end{small}
\end{proposition}

\begin{proof}
It has been known that:
\begin{footnotesize}
\begin{align*}
a^{i+1}\!\left(\boldsymbol{p}^{i+1}, \boldsymbol{v}^{i+1}\right)\!=\!\frac{1}{2}\left(a^{i+1}\left(\boldsymbol{u}^{i+1}, \boldsymbol{v}^{i+1}\right) \!-\! L^{i+1}\left(\boldsymbol{v}^{i+1}\right) \!+\! a^i\left(\boldsymbol{u}^i, \chi_i^i \gamma_i^{i+1} \boldsymbol{v}^{i+1}\right) \!-\! L^i\left(\chi_i^i \gamma_i^{i+1} \boldsymbol{v}^{i+1}\right)\right) \\
a^{i+1}\!\left(\tilde{\boldsymbol{p}}^{i+1}, \boldsymbol{v}^{i+1}\right)\!=\!\frac{1}{2}\left(a^{i+1}\left(\tilde{\boldsymbol{u}}^{i+1}, \boldsymbol{v}^{i+1}\right) \!-\! L^{i+1}\left(\boldsymbol{v}^{i+1}\right) \!+\! a^i\left(\tilde{\boldsymbol{u}}^i, \chi_i^i \gamma_i^{i+1} \boldsymbol{v}^{i+1}\right) \!-\! L^i\left(\chi_i^i \gamma_i^{i+1} \boldsymbol{v}^{i+1}\right)\right)
\end{align*}
\end{footnotesize}
Therefore, the difference $\boldsymbol{p}^{i+1}-\tilde{\boldsymbol{p}}^{i+1} \in V^{i+1}_{3}(0)$ satisfies the equation:
\begin{small}
$$
a^{i+1}\left(\boldsymbol{p}^{i+1}-\tilde{\boldsymbol{p}}^{i+1}, \boldsymbol{v}^{i+1}\right) = \frac{1}{2}\left(a^{i+1}\left(\boldsymbol{u}^{i+1} - \tilde{\boldsymbol{u}}^{i+1}, \boldsymbol{v}^{i+1}\right) +a^i\left(\boldsymbol{u}^i - \tilde{\boldsymbol{u}}^{i}, \chi_i^i \gamma_i^{i+1} \boldsymbol{v}^{i+1}\right) \right)
$$
\end{small}
Hence
\begin{small}
\begin{align*}
&\left\|\boldsymbol{p}^{i+1}-\tilde{\boldsymbol{p}}^{i+1}\right\|^{2}_{1, \Omega^{i+1}} \leqslant \frac{1}{2}\Big( \left\|\boldsymbol{u}^{i+1} - \tilde{\boldsymbol{u}}^{i+1}\right\|_{1,\Omega^{i+1}} \cdot \left\|\boldsymbol{p}^{i+1}-\tilde{\boldsymbol{p}}^{i+1}\right\|_{1, \Omega^{i+1}}\\ 
& + \left\|\boldsymbol{u}^i - \tilde{\boldsymbol{u}}^{i}\right\|_{1,\Omega^{i}} \cdot  c^{i}_{2} c^{i}_{t} \left\|\boldsymbol{p}^{i+1}-\tilde{\boldsymbol{p}}^{i+1} \right\|_{1,\Omega^{i+1}} \Big),
\end{align*}
\end{small}
where the equivalence of norm and trace theorem were used. Then, the following inequality is easily verified:
\begin{small}
\begin{align*}
\left\|\boldsymbol{p}^{i+1}-\tilde{\boldsymbol{p}}^{i+1}\right\|_{1, \Omega^{i+1}} & \leqslant \frac{1}{2}\Big( \left\|\boldsymbol{u}^{i+1} - \tilde{\boldsymbol{u}}^{i+1}\right\|_{1,\Omega^{i+1}} + c^{i}_{2} c^{i}_{t}\cdot \left\|\boldsymbol{u}^i - \tilde{\boldsymbol{u}}^{i}\right\|_{1,\Omega^{i}} \Big).
\end{align*}
\end{small}
Finally, according to Proposition \ref{prp:3.2} and Consequence \ref{con:3.1}, formulations (\ref{ieq:prp:3.3:1}), (\ref{ieq:prp:3.3:2}), (\ref{ieq:prp:3.3:3}) can be verified.
\end{proof}

Analogously, the following result can be proved.

\begin{proposition}\label{prp:3.4}
Let $\boldsymbol{q}^{i}, \tilde{\boldsymbol{q}}^{i} \in {V}^{i}_{2}(0)$ be the solutions of $\left({P}_q^{i}(\boldsymbol{\lambda})\right)$ and $({P}_q^{i}(\tilde{\boldsymbol{\lambda}}))$ with $\boldsymbol{\lambda}$, $\tilde{\boldsymbol{\lambda}} \in H^{1/2}\left( \Gamma^{1}_c\right) \times \cdots \times H^{1/2} \left(\Gamma^{n-1}_c\right)$, respectively. Then, $\exists{c}_{q}^{i}>1/2$ which is independent of $\boldsymbol{\lambda}$ and $\tilde{\boldsymbol{\lambda}}$ such that for $i=1$, we have:
\begin{small}
\begin{equation}\label{ieq:prp:3.4:1}
\left\|\boldsymbol{q}^{1} - \tilde{\boldsymbol{q}}^{1} \right\|_{1, \Omega^1} \leqslant {c}_{q}^{1} \left(  \left\|\lambda^1-\tilde{\lambda}^1\right\|_{\Gamma_c^{1, 2}} + \left\|\lambda^{2} -\tilde{\lambda}^{2} \right\|_{\Gamma_c^{2, 3}} \right),
\end{equation}
\end{small}
for $i=2,\ldots,n-2$, we have:
\begin{small}
\begin{equation}\label{ieq:prp:3.4:2}
\left\|\boldsymbol{q}^{i} - \tilde{\boldsymbol{q}}^{i} \right\|_{1, \Omega^i} \leqslant {c}_{q}^{i} \left( \left\|\lambda^{i-1} -\tilde{\lambda}^{i-1} \right\|_{\Gamma_c^{i-1, i}} + \left\|\lambda^i-\tilde{\lambda}^i\right\|_{\Gamma_c^{i, i+1}} + \left\|\lambda^{i+1} -\tilde{\lambda}^{i+1} \right\|_{\Gamma_c^{i+1, i+2}} \right),
\end{equation}
\end{small}
and for $i=n-1$, we have:
\begin{small}
\begin{equation}\label{ieq:prp:3.4:3}
\left\|\boldsymbol{q}^{n-1}-\tilde{\boldsymbol{q}}^{n-1}\right\|_{1, \Omega^{n-1}} \leqslant {c}_{q}^{n-1} \left( \left\|\lambda^{n-2} -\tilde{\lambda}^{n-2} \right\|_{\Gamma_c^{n-2, n-1}} + \left\|\lambda^{n-1}-\tilde{\lambda}^{n-1}\right\|_{\Gamma_c^{n-1, n}} \right).
\end{equation}
\end{small}
\end{proposition}

\begin{proof}
It has been known that:
\begin{small}
\begin{align*}
& a^i\left(\boldsymbol{q}^i, \boldsymbol{v}^i\right)=\frac{1}{2}\left(a^{i+1}\left(\boldsymbol{u}^{i+1}, \chi_{i+1}^i \gamma_i^i \boldsymbol{v}^i\right)-L^{i+1}\left(\chi_{i+1}^i \gamma_i^i v^i\right)+a^i\left(\boldsymbol{u}^i, \boldsymbol{v}^i\right)-L^i\left(\boldsymbol{v}^i\right)\right) \\
& a^i\left(\tilde{\boldsymbol{q}}^i, \boldsymbol{v}^i\right)=\frac{1}{2}\left(a^{i+1}\left(\tilde{\boldsymbol{u}}^{i+1}, \chi_{i+1}^i \gamma_i^i \boldsymbol{v}^i\right) -L^{i+1}\left(\chi_{i+1}^i \gamma_i^i v^i\right) +a^i\left(\tilde{\boldsymbol{u}}^i, \boldsymbol{v}^i\right) -L^i\left(\boldsymbol{v}^i\right)\right)
\end{align*}
\end{small}
Therefore, the difference $\boldsymbol{q}^{i}-\tilde{\boldsymbol{q}}^{i} \in V^{i}_{2}(0)$ satisfies the equation:
\begin{small}
$$
a^i\left(\boldsymbol{q}^{i}-\tilde{\boldsymbol{q}}^{i}, \boldsymbol{v}^i\right)=\frac{1}{2}\left(a^{i+1}\left(\boldsymbol{u}^{i+1} - \tilde{\boldsymbol{u}}^{i+1}, \chi_{i+1}^i \gamma_i^i \boldsymbol{v}^i\right) + a^i\left(\boldsymbol{u}^i - \tilde{\boldsymbol{u}}^{i} , \boldsymbol{v}^i\right) \right)
$$
\end{small}
Hence
\begin{small}
\begin{align*}
\left\|\boldsymbol{q}^{i}-\tilde{\boldsymbol{q}}^{i}\right\|^{2}_{1, \Omega^{i}} 
\leqslant \frac{1}{2}\Big( \left\|\boldsymbol{u}^{i+1} - \tilde{\boldsymbol{u}}^{i+1}\right\|_{1,\Omega^{i+1}} \cdot \frac{c_{t}^{i}}{c_{1}^{i}} \left\| \boldsymbol{q}^{i}-\tilde{\boldsymbol{q}}^{i} \right\|_{1,\Omega^{i}} + \left\|\boldsymbol{u}^i - \tilde{\boldsymbol{u}}^{i}\right\|_{1,\Omega^{i}} \cdot  \left\|\boldsymbol{q}^{i}-\tilde{\boldsymbol{q}}^{i}\right\|_{1, \Omega^{i}} \Big)
\end{align*}
\end{small}
where the equivalence of norm and trace theorem were used. Then, the following inequality is easily verified:
\begin{small}
$$
\left\|\boldsymbol{q}^{i}-\tilde{\boldsymbol{q}}^{i}\right\|_{1, \Omega^{i}} \leqslant \frac{1}{2}\Big( \frac{c_{t}^{i}}{c_{1}^{i}} \cdot \left\|\boldsymbol{u}^{i+1} - \tilde{\boldsymbol{u}}^{i+1} \right\|_{1,\Omega^{i+1}} + \left\|\boldsymbol{u}^i - \tilde{\boldsymbol{u}}^{i}\right\|_{1,\Omega^{i}} \Big).
$$
\end{small}
Finally, according to Proposition \ref{prp:3.2} and Consequence \ref{con:3.1}, formulations (\ref{ieq:prp:3.4:1}), (\ref{ieq:prp:3.4:2}), (\ref{ieq:prp:3.4:3}) can be proved.
\end{proof}

\begin{proposition}\label{prp:3.5}
Let $\lambda^{i}, \lambda_k^{i} \in H^{1/2}\left(\Gamma_c^{i}\right)$ be such that $\lambda_k^{i} \rightarrow \lambda^{i}$ in $H^{1/2}\left(\Gamma_c^{i}\right)$ as $k \rightarrow \infty$, $i=1,\ldots,n-1$. Then for every $\boldsymbol{v}^i \in {K}^i\left(\lambda^{i} \cdot \boldsymbol{\alpha}^{i+1}\right)$ there exists a sequence $\left\{\boldsymbol{v}_k^i\right\}, \boldsymbol{v}_k^i \in {K}^i\left(\lambda^{i}_{k} \cdot \boldsymbol{\alpha}^{i+1}\right)$ such that
\begin{equation}\label{eq:prp:3.5.1}
\boldsymbol{v}_k^{i} \rightarrow \boldsymbol{v}^{i} \text { in } V^i \text { as } k \rightarrow \infty
\end{equation}
\end{proposition}

\begin{proof}
It is easy to see that ${K}^i\left(\lambda^{i} \cdot \boldsymbol{\alpha}^{i+1}\right)= \chi^i_{i} \lambda^{i} +K^i(0)$, i.e., if $\boldsymbol{v}^i \in {K}^i(\lambda^{i} \cdot \boldsymbol{\alpha}^{i+1})$ then $\boldsymbol{v}^i =\chi^i_{i} \lambda^{i} + \boldsymbol{w}^i$ for some $\boldsymbol{w}^i \in {V}^i$ such that $\boldsymbol{w}^i \cdot \boldsymbol{\beta}^{i} \leq 0$ on $\Gamma_c^{i}$. Then $\boldsymbol{v}^i_{k} := \chi^i_{i} \lambda^{i}_{k} +\boldsymbol{w}^i \in {K}^i \left(\lambda_k^{i} \cdot \boldsymbol{\alpha}^{i+1} \right)$ and the sequence $\left\{\boldsymbol{v}_k^i\right\}$ has the property (\ref{eq:prp:3.5.1}).
\end{proof}

To continue the proof of the theorem, the following decomposition is defined:
\begin{definition}\label{def:3.1}
For $\boldsymbol{v}^i \in V^{i}$, it can be decomposed as
\begin{small}
$$
\boldsymbol{v}^i = \boldsymbol{v}^i_{+} + \boldsymbol{v}^i_{-},
$$ 
\end{small}
where 
\begin{small}
$$
\boldsymbol{v}^i_{-}= \chi^{i-1}_{i} \gamma^{i}_{i-1} \boldsymbol{v}^i \text{ and } \boldsymbol{v}^i_{+}=\boldsymbol{v}^i - \boldsymbol{v}^i_{-}.
$$
\end{small}
\end{definition}

In order to proof the convergence of algorithm, the following proposition should be verified.

\begin{proposition}\label{prp:3.6}
Let $\left\{\left(\boldsymbol{p}_k^{i+1}, \boldsymbol{q}_k^i, \lambda^{i}_{k}\right)\right\}, i=1,\ldots,n-1$ and $\left\{\boldsymbol{u}_k^i\right\}, i=1,\ldots,n$ be the sequence generated by Algorithm \ref{algorithm:3.1} (LD) and suppose that
\begin{small}
\begin{equation}\label{eq:prp:3.6:1}
\lambda_k^{i} \rightarrow \lambda_*^{i} \text { in } H^{1 / 2}\left(\Gamma_c^{i}\right) \text { as } k \rightarrow \infty, i=1,\ldots,n-1.
\end{equation}
\end{small}
Then there exist functions $\boldsymbol{u}_*^i, \boldsymbol{p}_*^i, \boldsymbol{q}_*^i \in {V}^i$ such that
\begin{small}
\begin{equation}\label{eq:prp:3.6:2}
\left.\begin{array}{l}
\boldsymbol{u}_k^i \rightarrow \boldsymbol{u}_*^i ,i=1,\ldots,n\\
\boldsymbol{p}_k^{i+1} \rightarrow \boldsymbol{p}_*^{i+1} ,i=1,\ldots,n-1\\
\boldsymbol{q}_k^i \rightarrow \boldsymbol{q}_*^i ,i=1,\ldots,n-1
\end{array}\right\} \text { in } {V}^i \text { as } k \rightarrow \infty.
\end{equation}
\end{small}
In addition, the function $\boldsymbol{u} = \left(\boldsymbol{u}^1_{*}, \ldots, \boldsymbol{u}^n_{*}\right)$ solves ($P_{0}$).
\end{proposition}

\begin{proof}
Based on the (\ref{eq:prp:3.2:3}), (\ref{ieq:con:3.1:1}), (\ref{ieq:con:3.1:2}), it can be known that $\left\{\boldsymbol{u}_k^i\right\}$ converges to $\left\{\boldsymbol{u}_*^i\right\} \in V^{i}$, $i=1,\ldots,n$. And according to inequalities (\ref{ieq:prp:3.3:1}), (\ref{ieq:prp:3.3:2}), (\ref{ieq:prp:3.3:3}), (\ref{ieq:prp:3.4:1}), (\ref{ieq:prp:3.4:2}) and (\ref{ieq:prp:3.4:3}), the convergence of $\left\{\left(\boldsymbol{p}_k^{i+1}, \boldsymbol{q}_k^i\right)\right\}$ to some $\left\{\left(\boldsymbol{p}_*^{i+1}, \boldsymbol{q}_*^i\right)\right\}$ in $V^{i+1}_{3}(0) \times V^{i}_{2}(0)$ can be verified. When $k \rightarrow \infty$ in ($P^{n}_{d}(\boldsymbol{\lambda}_{k})$), ($P^{i+1}_{p}(\boldsymbol{\lambda}_{k})$) and ($P^{i}_{q}(\boldsymbol{\lambda}_{k})$), $i=1,\ldots,n-1$, it can be known that the limit functions $\left\{\boldsymbol{u}_*^n\right\} \in V^{n}$ and $\left\{\left(\boldsymbol{p}_*^{i+1}, \boldsymbol{q}_*^i\right)\right\} \in V^{i+1}_{3}(0) \times V^{i}_{2}(0)$ can solve ($P^{n}_{d}(\boldsymbol{\lambda}_{*})$), ($P^{i+1}_{p}(\boldsymbol{\lambda}_{*})$) and ($P^{i}_{q}(\boldsymbol{\lambda}_{*})$). Based on the Proposition \ref{prp:3.5}, $\boldsymbol{u}_*^i \in K^i\left(\boldsymbol{\lambda}^i_{*} \cdot \boldsymbol{\alpha}^{i+1}\right)$, $i=1,\ldots,n-1$, i.e. $\boldsymbol{u}_*^i$ is the solution of ($P^{i}_{d}(\boldsymbol{\lambda}_{*})$).

Then, it should be verified that $\boldsymbol{u} = \left(\boldsymbol{u}^1_{*}, \ldots, \boldsymbol{u}^n_{*}\right)$ solves ($P_{0}$). By passing to the limit with $k\to +\infty$ in step 4 of Algorithm \ref{algorithm:3.1}, it is easily to find that $\boldsymbol{p}_*^{i+1} = - \boldsymbol{q}_*^{i}$ on $\Gamma_{c}^{i}$, which implies $\boldsymbol{p}_*^{i+1} = - \boldsymbol{q}_*^{i} \equiv 0$ in $\Omega^{i}, i=1,\ldots,n-1$. 
Indeed, choosing $\pm\boldsymbol{v}^{i}+\boldsymbol{u}^{i}_{*}$, $\boldsymbol{v}^{i}\in V^{i}_{3}(0)\cap V^{i}_{2}(0)$ as a test function in $P_{d}^{i}(\boldsymbol{\lambda}_{*})$, which satisfy $\pm\boldsymbol{v}^{i} + \boldsymbol{u}^{i}_{*} \in K^{i}(\lambda^{i}_{*} \cdot\boldsymbol{\alpha}^{i+1}) \cap V^{i}_{2}(\boldsymbol{u}^i_{*})$, according to \ref{prb:3.4:di} it can be obtained that:
\begin{small}
\begin{align*}
& a^i\left(\boldsymbol{u}^i_{*}, \pm\boldsymbol{v}^i\right) +j^i\left( \pm\boldsymbol{v}^{i} + \boldsymbol{u}^{i}_{*}, \lambda^i_{*} \right) -j^i\left( \boldsymbol{u}^i_{*}, \lambda^i_{*}\right) \geqslant L^i\left(\pm\boldsymbol{v}^i\right)\\
\Rightarrow & a^i\left(\boldsymbol{u}^i_{*}, \boldsymbol{v}^i\right) = L^i\left(\boldsymbol{v}^i\right),~\forall \boldsymbol{v}^{i}\in V^{i}_{3}(0)\cap V^{i}_{2}(0) ,~i=1,\dots,n-1,
\end{align*}
\end{small}
and from the definition of \ref{prb:3.4:dn},
\begin{small}
$$
a^n\left(\boldsymbol{u}^n_{*}, \boldsymbol{v}^n\right) = L^n\left(\boldsymbol{v}^n\right),~\forall \boldsymbol{v}^{n}\in V^{n}_{2}(0).
$$
\end{small}
Therefore, by definition of \ref{prb:3.5:pi}:
\begin{small}
$$
\left\{
\begin{aligned}
&a^{i+1}\left(\boldsymbol{p}^{i+1}_{*}, \boldsymbol{v}^{i+1}\right) = 0,~ &&\forall \boldsymbol{v}^{i+1}\in V^{i+1}_{3}(0)\cap V^{i+1}_{2}(0)\\
&a^{i}\left(\boldsymbol{q}^{i}_{*}, \boldsymbol{v}^{i}\right) = 0,~ &&\forall \boldsymbol{v}^{i}\in V^{i}_{3}(0)\cap V^{i}_{2}(0)
\end{aligned}
\right. ~i=1,\ldots,n-1,
$$
\end{small}
i.e. $\boldsymbol{p}^{i+1}_{*}=0$ on $\Omega^{i+1}/\Gamma_2^{i+1}$ and $\boldsymbol{q}^{i}_{*}=0$ on $\Omega^{i}/\Gamma_{3}^{i}$, $i=1,\ldots,n-1$. Based on this property, the definition of $\chi^{i}_{j}$ and the fact that $\boldsymbol{p}_*^{i+1} = - \boldsymbol{q}_*^{i}$ on $\Gamma_{c}^{i}$, it can be found that $\chi^{i}_{i}\gamma^{i+1}_{i}\boldsymbol{p}^{i+1}_{*} = - \boldsymbol{q}^{i}_{*}$ and $\chi^{i}_{i+1}\gamma^{i}_{i}\boldsymbol{q}^{i}_{*} = - \boldsymbol{p}^{i+1}_{*}$. 
Inserting $\boldsymbol{v}^{i+1} = \boldsymbol{p}^{i+1}_{*}$ into \ref{prb:3.5:pi} and $\boldsymbol{v}^{i} = \boldsymbol{q}^{i}_{*}$ into \ref{prb:3.5:qi} with $\boldsymbol{\lambda}=\boldsymbol{\lambda}_{*}$, we have
\begin{small}
\begin{align*}
a^{i+1}\left(\boldsymbol{p}^{i+1}_{*}, \boldsymbol{p}^{i+1}_{*}\right) 
=& \frac{1}{2} \Big(a^{i+1}\left(\boldsymbol{u}^{i+1}_{*}, \boldsymbol{p}^{i+1}_{*}\right)- L^{i+1}\left(\boldsymbol{p}^{i+1}_{*}\right)- a^{i}\left(\boldsymbol{u}^{i}_{*}, \boldsymbol{q}^{i}_{*}\right) + L^{i}\left(\boldsymbol{q}^{i}_{*}\right) \Big)\\
a^{i}\left(\boldsymbol{q}^{i}_{*}, \boldsymbol{q}^{i}_{*}\right) 
=&\frac{1}{2} \Big(-a^{i+1}\left(\boldsymbol{u}^{i+1}_{*}, \boldsymbol{p}^{i+1}_{*}\right)+ L^{i+1}\left( \boldsymbol{p}^{i+1}_{*}\right) + a^{i}\left(\boldsymbol{u}^{i}_{*},  \boldsymbol{q}^{i}_{*}\right) - L^{i}\left(\boldsymbol{q}^{i}_{*}\right) \Big)
\end{align*}
\end{small}
which means that 
\begin{small}
$$
a^{i+1}\left(\boldsymbol{p}^{i+1}_{*}, \boldsymbol{p}^{i+1}_{*}\right) = - a^{i}\left(\boldsymbol{q}^{i}_{*}, \boldsymbol{q}^{i}_{*}\right).
$$
\end{small}
Hence $\boldsymbol{p}^{i+1}_{*}=0$ on $\Omega^{i+1}$ and $\boldsymbol{q}^{i}_{*}=0$ on $\Omega^{i}$, $i=1,\ldots,n-1$. Then, based on \ref{prb:3.5:pi}, when $\boldsymbol{\lambda}=\boldsymbol{\lambda}_{*}$, the right hand side of equation vanishes:
\begin{small}
$$
a^{i}\left(\boldsymbol{u}^{i}_{*}, \boldsymbol{v}^{i}\right) = L^{i}\left(\boldsymbol{v}^{i}\right)- a^{i-1}\left(\boldsymbol{u}^{i-1}_{*}, \chi_{i-1}^{i-1} \gamma^{i}_{i-1}\boldsymbol{v}^{i}\right) + L^{i-1}\left(\chi_{i-1}^{i-1} \gamma^{i}_{i-1}\boldsymbol{v}^{i}\right), ~\forall \boldsymbol{v}^{i}\in V_{3}^{i}(0).
$$
\end{small}
When $i=n$, we get
\begin{small}
$$
a^{n}\left(\boldsymbol{u}^{n}_{*}, \boldsymbol{v}^{n}\right) = L^{n}\left(\boldsymbol{v}^{n}\right)- a^{n-1}\left(\boldsymbol{u}^{n-1}_{*}, \chi_{n-1}^{n-1} \gamma^{n}_{n-1}\boldsymbol{v}^{n}\right) + L^{n-1}\left(\chi_{n-1}^{n-1} \gamma^{n}_{n-1}\boldsymbol{v}^{n}\right), ~\forall \boldsymbol{v}^{n}\in V^{n},
$$
\end{small}
which is equivalent to (\ref{ieq:prp:3.1:1}) whit $i=n$.
When $i\ne 1$ and $n$, on $\Gamma_{c}^{i-1}$, we have
\begin{small}
$$
a^{i}\left(\boldsymbol{u}^{i}_{*}, \boldsymbol{v}^{i}\right) = L^{i}\left(\boldsymbol{v}^{i}\right)- a^{i-1}\left(\boldsymbol{u}^{i-1}_{*}, \chi_{i-1}^{i-1} \gamma^{i}_{i-1}\boldsymbol{v}^{i}\right) + L^{i-1}\left(\chi_{i-1}^{i-1} \gamma^{i}_{i-1}\boldsymbol{v}^{i}\right), ~\forall \boldsymbol{v}^{i}\in V^{i}
$$
\end{small}
and on $\Omega^{i}/\Gamma_{c}^{i-1}$, we have
\begin{small}
$$
a^{i-1}\left(\boldsymbol{u}^{i-1}_{*}, \chi_{i-1}^{i-1} \gamma^{i}_{i-1}\boldsymbol{v}^{i}\right) = L^{i-1}\left(\chi_{i-1}^{i-1} \gamma^{i}_{i-1}\boldsymbol{v}^{i}\right) = 0, ~\forall \boldsymbol{v}^{i}\in V^{i}.
$$
\end{small}
In \ref{prb:3.4:di} with $\boldsymbol{\lambda} = \boldsymbol{\lambda}_{*}$, it should be noticed that $\boldsymbol{v}^{i}\in K^{i}(\lambda^{i}_{*}\cdot \boldsymbol{\alpha}^{i+1}) \cap V_{2}^{i}(\boldsymbol{u}^{i}_{*})$, $\boldsymbol{u}^{i}_{*} = \lambda^{i-1}_{*}$ on $\Gamma_{c}^{i-1}$. Therefore, on $\Gamma_{c}^{i-1}$, we have
\begin{small}
$$
a^i\left(\boldsymbol{u}^i_{*}, \boldsymbol{v}^i-\boldsymbol{u}^i_{*}\right) = j^i\left( \boldsymbol{v}^i, \lambda^i_{*}\right)-j^i\left( \boldsymbol{u}^i_{*}, \lambda^i_{*}\right) = L^i\left(\boldsymbol{v}^i-\boldsymbol{u}^i_{*}\right) = 0,
$$
\end{small}
and on $\Omega^{i}/\Gamma_{c}^{i-1}$, we have:
\begin{small}
$$
a^i\left(\boldsymbol{u}^i_{*}, \boldsymbol{v}^i-\boldsymbol{u}^i_{*}\right)+j^i\left( \boldsymbol{v}^i, \lambda^i_{*}\right) -j^i\left( \boldsymbol{u}^i_{*}, \lambda^i_{*}\right) \geqslant L^i\left(\boldsymbol{v}^i-\boldsymbol{u}^i_{*}\right).
$$
\end{small}
Then, $\forall \boldsymbol{v}^i \in K^{i}(\lambda^{i}_{*}\cdot \boldsymbol{\alpha}^{i+1})$, it can be decomposed as $\boldsymbol{v}^i = \boldsymbol{v}^i_{+} + \boldsymbol{v}^i_{-}$, where $\boldsymbol{v}^i_{-}= \chi^{i-1}_{i} \gamma^{i}_{i-1} \boldsymbol{v}^i$. Therefore, we have 
\begin{small}
\begin{align*}
& a^{i}\left(\boldsymbol{u}^{i}_{*}, \left(\boldsymbol{v}^i-\boldsymbol{u}^i_{*}\right)_{-}\right) = L^{i}\left(\left(\boldsymbol{v}^i-\boldsymbol{u}^i_{*}\right)_{-}\right)- a^{i-1}\left(\boldsymbol{u}^{i-1}_{*}, \chi_{i-1}^{i-1} \gamma^{i}_{i-1} \left(\boldsymbol{v}^i-\boldsymbol{u}^i_{*}\right)\right)\\ 
&+ L^{i-1}\left(\chi_{i-1}^{i-1} \gamma^{i}_{i-1} \left(\boldsymbol{v}^i-\boldsymbol{u}^i_{*}\right) \right), \\
& a^i\left(\boldsymbol{u}^i_{*}, \left(\boldsymbol{v}^i-\boldsymbol{u}^i_{*}\right)_{+} \right) +j^i\left( \boldsymbol{v}^i, \lambda^i_{*}\right) -j^i\left( \boldsymbol{u}^i_{*}, \lambda^i_{*}\right) \geqslant L^i\left(\left(\boldsymbol{v}^i-\boldsymbol{u}^i_{*}\right)_{+}\right).
\end{align*}
\end{small}
By adding the two formulas above, we get:
\begin{small}
\begin{align*}
&a^i\left(\boldsymbol{u}^i_{*}, \boldsymbol{v}^i-\boldsymbol{u}^i_{*} \right) +j^i\left( \boldsymbol{v}^i, \lambda^i_{*}\right) -j^i\left( \boldsymbol{u}^i_{*}, \lambda^i_{*}\right)\geqslant \\
&L^i\left(\boldsymbol{v}^i-\boldsymbol{u}^i_{*}\right) - a^{i-1}\left(\boldsymbol{u}^{i-1}_{*}, \chi_{i-1}^{i-1} \gamma^{i}_{i-1} \left(\boldsymbol{v}^i-\boldsymbol{u}^i_{*}\right)\right) + L^{i-1}\left(\chi_{i-1}^{i-1} \gamma^{i}_{i-1} \left(\boldsymbol{v}^i-\boldsymbol{u}^i_{*}\right) \right),
\end{align*}
\end{small}
which is equivalent to (\ref{ieq:prp:3.1:1}) with $i\neq 1$ and $n$, because $\lambda^i_{*} = \boldsymbol{u}^{i+1}_{*}$ on $\Gamma_{c}^{i}$.

Finally, based on Proposition \ref{prp:3.1}, \ref{prb:3.4:d1} is equivalent to (\ref{ieq:prp:3.1:1}) with $i=1$.
\end{proof}

\begin{proof}[Proof of Theorem \ref{thm:3.1}]

Based on Proposition \ref{prp:3.6}, it is sufficient to verify that the sequence $\{\lambda_{k}^{i}\}$ calculated from step 4 of Algorithm \ref{algorithm:3.1} is convergent if $\theta\in(0,\theta_{*})$. For this end, some relevant definitions need to be introduced: the mapping  $T^{i}_{\theta}: H^{1/2}(\Gamma_{c}^{i}) \to H^{1/2}(\Gamma_{c}^{i})$ $(i=1,\ldots,n-1)$ defined by $T^{i}_{\theta}\lambda^{i} = \lambda^{i} - \theta T^{i}\lambda^{i}$, where $\theta>0$ and $T^{i}: H^{1/2}(\Gamma_{c}^{i}) \to H^{1/2}(\Gamma_{c}^{i})$ defined by $T^{i} \lambda^{i} = \gamma^{i+1}_{i}\boldsymbol{p}^{i+1} + \gamma^{i}_{i}\boldsymbol{q}^{i}$ on $\Gamma_{c}^{i}$ with $\boldsymbol{p}^{i+1}$, $\boldsymbol{q}^{i}$ being the solutions of \ref{prb:3.5:pi} and \ref{prb:3.5:qi}, respectively. Let $\boldsymbol{p}=\left(\boldsymbol{p}^{2},\ldots,\boldsymbol{p}^{n}\right)$, $\boldsymbol{q}=\left(\boldsymbol{q}^{1},\ldots,\boldsymbol{q}^{n-1}\right)$, $\boldsymbol{\lambda} =\left(\lambda^{1}, \ldots, \lambda^{n-1}\right)$. Then, the mapping $T_{\theta}: H^{1/2}(\Gamma_{c}) \to H^{1/2}(\Gamma_{c})$ defined by 
\begin{small}
$$
T_{\theta}\boldsymbol{\lambda} = \boldsymbol{\lambda} - \theta T\boldsymbol{\lambda} = (T_{\theta}^{1}\lambda^{1}, \ldots, T_{\theta}^{n-1}\lambda^{n-1}),
$$ 
\end{small}
where $T: H^{1/2}(\Gamma_{c}) \to H^{1/2}(\Gamma_{c})$ defined by 
\begin{small}
$$
T \boldsymbol{\lambda} = \left(T^{1} \lambda^{1},\ldots, T^{n-1} \lambda^{n-1}\right) = \left( \gamma^{2}_{1}\boldsymbol{p}^{2} + \gamma^{1}_{1}\boldsymbol{q}^{1}, \ldots, \gamma^{n}_{n-1} \boldsymbol{p}^{n} + \gamma^{n-1}_{n-1}\boldsymbol{q}^{n-1} \right).
$$ 
\end{small}
Furthermore, norm $\left\|\cdot\right\|_{\Gamma_{c}}$ is defined as $\left\|\boldsymbol{\lambda}\right\|_{\Gamma_{c}}=\sqrt{\sum_{i=1}^{n-1} \left\|\lambda^{i}\right\|_{\Gamma_{c}^{i,i+1}}^{2}}$. In what follows, it will be proved that $T_\theta$ is compression mapping in $H^{\frac{1}{2}}\left(\Gamma_C\right)$ for $\theta$ small enough.

Let $\boldsymbol{\lambda},\bar{\boldsymbol{\lambda}} \in H^{1/2}(\Gamma_{c}^{i})$ and $\boldsymbol{p}^{i+1}$, $\bar{\boldsymbol{p}}^{i+1}$ be solutions of \ref{prb:3.5:pi} with $\boldsymbol{\lambda}=\boldsymbol{\lambda}$ and $\bar{\boldsymbol{\lambda}}$, respectively. Similarly, let $\boldsymbol{q}^{i}$, $\bar{\boldsymbol{q}}^{i}$ be solutions of \ref{prb:3.5:qi} with $\boldsymbol{\lambda}=\boldsymbol{\lambda}$ and $\bar{\boldsymbol{\lambda}}$, respectively. Then
\begin{small}
\begin{align}
&\left\| T_{\theta} \boldsymbol{\lambda} - T_{\theta} \bar{\boldsymbol{\lambda}} \right\|_{\Gamma_{c}}^{2} =
\sum_{i=1}^{n-1} \left\| T_{\theta}^{i} \lambda^{i} - T_{\theta}^{i} \bar{\lambda}^{i} \right\|_{\Gamma_{c}^{i,i+1}}^{2} \nonumber\\
= & \sum_{i=1}^{n-1} a^{i+1} \Big( \chi^{i}_{i+1}\left( \lambda^{i} - \bar{\lambda}^{i} - \theta \left( T^{i}\lambda^{i} - T^{i}\bar{\lambda}^{i} \right) \right),\chi^{i}_{i+1}\left( \lambda^{i} - \bar{\lambda}^{i} - \theta \left( T^{i}\lambda^{i} - T^{i}\bar{\lambda}^{i} \right) \right)\Big) \label{ieq:thm:3.1:prf:1}\\
= & \sum_{i=1}^{n-1} \Big( \left\|\lambda^{i} - \bar{\lambda}^{i} \right\|_{\Gamma_{c}^{i,i+1}}^{2} + \theta^{2} \left\| T^{i} \lambda^{i} - T^{i} \bar{\lambda}^{i}  \right\|_{\Gamma_{c}^{i,i+1}}^{2} \nonumber \\
&~~~~~~~~~~~~ -2\theta a^{i+1} \left( \chi^{i}_{i+1}\left( \lambda^{i} - \bar{\lambda}^{i} \right), \chi^{i}_{i+1}\left( T^{i} \lambda^{i} - T^{i} \bar{\lambda}^{i} \right)\right) \Big). \nonumber
\end{align}
\end{small}
The second and third term on right hand side of (\ref{ieq:thm:3.1:prf:1}) should be estimated. From the definition of $T^{i}$, formulations (\ref{ieq:prp:3.3:1}) - (\ref{ieq:prp:3.4:3}), the following inequalities can be reached:
\begin{small}
\begin{align}
&\left\| T^{i} \lambda^{i} - T^{i} \bar{\lambda}^{i} \right\|_{\Gamma_{c}^{i,i+1}}
= \left\| \gamma^{i+1}_{i}\boldsymbol{p}^{i+1} - \gamma^{i+1}_{i} \bar{ \boldsymbol{p}}^{i+1} + \gamma^{i}_{i}\boldsymbol{q}^{i} - \gamma^{i}_{i} \bar{\boldsymbol{q}}^{i} \right\|_{\Gamma_{c}^{i,i+1}} \nonumber \\
\leqslant& \left\| \gamma^{i+1}_{i}\boldsymbol{p}^{i+1} - \gamma^{i+1}_{i} \bar{ \boldsymbol{p}}^{i+1}\right\|_{\Gamma_{c}^{i,i+1}} + \frac{1}{c^{i}_{1}}
\left\|\gamma^{i}_{i}\boldsymbol{q}^{i} - \gamma^{i}_{i} \bar{\boldsymbol{q}}^{i} \right\|_{\Gamma_{c}^{i,i}}  \label{ieq:thm:3.1:prf:2}\\
\leqslant& \bar{c}^{i} \left( \left\|\lambda^{i-1} -\bar{\lambda}^{i-1} \right\|_{\Gamma_c^{i-1, i}} + \left\|\lambda^i-\tilde{\lambda}^i\right\|_{\Gamma_c^{i, i+1}} + \left\|\lambda^{i+1} -\tilde{\lambda}^{i+1} \right\|_{\Gamma_c^{i+1, i+2}} \right) \nonumber
\end{align}
\end{small}
where $\bar{c}^{i}\geqslant \frac{c^{min}_{t}}{2}$ as follows from proposition \ref{prp:3.3} and \ref{prp:3.4}.

The estimate of the last term of (\ref{ieq:thm:3.1:prf:1}) is more involved. Select $\pm \boldsymbol{v}^{i}+\boldsymbol{u}^{i}$ with $\boldsymbol{v}^{i}\in V^{i}_{3}(0)\cap V^{i}_{2}(0)$ as a test function in \ref{prb:3.4:di}, which satisfy $\pm \boldsymbol{v}^{i}+\boldsymbol{u}^{i} \in K^{i}(\lambda^{i}\cdot\boldsymbol{\alpha}^{i+1}) \cap V^{i}_{2}(\boldsymbol{u}^{i})$. According to definitions of \ref{prb:3.4:di} and \ref{prb:3.4:dn}, we have
\begin{small}
$$
a^{i}(\boldsymbol{u}^{i},\boldsymbol{v}^{i}) = L^{i}(\boldsymbol{v}^{i}),~~\forall \boldsymbol{v}^{i}\in V^{i}_{3}(0)\cap V^{i}_{2}(0),~~i=1,\ldots,n.
$$
\end{small}
From the definition of \ref{prb:3.5:pi}, it can be known that $a^{i+1}(\boldsymbol{p}^{i+1}, \boldsymbol{v}^{i+1}) = 0$, $\forall\boldsymbol{v}^{i+1}\in V^{i+1}_{3}(0)\cap V^{i+1}_{2}(0)$. Since $\boldsymbol{p}^{i+1}\in V_{3}^{i+1}(0)$, the above formulation can be reduced to
\begin{small}
$$
a^{i+1}(\boldsymbol{p}^{i+1}, \boldsymbol{v}^{i+1}) = 0 \Rightarrow a^{i+1}(\boldsymbol{p}^{i+1}- \bar{\boldsymbol{p}}^{i+1}, \boldsymbol{v}^{i+1}) = 0,~~\forall\boldsymbol{v}^{i+1}\in V^{i+1}_{2}(0).
$$
\end{small}
Therefore
\begin{small}
$$
\chi^{i}_{i+1}\gamma^{i+1}_{i} (\boldsymbol{p}^{i+1}- \bar{\boldsymbol{p}}^{i+1}) = \boldsymbol{p}^{i+1}- \bar{\boldsymbol{p}}^{i+1}.
$$
\end{small}
Then, the third term of (\ref{ieq:thm:3.1:prf:1}) can be expressed as: 
\begin{small}
\begin{align}
&a^{i+1} \left( \chi^{i}_{i+1}\left( \lambda^{i} - \bar{\lambda}^{i} \right), \chi^{i}_{i+1}\left( T^{i} \lambda^{i} - T^{i} \bar{\lambda}^{i} \right)\right) = a^{i+1} \left( \chi^{i}_{i+1}\left( \lambda^{i} - \bar{\lambda}^{i} \right), \boldsymbol{p}^{i+1} -\bar{\boldsymbol{p}}^{i+1} \right)\label{ieq:thm:3.1:prf:3}\\
&+ a^{i+1} \left( \chi^{i}_{i+1}\left( \lambda^{i} - \bar{\lambda}^{i} \right), \chi^{i}_{i+1}\left( \gamma^{i}_{i}\boldsymbol{q}^{i} - \gamma^{i}_{i} \bar{\boldsymbol{q}}^{i} \right)\right). \nonumber
\end{align}
\end{small}

Since $\chi^{i}_{i+1}\left( \lambda^{i} - \bar{\lambda}^{i} \right) \in V^{i+1}_{3}(0)$, it can be used as a test function in \ref{prb:3.5:pi}. So, the first term of (\ref{ieq:thm:3.1:prf:3}) satisfy:
\begin{small}
\begin{equation}\label{ieq:thm:3.1:prf:4}
\begin{aligned}
& a^{i+1} \left( \chi^{i}_{i+1}\left( \lambda^{i} - \bar{\lambda}^{i} \right), \boldsymbol{p}^{i+1} -\bar{\boldsymbol{p}}^{i+1} \right) = a^{i+1} \left( \boldsymbol{p}^{i+1} -\bar{\boldsymbol{p}}^{i+1}, \chi^{i}_{i+1}\left( \lambda^{i} - \bar{\lambda}^{i} \right) \right)\\
=& \frac{1}{2} \Big(a^{i+1}\left(\boldsymbol{u}^{i+1}- \bar{\boldsymbol{u}}^{i+1}, \chi^{i}_{i+1}\left( \lambda^{i} - \bar{\lambda}^{i} \right)\right) + a^{i}\left(\boldsymbol{u}^{i}- \bar{\boldsymbol{u}}^{i}, \chi_{i}^{i} \left( \lambda^{i} - \bar{\lambda}^{i} \right) \right) \Big)
\end{aligned}
\end{equation}
\end{small}
From the definition of \ref{prb:3.4:di}, we can know that $\boldsymbol{u}^{i+1}=\lambda^{i}$ on $\Gamma_{c}^{i}$. Therefore, the first term of (\ref{ieq:thm:3.1:prf:4}) satisfy:
\begin{small}
$$
a^{i+1}\left(\boldsymbol{u}^{i+1}- \bar{\boldsymbol{u}}^{i+1}, \chi^{i}_{i+1}\left( \lambda^{i} - \bar{\lambda}^{i} \right)\right) = \left\| \lambda^{i} - \bar{\lambda}^{i} \right\|^{2}_{\Gamma_{c}^{i,i+1}}.
$$
\end{small}
Then, based on proposition \ref{prp:3.2}, the second term of (\ref{ieq:thm:3.1:prf:4}) can be represented as follows: 
\begin{small}
\begin{align*}
a^{i}\left(\boldsymbol{u}^{i}- \bar{\boldsymbol{u}}^{i}, \chi_{i}^{i} \left( \lambda^{i} - \bar{\lambda}^{i} \right) \right)
\geqslant  a^{i}\left(\boldsymbol{u}^{i}- \bar{\boldsymbol{u}}^{i}, \left(\boldsymbol{u}^{i} - \chi^{i-1}_{i} \lambda^{i-1} \right) - \left(\bar{\boldsymbol{u}}^{i} - \chi^{i-1}_{i} \bar{\lambda}^{i-1} \right) \right)
\end{align*}
\end{small}
where $i=2,\ldots,n-1$. Since $\boldsymbol{u}^{i}=\lambda^{i-1}$ on $\Gamma_{c}^{i-1}$, 
\begin{small}
$$
\boldsymbol{u}^{i} - \chi^{i-1}_{i} \lambda^{i-1} = \boldsymbol{u}^{i} - \chi^{i-1}_{i} \gamma^{i}_{i-1} \boldsymbol{u}^{i} = \boldsymbol{u}^{i}_{+}.
$$
\end{small}
Then 
\begin{small}
$$
a^{i}\left(\boldsymbol{u}^{i}- \bar{\boldsymbol{u}}^{i}, \chi_{i}^{i} \left( \lambda^{i} - \bar{\lambda}^{i} \right) \right) \geqslant a^{i}\left(\boldsymbol{u}^{i}- \bar{\boldsymbol{u}}^{i}, \boldsymbol{u}^{i}_{+}- \bar{\boldsymbol{u}}^{i}_{+} \right) = 
a^{i}\left(\boldsymbol{u}^{i}_{+}- \bar{\boldsymbol{u}}^{i}_{+}, \boldsymbol{u}^{i}_{+}- \bar{\boldsymbol{u}}^{i}_{+} \right).
$$
\end{small}
For case $i=1$, according to proposition \ref{prp:3.2}, there are:
\begin{small}
$$
a^{1}\left(\boldsymbol{u}^{1}- \bar{\boldsymbol{u}}^{1}, \chi_{1}^{1} \left( \lambda^{1} - \bar{\lambda}^{1} \right) \right) \geqslant a^{1}\left(\boldsymbol{u}^{1}- \bar{\boldsymbol{u}}^{1}, \boldsymbol{u}^{1}- \bar{\boldsymbol{u}}^{1} \right) \geqslant
a^{1}\left(\boldsymbol{u}^{1}_{+}- \bar{\boldsymbol{u}}^{1}_{+}, \boldsymbol{u}^{1}_{+}- \bar{\boldsymbol{u}}^{1}_{+} \right).
$$
\end{small}
Based on the above analysis, the first term (\ref{ieq:thm:3.1:prf:3}) satisfies the following inequality relations:
\begin{small}
\begin{align}
a^{i+1} \left( \chi^{i}_{i+1}\left( \lambda^{i} - \bar{\lambda}^{i} \right), \boldsymbol{p}^{i+1} -\bar{\boldsymbol{p}}^{i+1} \right) \geqslant \frac{1}{2} \left( \left\| \lambda^{i} - \bar{\lambda}^{i} \right\|^{2}_{\Gamma_{c}^{i,i+1}} + a^{i}\left(\boldsymbol{u}^{i}_{+}- \bar{\boldsymbol{u}}^{i}_{+}, \boldsymbol{u}^{i}_{+}- \bar{\boldsymbol{u}}^{i}_{+} \right)
\right). \label{ieq:thm:3.1:prf:5}
\end{align}
\end{small}

Then consider the second term of formula (\ref{ieq:thm:3.1:prf:3}). Since 
\begin{small}
$$
a^{i+1}\left(\chi^{i}_{i+1}\left( \gamma^{i}_{i}\boldsymbol{q}^{i} - \gamma^{i}_{i} \bar{\boldsymbol{q}}^{i} \right), \boldsymbol{v}^{i+1} \right)= 0,~~\forall \boldsymbol{v}^{i+1} \in V_{2}^{i+1}(0),
$$
\end{small}
it is easy to verify:
\begin{small}
$$
a^{i+1} \left( \chi^{i}_{i+1}\left( \lambda^{i} - \bar{\lambda}^{i} \right), \chi^{i}_{i+1}\left( \gamma^{i}_{i}\boldsymbol{q}^{i} - \gamma^{i}_{i} \bar{\boldsymbol{q}}^{i} \right)\right) = a^{i+1} \left( \boldsymbol{u}^{i+1} - \bar{\boldsymbol{u}}^{i+1}, \chi^{i}_{i+1}\left( \gamma^{i}_{i}\boldsymbol{q}^{i} - \gamma^{i}_{i} \bar{\boldsymbol{q}}^{i} \right)\right)
$$
\end{small}
By inserting $\boldsymbol{v}^{i} = \boldsymbol{q}^{i} - \bar{\boldsymbol{q}}^{i} \in V_{2}^{i}(0)$ into \ref{prb:3.5:qi} with $\boldsymbol{\lambda}=\boldsymbol{\lambda}$ and $\bar{\boldsymbol{\lambda}}$, we have
\begin{footnotesize}
\begin{equation*}
a^{i+1} \left( \chi^{i}_{i+1}\left( \lambda^{i} - \bar{\lambda}^{i} \right), \chi^{i}_{i+1}\left( \gamma^{i}_{i}\boldsymbol{q}^{i} - \gamma^{i}_{i} \bar{\boldsymbol{q}}^{i} \right)\right) 
= 2 a^{i}\left(\boldsymbol{q}^{i} - \bar{\boldsymbol{q}}^{i}, \boldsymbol{q}^{i} - \bar{\boldsymbol{q}}^{i}\right) - a^{i}\left(\boldsymbol{u}^{i} - \bar{\boldsymbol{u}}^{i+1},  \boldsymbol{q}^{i} - \bar{\boldsymbol{q}}^{i} \right)
\end{equation*}
\end{footnotesize}
Since $\boldsymbol{q}^{i},\bar{\boldsymbol{q}}^{i} \in V_{2}^{i}(0)$, 
\begin{small}
$$
\boldsymbol{q}^{i}_{-} = \chi^{i-1}_{i} \gamma^{i}_{i-1}\boldsymbol{q}^{i} =0 \Rightarrow 
\boldsymbol{q}^{i}_{+} = \boldsymbol{q}^{i}.
$$
\end{small}
Similarly, $\bar{\boldsymbol{q}}^{i}_{+} = \bar{\boldsymbol{q}}^{i}$. Then
\begin{small}
\begin{equation}\label{eq:thm:3.1:prf:6}
\begin{aligned}
& a^{i+1} \left( \chi^{i}_{i+1}\left( \lambda^{i} - \bar{\lambda}^{i} \right), \chi^{i}_{i+1}\left( \gamma^{i}_{i}\boldsymbol{q}^{i} - \gamma^{i}_{i} \bar{\boldsymbol{q}}^{i} \right)\right) \\
= &2 a^{i}\left(\boldsymbol{q}^{i}_{+} - \bar{\boldsymbol{q}}^{i}_{+}, \boldsymbol{q}^{i}_{+} - \bar{\boldsymbol{q}}^{i}_{+} \right) - a^{i}\left(\boldsymbol{u}^{i}_{+} - \bar{\boldsymbol{u}}^{i+1}_{+},  \boldsymbol{q}^{i}_{+} - \bar{\boldsymbol{q}}^{i}_{+} \right).
\end{aligned}
\end{equation}
\end{small}
Therefore, by substituting formula (\ref{ieq:thm:3.1:prf:5}) and (\ref{eq:thm:3.1:prf:6}) into formula (\ref{ieq:thm:3.1:prf:3}), the following inequality relationship is established:
\begin{small}
\begin{equation}\label{ieq:thm:3.1:prf:7}
a^{i+1} \left( \chi^{i}_{i+1}\left( \lambda^{i} - \bar{\lambda}^{i} \right), \chi^{i}_{i+1}\left( T^{i} \lambda^{i} - T^{i} \bar{\lambda}^{i} \right)\right) 
\geqslant \frac{1}{2} \left\| \lambda^{i} - \bar{\lambda}^{i} \right\|^{2}_{\Gamma_{c}^{i,i+1}}
\end{equation}
\end{small}
Finally, based on formula (\ref{ieq:thm:3.1:prf:2}) and (\ref{ieq:thm:3.1:prf:7}), formula (\ref{ieq:thm:3.1:prf:1}) satisfies:
\begin{small}
\begin{equation*}
\left\| T_{\theta} \boldsymbol{\lambda} - T_{\theta} \bar{\boldsymbol{\lambda}} \right\|_{\Gamma_{c}}^{2}  
\leqslant \left( 1 + 9\theta^{2} \bar{c}^{2} - \theta \right) \left\|\boldsymbol{\lambda} - \bar{\boldsymbol{\lambda}} \right\|_{\Gamma_{c}}^{2}.
\end{equation*}
\end{small}
where $\bar{c}=\max\{\bar{c}^{1},\ldots,\bar{c}^{n-1}\} \geqslant \frac{c^{\min}_{t}}{2}$ and $\theta>0$. Thus, if $0<\theta<\theta_{*} := \frac{1}{9\bar{c}^{2}}$, the mapping $T_{\theta}$ is contraction. Since  $\bar{c}\geqslant \frac{c^{\min}_{t}}{2}$, it's easy to prove that $\theta_{*}\in (0, 4/(3c_{t}^{\min})^2)$.
The method of successive approximations
\begin{small}
\begin{equation}\label{ieq:thm:3.1:prf:8}
\left\{\begin{array}{l}
\boldsymbol{\lambda}_0 \in H^{1 / 2}\left(\Gamma_c\right) \\
\boldsymbol{\lambda}_{k+1}=T_\theta \boldsymbol{\lambda}_k, \quad k=0,1, \ldots
\end{array}\right.
\end{equation}
\end{small}
is convergent to a unique fixed-point of $T_\theta$ in $H^{1 / 2}\left(\Gamma_c\right)$ for arbitrary $\boldsymbol{\lambda}_0 \in H^{1 / 2}\left(\Gamma_c\right)$, and (\ref{ieq:thm:3.1:prf:8}) is equivalent to the Algorithm \ref{algorithm:3.1}.

\end{proof}

Now the main result of this section can be proved.

\section{Discrete form of the layer decomposition algorithm}

In this section, discrete form of Algorithm \ref{algorithm:3.1} will be constructed. Then, the convergence theorem of the algorithm will be proved, where the choice of the upper bound $\theta_{*}$ of the parameter $\theta$ can be independent of the mesh size $h$.

According to the physical model in Fig.\ref{3:fig1}, the multi-layer elastic contact system is composed of $n$-layer elastic bodies. Then a non-degenerate quasi-uniform subdivision (tetrahedrons) of the $i$-th layer $\Omega^i$ is denoted as $\mathcal{E}_h^i=\left\{E_1^i, E_2^i, \ldots, E_{N_h^i}^i\right\}$, where $E_j^i$ $\left(i=1,2, \ldots, n, j=1,2, \ldots, N_h^i\right)$ represents a regular finite element, $h$ is the maximum diameter of all elements and $N_h^i$ is the number of finite elements. 
Therefore, the finite element division of the whole system is expressed as: $\mathcal{E}_h=\cup_{i=1}^n \mathcal{E}_h^i$. 
Furthermore, it can be supposed that $\mathcal{E}_h^i$ and $\mathcal{E}_h^{i+1}$ are compatible on $\Gamma_{c}^{i}$, $\forall h>0$ and $i=1,\ldots,n-1$, i.e., the nodes of $\mathcal{E}_h^i$, $\mathcal{E}_h^{i+1}$ on $\Gamma_{c}^{i}$ coincide.
For each elastic body $\Omega^i$, the finite element division of its boundary is denoted as: ${\Gamma}_{rh}^i=\cup_{j=1}^{I(i, r)} \bar{\Gamma}_{rh, j}^i(r=2,3, i=1,2, \ldots, n$ and $j=1,2, \ldots, I(i, r))$, where each piece $\bar{\Gamma}_{rh, j}^i$ represented by an affine function. The finite dimension space consisting of piecewise linear polynomial functions, denoted as $V_h^{i} \subset V^i$, are defined on subdivision $\mathcal{E}_h^i$, i.e.
\begin{small}
$$
V_{h}^{i} = \left\{ \boldsymbol{v}^{i}_{h} \in \left( C(\bar{\Omega}^{i}) \right)^{d} ~\Big|~ \boldsymbol{v}^{i}_{h}|_{E_j^i} \in \left(P_{1}(E_j^i)\right)^{d},~\forall j=1, \ldots, N_h^i, ~ \boldsymbol{v}^{i}_{h} = 0 \text{ on } \Gamma^{i}_{1}\right\}
$$
\end{small}
where $i=1, \ldots, n$, and $V_h=V_h^{1} \times V_h^{2} \times \cdots \times V_h^{n}$.
According to the previous definition, the finite element space of $\mathcal{K}$ is defined as:
\begin{small}
$$
\mathcal{K}_{h}=\left\{\boldsymbol{v}_{h} \in V_{h} ~\big|~ [(v_{h}^i)_{N}] \leqslant 0 \text {, a.e. on } \Gamma_c^i, i=1,2, \ldots, n-1\right\} \subset V_{h} \text {. }
$$
\end{small}
Since $\mathcal{E}_h^{i}$, $\mathcal{E}_h^{i+1}$ are compatible on $\Gamma_c^{i}$, if the contact conditions are satisfied at all contact nodes, it also hold on $\Gamma_c^{i}$. Then, approximation problem in the finite element space of Problem \ref{prb:3.2:P_0_v} can be expressed as follows:
\begin{problem}[ $P_{h0}^{v}$]\label{prb:3.6:P_h0_v}
Find a displacement $\boldsymbol{u}_{h}\in \mathcal{K}_{h}$ which satisfies:
\begin{small}
\begin{equation}\label{ieq:prb:3.6:var_h}
a(\boldsymbol{u}_{h},\boldsymbol{v}_{h}-\boldsymbol{u}_{h}) + j(\boldsymbol{v}_{h}) - j(\boldsymbol{u}_{h}) \geqslant L(\boldsymbol{v}_{h}-\boldsymbol{u}_{h}), ~\forall \boldsymbol{v}_{h} \in \mathcal{K}_{h}.
\end{equation}
\end{small}
\end{problem}
Similar to the proof of property \ref{prp:3.1}, the Problem \ref{prb:3.6:P_h0_v} can be equivalently represented as $n$ coupled problems in their own $\Omega^i, i=1,\ldots,n$.
However, before this property can be given, discretization of the trace spaces ${H}^{1 / 2}\left(\Gamma_c^{i}\right)$ and discretization of the extension mappings $\chi^{i}_{j}$ need to be introduced.
Let 
\begin{small}
$$
{W}_h^{i,j}=\left\{\varphi_h^i \in\left(C\left(\bar{\Gamma}_c^{i}\right)\right)^d ~\Big|~ \exists \boldsymbol{v}_h^j \in {V}_h^j: \boldsymbol{v}_h^j=\varphi_h^i \text{ on } \Gamma_c^{i} \right\},~j=i \text{ or } i+1, ~  i=1,\ldots,n-1.
$$
\end{small}
From the assumption on $\mathcal{E}_h^i, \mathcal{E}_h^{i+1}$, it follows that $W_{h}^{i,i} = W_{h}^{i,i+1}$. Then, $W_{h}^{i}$ can be defined as $W_{h}^{i} := W_{h}^{i,i+1}$.
Therefore, the discrete extension mappings $\hat{\chi}^{i}_{j}: W_{h}^{i} \to V_{h}^{j}$, $j=i$ or $i+1$ can be defined by:
\begin{small}
\begin{equation}\label{def:3:chi_h}
\left.\begin{array}{rl}
\hat{\chi}^{i}_{j} \varphi^{i}_{h} \in {V}^j_{h}: & a^j\left(\hat{\chi}^i_j \varphi^{i}_{h}, \boldsymbol{v}^j_{h}\right)=0 \quad \forall \boldsymbol{v}_{h}^j \in {V}_{(i-j+3)h}^j(0), \\
& \hat{\chi}^{i}_{j} \varphi^{i}_h =\varphi^{i}_{h} \text { on } \Gamma_c^{i}, ~ j =i \text{ or } i+1
\end{array}\right\},
\end{equation}
\end{small}
where the space $V^{i}_{sh}(0)$ is defined as: 
\begin{small}
$$
V^{i}_{sh}(0)=\left\{\boldsymbol{v}^{i}_{h}\in V^{i}_{h} ~\Big|~ \boldsymbol{v}^{i}_{h}=0,~ \text{ on } \Gamma^{i+s-3}_{c}\right\},~s=2 \text{ or } 3.
$$
\end{small}
Analogously to the norm $\| \cdot \|_{\Gamma_c^{i,j}}$, on the space ${W}_h^{i}$, the following mesh-dependent norms can be defined:
\begin{small}
\begin{equation}\label{def:3:norm_Gamma_h}
\| \varphi^{i}_{h} \|_{\Gamma_{ch}^{i,j}}:=\inf_{\boldsymbol{v}^j_{h} \in V^j_{h}} \left\{\left\|\boldsymbol{v}^j_{h}\right\|_{1, \Omega^j} ~\Big|~ \boldsymbol{v}^j_{h}=\varphi^{i}_{h} \text { on } \Gamma_c^{i}\right\}=\| \hat{\chi}^i_{j} \varphi^{i}_{h} \|_{1, \Omega^j},~ j=i \text{ or } i+1.
\end{equation}
\end{small}

Then, the following spaces need to be defined to facilitate subsequent processes:
\begin{small}
\begin{align*}
K^{i}_{h}(\phi)&=\left\{\boldsymbol{v}^i_{h} \in {V}^i_{h} ~\Big|~ \boldsymbol{v}^i_{h} \cdot \boldsymbol{\beta}^i \leq-\phi \text { on } \Gamma_c^{i}\right\}, \\
V^{i}_{sh}(\boldsymbol{w}^{j}_{h}) &=\left\{\boldsymbol{v}^{i}_{h} \in V_{h}^{i} ~\Big|~ \left\|\gamma^{i}_{k} \boldsymbol{v}^{i}_{h} - \gamma^{j}_{k} \boldsymbol{w}^{j}_{h} \right\|_{\Gamma^{k,k+1}_{ch}}=0,~k=i+s-3 \right\},
\end{align*}
\end{small}
where $\boldsymbol{w}^{j}_{h} \in V^{j}_{h}$, $s=2$ or $3$, $j=i+s-3$ or $i+s-2$ and if $i=1$ and $s=2$ or $i=n$ and $s=3$, we define that $V^{1}_{2h}(\boldsymbol{w}_{h}) = V^{1}_{h}$ and $V^{n}_{3h}(\boldsymbol{w}_{h}) = V^{n}_{h}$ for any $\boldsymbol{w}_{h}$, respectively.

Based on the definitions of $\| \cdot \|_{\Gamma_c^{i,j}}$ and $\| \cdot \|_{\Gamma_{ch}^{i,j}}$, the following inequality is obvious:
\begin{small}
\begin{equation}\label{ieq:3:norm_Gamma_h_1}
\| \varphi^{i}_{h} \|_{\Gamma_{ch}^{i,j}} \geqslant \| \varphi^{i}_{h} \|_{\Gamma_c^{i,j}}, ~\forall \varphi^{i}_{h} \in W^{i}_{h},~ \forall h>0, ~i=1,\ldots,n-1.
\end{equation}
\end{small}
In order to construct an equivalence relationship between $\| \cdot \|_{\Gamma_c^{i,j}}$ and $\| \cdot \|_{\Gamma_{ch}^{i,j}}$, the following assumptions are introduced:

\begin{assumption}\label{assum:3.1}
For $i=1,\ldots,n-1$, there exists a constant $c^{i}>0$ independent of $h$ and such that
\begin{small}
\begin{equation}\label{ieq:assum:3:norm_Gamma_h_2}
\| \varphi^{i}_{h} \|_{\Gamma_{ch}^{i,j}} \leqslant c^{i} \| \varphi^{i}_{h} \|_{\Gamma_c^{i,j}}, ~\forall \varphi^{i}_{h} \in W^{i}_{h},~ \forall h>0, ~i=1,\ldots,n-1.
\end{equation}
\end{small}
\end{assumption}

\begin{remark}
Under some regularity conditions, Assumption \ref{assum:3.1} can be verified. Referring to the methods in Refs.\cite{3haslinger2014domain,bjorstad1986iterative}, the following verification can be performed. For $\varphi^i_h \in W^{i}_h$, we assume that the regularity condition $\chi_j^i \varphi_h^i \in H^{1+\epsilon}(\Omega^{j})^d$ holds, where $\epsilon \in (0,\frac{1}{2})$.
According to the definition of norm $\|\cdot\|_{\Gamma_{ch}^{i, j}}$, we have
$$
\| \varphi^{i}_h \|_{\Gamma_{ch}^{i, j}} = \left\|\hat{\chi}_j^i \varphi_h^i\right\|_{1, \Omega^j} \leqslant \left\| \chi_j^i \varphi_h^i\right\|_{1, \Omega^j} + \left\| \hat{\chi}_j^i \varphi_h^i - \chi_j^i \varphi_h^i\right\|_{1, \Omega^j},~ \forall \varphi^i_h \in W^i_h.
$$
Since $\hat{\chi}_j^i \varphi_h^i$ is the Galerkin approximation of $\chi_j^i \varphi_h^i$, based on the interpolation theorem \cite{chouly2023finite}, it can be verified that 
$$
\left\| \hat{\chi}_j^i \varphi_h^i - \chi_j^i \varphi_h^i\right\|_{1, \Omega^j} = c^i_1 h^{\epsilon} \|\chi_j^i \varphi_h^i\|_{1+\epsilon, \Omega^j}, ~ c^i_1>0.
$$
Then, taking into account the trace theorem in Sobolev space \cite{3adams2003sobolev}, it can be verified that $\|\chi_j^i \varphi_h^i\|_{1+\epsilon, \Omega^j} \leqslant c^i_2 \|\varphi_h^i\|_{\frac{1}{2} + \epsilon, \Gamma^i_c}$, where the constant $c^i_2>0$ and is independent of $\varphi^i_h$. Therefore, if the following inequality holds:
$$
\|\varphi_h^i\|_{\frac{1}{2} + \epsilon, \Gamma^i_c} \leqslant c^i_3 h^{-\epsilon} \|\varphi_h^i\|_{\frac{1}{2}, \Gamma^i_c}, \forall \varphi_h^i\in W^i_h,
$$
where $c^i_3 > 0$ and is independent of $h$, it can be deduced that
$$
\| \varphi^{i}_h \|_{\Gamma_{ch}^{i, j}} \leqslant \left\| \chi_j^i \varphi_h^i\right\|_{1, \Omega^j} + c^i_1c^i_2c^i_3\|\varphi_h^i\|_{\frac{1}{2}, \Gamma^i_c} = c^i \| \varphi^{i}_h \|_{\Gamma_{c}^{i, j}}.
$$
In an appropriate finite element space $W^i_h$, it can be verified that the above conditions hold.
\end{remark}

Therefore, under the Assumption \ref{assum:3.1}, $\| \cdot \|_{\Gamma_c^{i,j}}$ and $\| \cdot \|_{\Gamma_{ch}^{i,j}}$ are equivalent. Then, based on the equivalence property of norm $\| \cdot \|_{\Gamma_c^{i,i}}$ and norm $\| \cdot \|_{\Gamma_c^{i,i+1}}$, it can be obtained that norm $\| \cdot \|_{\Gamma_{ch}^{i,i}}$ and norm $\| \cdot \|_{\Gamma_{ch}^{i,i+1}}$ are also equivalent, that is, $\exists c^{i}_{1} , c^{i}_{2}>0$ independent of mesh size $h$ such that:
\begin{small}
$$
c^{i}_{1} \|\varphi^i_h\|_{\Gamma_{ch}^{i,i+1}} \leqslant \|\varphi^i_{h}\|_{\Gamma_{ch}^{i,i}} \leqslant c^{i}_{2} \|\varphi^i_{h}\|_{\Gamma_{ch}^{i,i+1}},~ \forall \varphi_{h}^i \in W_{h}^{i}, ~ i=1\ldots,n-1.
$$
\end{small}

Similar to property \ref{prp:3.1}, Problem \ref{prb:3.6:P_h0_v} is equipped with the following property:

\begin{proposition}\label{prp:3.7}
A function $\boldsymbol{u}_{h}=(\boldsymbol{u}^1_{h},\ldots,\boldsymbol{u}^{n}_{h}) \in \mathcal{K}_{h}$ is a solution of $P_{h0}^{v}$ if and only if $\boldsymbol{u}^{i}_{h}\in K^{i}_{h}(\boldsymbol{u}^{i+1}_{h} \cdot\boldsymbol{\alpha}^{i+1})$, $i=1,\ldots,n-1$ and $\boldsymbol{u}^{n}_{h}\in V^{n}_{h}$ solve the following inequalities:
\begin{small}
\begin{equation}\label{ieq:prp:3.7:1}
\left\{
\begin{aligned}
& a^{1}(\boldsymbol{u}^{1}_{h}, \boldsymbol{v}^{1}_{h}-\boldsymbol{u}^{1}_{h}) + j^{1}( \boldsymbol{v}^{1}_{h}, \boldsymbol{u}^{2}_{h}) -j^{1}( \boldsymbol{u}^{1}_{h}, \boldsymbol{u}^{2}_{h}) \geqslant L^{1}(\boldsymbol{v}^{1}_{h} -\boldsymbol{u}^{1}_{h}),\\
& a^{i}\left(\boldsymbol{u}^{i}_{h}, \boldsymbol{v}^{i}_{h} -\boldsymbol{u}^{i}_{h}\right) + j^{i}\left( \boldsymbol{v}^{i}_{h}, \boldsymbol{u}^{i+1}_{h}\right) - j^{i}\left( \boldsymbol{u}^{i}_{h}, \boldsymbol{u}^{i+1}_{h}\right) \geqslant  L^{i}\left(\boldsymbol{v}^{i}_{h} -\boldsymbol{u}^{i}_{h}\right)\\ 
& +L^{i-1}\left(\hat{\chi}^{i-1}_{i-1}\gamma^{i}_{i-1} \left(\boldsymbol{v}^{i}_{h} -\boldsymbol{u}^{i}_{h} \right)\right) - a^{i-1}\left(\boldsymbol{u}^{i-1}_{h}, \hat{\chi}^{i-1}_{i-1}\gamma^{i}_{i-1}\left(\boldsymbol{v}^{i}_{h} -\boldsymbol{u}^{i}_{h} \right)\right), \\
& a^{n}\left(\boldsymbol{u}^{n}_{h}, \boldsymbol{v}^{n}_{h}\right) = L^{n}\left(\boldsymbol{v}^{n}_{h}\right) + L^{n-1}\left(\hat{\chi}^{n-1}_{n-1}\gamma^{n}_{n-1} \boldsymbol{v}^{n}_{h}\right) - a^{n-1}\left(\boldsymbol{u}^{n-1}_{h}, \hat{\chi}^{n-1}_{n-1}\gamma^{n}_{n-1}\boldsymbol{v}^{n}_{h} \right). 
\end{aligned}
\right. 
\end{equation}
\end{small}
where $\forall \boldsymbol{v}^{i}_{h} \in K^{i}(\boldsymbol{u}^{i+1}_{h} \cdot\boldsymbol{\alpha}^{i+1}),~i=1,\ldots,n-1$, $\forall \boldsymbol{v}^{n}_{h} \in V^{n}_{h}$.
\end{proposition}

It should be noticed that to satisfy the condition on $\Gamma_c^{i}$ in $K^{i}(\boldsymbol{u}^{i+1}_{h} \cdot\boldsymbol{\alpha}^{i+1})$ it is sufficient to ensure that it holds at all finite element nodes. In order to give the discrete form of Algorithm \ref{algorithm:3.1} (LD), the discretization problems of \ref{prb:3.4:d1}, \ref{prb:3.4:di}, \ref{prb:3.4:dn}, \ref{prb:3.5:pi} and \ref{prb:3.5:qi} $(i=1,\ldots,n-1)$ should be given.

Then, for given $\boldsymbol{\lambda}_{h} = \left(\lambda_{h}^1, \ldots, \lambda_{h}^{n-1}\right) \in W_{h}^{1} \times \cdots \times W_{h}^{n-1}$, the problems $P^{i}_{dh}(\boldsymbol{\lambda}_{h})$ can be defined by:
\begin{problem}[ $P_{dh}^i(\boldsymbol{\lambda}_{h})$]\label{prb:3.7:dh}
\begin{small}
\begin{align}
&\left\{
\begin{aligned}
&\text{Find } \boldsymbol{u}^1_{h}:=\boldsymbol{u}^1_{h}(\boldsymbol{\lambda}_{h}) \in K^{1}_{h}(\lambda^{1}_{h} \cdot \boldsymbol{\alpha}^{2}) \text{ such that:}\\
&a^1\left(\boldsymbol{u}^1_{h}, \boldsymbol{v}^1_{h} -\boldsymbol{u}^1_{h} \right) + j^1\left( \boldsymbol{v}^1_{h}, \lambda^1_{h} \right)-j^1\left( \boldsymbol{u}^1_{h}, \lambda^1_{h}\right) \geqslant L^1\left(\boldsymbol{v}^1_{h}-\boldsymbol{u}^1_{h} \right)\\
& \forall \boldsymbol{v}^1_{h} \in K^{1}_{h}(\lambda^{1}_{h} \cdot\boldsymbol{\alpha}^{2})
\end{aligned}
\right.  \tag{\text{$P^{1}_{dh}(\boldsymbol{\lambda}_{h})$}}
\label{prb:3.7:dh1}\\
&\left\{
\begin{aligned}
&\text{Find } \boldsymbol{u}^i_{h} :=\boldsymbol{u}^i_{h}(\boldsymbol{\lambda}_{h}) \in K^{i}_{h} (\lambda^{i}_{h} \cdot \boldsymbol{\alpha}^{i+1}), ~i=2,\ldots,n-1 \text{ such that:}\\
&a^i\left(\boldsymbol{u}^i_{h}, \boldsymbol{v}^i_{h}-\boldsymbol{u}^i_{h}\right) + j^i\left( \boldsymbol{v}^i_{h}, \lambda^i_{h} \right) - j^i\left( \boldsymbol{u}^i_{h}, \lambda^i_{h} \right) \geqslant L^i\left(\boldsymbol{v}^i_{h} -\boldsymbol{u}^i_{h} \right)\\
& \forall \boldsymbol{v}^i_{h} \in K^{i}_{h}(\lambda^{i}_{h} \cdot \boldsymbol{\alpha}^{i+1}) \cap V^{i}_{2h}(\boldsymbol{u}^i_{h}), ~\boldsymbol{u}^i_{h} = \lambda^{i-1}_{h} \text{ on } \Gamma^{i-1}_{c}
\end{aligned}
\right.  \tag{\text{$P^{i}_{dh}(\boldsymbol{\lambda}_h)$}}
\label{prb:3.7:dhi}\\
&\left\{
\begin{aligned}
&\text{Find } \boldsymbol{u}^n_{h} := \boldsymbol{u}^n_{h}(\boldsymbol{\lambda}_{h}) \in V^{n}_{h} \text{ such that:}\\
&a^n\left(\boldsymbol{u}^n_{h}, \boldsymbol{v}^n_{h}\right) = L^n\left(\boldsymbol{v}^n_{h} \right)\\
& \forall \boldsymbol{v}^n_{h} \in V^{n}_{2h}(0), ~ \boldsymbol{u}^n_{h} = \lambda^{n-1}_{h} \text{ on } \Gamma^{n-1}_{c}
\end{aligned}
\right.  \tag{\text{$P^{n}_{dh}(\boldsymbol{\lambda}_h)$}} \label{prb:3.7:dhn}
\end{align}
\end{small}
\end{problem}

To mach the discrete contact stress on $\Gamma^{i}_{c}$, the auxiliary problems $P_{ph}^{i+1}(\boldsymbol{\lambda}_{h})$ and $P_{qh}^{i}(\boldsymbol{\lambda}_h)$ can be defined by: 
\begin{problem}[$P_{ph}^{i+1}(\boldsymbol{\lambda})$ and $P_{qh}^{i}(\boldsymbol{\lambda})$]\label{prb:3.8:P_pq_h}
Find $\boldsymbol{p}^{i+1}_{h}\in V^{i+1}_{3h}(0)$, $\boldsymbol{q}^{i}_{h} \in V^{i}_{2h}(0) (i=1,\ldots,n-1)$  such that:
\begin{small}
\begin{align}
a^{i+1}\left(\boldsymbol{p}^{i+1}_{h}, \boldsymbol{v}^{i+1}_{h}\right) =& \frac{1}{2} \Big(a^{i+1}\left(\boldsymbol{u}^{i+1}_{h}, \boldsymbol{v}^{i+1}_{h}\right)- L^{i+1}\left(\boldsymbol{v}^{i+1}_{h}\right)+ a^{i}\left(\boldsymbol{u}^{i}_{h}, \hat{\chi}_{i}^{i} \gamma^{i+1}_{i}\boldsymbol{v}^{i+1}_{h}\right) \nonumber\\ 
&- L^{i}\left(\hat{\chi}_{i}^{i} \gamma^{i+1}_{i}\boldsymbol{v}^{i+1}_{h}\right) \Big),~ \forall \boldsymbol{v}^{i+1}_{h} \in V^{i+1}_{3h}(0) \tag{$P^{i+1}_{ph}(\boldsymbol{\lambda}_{h})$} \label{prb:3.8:phi} \\
a^{i}\left(\boldsymbol{q}^{i}_{h}, \boldsymbol{v}^{i}_{h}\right) =& \frac{1}{2} \Big(a^{i+1}\left(\boldsymbol{u}^{i+1}_{h}, \hat{\chi}^{i}_{i+1} \gamma^{i}_{i}\boldsymbol{v}^{i}_{h}\right)- L^{i+1}\left(\hat{\chi}^{i}_{i+1} \gamma^{i}_{i} \boldsymbol{v}^{i}_{h} \right) + a^{i}\left(\boldsymbol{u}^{i}_{h}, \boldsymbol{v}^{i}_{h}\right) \nonumber\\
&- L^{i}\left(\boldsymbol{v}^{i}_{h}\right) \Big),~ \forall \boldsymbol{v}^{i}_{h} \in V^{i}_{2h}(0) \tag{\text{$P^{i}_{qh}(\boldsymbol{\lambda})$}} 
\label{prb:3.8:qhi}
\end{align}
\end{small}
where $\boldsymbol{u}^{i}_{h}:= \boldsymbol{u}^{i}_{h}(\boldsymbol{\lambda}_{h})$ are the solutions of $P_{dh}^{i}(\boldsymbol{\lambda}_{h})$.
\end{problem}

The discrete layer decomposition algorithm can be defined as Algorithm \ref{algorithm:3.2}. Thus, the definition of convergence in the discrete case is summarized as follows:

\begin{algorithm}[!ht]
\caption{Discrete Layer Decomposition Algorithm}
\label{algorithm:3.2}
Let $\boldsymbol{\lambda}_{0h} = (\lambda^{1}_{0h},\ldots,\lambda^{n-1}_{0h}) \in W_{h}^{1} \times \cdots \times W_{h}^{n-1}$ and $\theta>0$ be given.

\textbf{Output:} $\boldsymbol{u}^{i}_{kh}$, $\boldsymbol{p}^{i}_{kh}$, $\boldsymbol{q}^{i}_{kh}$, $\boldsymbol{\lambda}_{kh}$

\textbf{For} $k\leq N$ \textbf{do}

~~~~$
\left\{
\begin{aligned}
& \boldsymbol{u}^{1}_{kh}(\lambda_{(k-1)h}) \in K^1_{h}\left(\lambda^1_{(k-1)h} \cdot \boldsymbol{\alpha}^2\right) \text{ solves } P^{1}_{dh}(\lambda_{(k-1)h})\\
& \boldsymbol{u}^{i}_{kh}(\lambda_{(k-1)h}) \in K^i_{h}\left(\lambda^i_{(k-1)h} \cdot \boldsymbol{\alpha}^{i+1} \right) \text{ solves } P^{i}_{dh}(\lambda_{(k-1)h}), ~i=2,\ldots,n-1\\
& \boldsymbol{u}^{n}_{kh}(\lambda_{(k-1)h}) \in V^{n}_{h} \text{ solves } P^{n}_{dh}(\lambda_{(k-1)h})
\end{aligned}
\right.
$

~~~~$
\left\{
\begin{aligned}
& \boldsymbol{p}^{i+1}_{kh} \in V^{i+1}_h \text{ solves } P^{i+1}_{ph}(\lambda_{(k-1)h})\\
& \boldsymbol{q}^{i}_{kh} \in V^{i}_{h} \text{ solves } P^{i}_{qh}(\lambda_{(k-1)h})
\end{aligned}
\right.
$

~~~~$\lambda_{kh}^i=\lambda_{(k-1)h}^{i} -\theta\left(\gamma^{i+1}_{i}\boldsymbol{p}^{i+1}_{kh} + \gamma^{i}_{i}\boldsymbol{q}^{i}_{kh}\right)$ on $\Gamma_c^i$.

\textbf{end}
\end{algorithm}

\begin{theorem}\label{thm:3.2}
If Assumption \ref{assum:3.1} hold, there exist $0<\theta_*<4/(3c_{t}^{\min})^2$ and functions $\boldsymbol{\lambda}_{*h} = (\lambda^{1}_{*h},\ldots,\lambda^{n-1}_{*h}) \in W^{1}_{h} \times \cdots \times W^{n-1}_{h}$, $\boldsymbol{u}^{i}_{*h}$, $\boldsymbol{p}^{i}_{*h}$, $\boldsymbol{q}^{i}_{*h}$, $i=1,\ldots,n$ such that $\forall \theta \in\left(0, \theta_*\right)$,
\begin{small}
\begin{align*}
&~\lambda_{kh}^i \rightarrow \lambda_{*h}^i \text { in } W^i_h, && i=1, \ldots, n-1\\
&\left.\begin{array}{l}
\boldsymbol{p}_{kh}^{i+1} \rightarrow \boldsymbol{p}_{*h}^{i+1} \text { in } V_{h}^{i+1}\\
\boldsymbol{q}_{kh}^{i} \rightarrow \boldsymbol{p}_{*h}^i \text { in } V_{h}^i
\end{array}\right\}, && i=1, \ldots, n-1 \\
&~\boldsymbol{u}_{kh}^i \rightarrow \boldsymbol{u}_{*h}^i \text { in } V_{h}^i, && i=1 \ldots, n
\end{align*}
\end{small}
as $k\to \infty$, where the sequence $\left\{\left(\boldsymbol{u}_{kh}^i, \boldsymbol{p}_{kh}^i, \boldsymbol{q}_{kh}^i, \lambda_{kh}^i\right)\right\}$ is generated by Algorithm \ref{algorithm:3.2}. In addition, the $\boldsymbol{u}_{*h}=\left(\boldsymbol{u}^1_{*h}, \ldots, \boldsymbol{u}^n_{*h}\right) \in \mathcal{K}_{h}$ solves $P^{v}_{h0}$.
\end{theorem}

\begin{proof}

First, in the discrete case, under Assumption \ref{assum:3.1} the solutions of problems \ref{prb:3.7:dh1}, \ref{prb:3.7:dhi}, \ref{prb:3.7:dhn} $(i=2,\ldots,n-1)$, \ref{prb:3.8:phi} $(i=1,\ldots,n-1)$ and \ref{prb:3.8:qhi} $(i=1,\ldots,n-1)$ still have propositions similar to propositions \ref{prp:3.2} - \ref{prp:3.5}, and their concrete proof is also similar to the proof process in the continuous case. So these discrete properties will be used directly and the proof will be omitted.

Then, similar to Proposition \ref{prp:3.6}, it can be proved that if $\boldsymbol{\lambda}_{kh} = \{ \lambda^{1}_{kh}, \ldots, \lambda^{n-1}_{kh} \}$ converges to $\boldsymbol{\lambda}_{*h} = \{ \lambda^{1}_{*h}, \ldots, \lambda^{n-1}_{*h} \}$, the sequence $\{ \boldsymbol{p}_{k h}^{i+1}, \boldsymbol{q}_{k h}^i, \boldsymbol{u}_{k h}^i \}$ generated by the Algorithm \ref{algorithm:3.2} (LD) converges to $\{ \boldsymbol{p}_{*h}^{i+1}, \boldsymbol{q}_{*h}^i, \boldsymbol{u}_{*h}^i \}$ and $\boldsymbol{u}_{*h}= \{\boldsymbol{u}_{*h}^1, \ldots, \boldsymbol{u}_{*h}^n \}$ is the solution of problem $P^{v}_{0h}$. Therefore, by proving that the sequence $\boldsymbol{\lambda}_{kh}$ converges to the sequence $\boldsymbol{\lambda}_{*h}$, the conclusion of Theorem \ref{thm:3.2} can be proved.

Finally, the convergence of $\boldsymbol{\lambda}_{h} = \{ \lambda^{1}_{h}, \ldots, \lambda^{n-1}_{h} \}$ can be proved by constructing the compressed mapping $T_{\theta h}$ according to the proof of Theorem \ref{thm:3.1}. It is worth noting that under Assumption \ref{assum:3.1}, the equivalence property of the norm guarantees that the convergence of the algorithm is independent of the mesh size $h$.
\end{proof}

It is worth emphasizing that although the convergence of the discrete Algorithm \ref{algorithm:3.2} LD dose not depends on the mesh size $h$, $h$ determines the accuracy at which the solution $\boldsymbol{u}_{*h}$ of the discrete problem $P^{v}_{0h}$ converges to the solution $\boldsymbol{u}_{*}$ of the problem $P^{v}_{0}$. However, in the actual implementation of the algorithm, the iteration cannot proceed infinitely, so the tolerance error $tol>0$ needs to be selected so that the algorithm stops when $\|\boldsymbol{\lambda}_{kh}- \boldsymbol{\lambda}_{k-1h}\|/\|\boldsymbol{\lambda}_{kh}\|<tol$. Moreover, when $tol$ is sufficiently small, the error between the finite element solution and the exact solution will be determined by the mesh diameter $h$.

\begin{remark}
In the mesh definition of the discrete problem, the meshes $\mathcal{E}^{i}_h$ and $\mathcal{E}^{i+1}_h$ on the contact zone $\Gamma^i_c$ are required to match. This condition can ensure that the non-penetration condition is strictly established, that is, $\mathcal{K}_h\subset \mathcal{K}$, but it also limits the applicability of this algorithm. In Refs.\cite{belgacem1998mortar, hild2000numerical, drouet2017accurate, fernandez2003numerical, wohlmuth2012abstract}, the mortar method was proposed to deal with the incompatibility of mesh nodes in the contact zone $\Gamma^{i}_c$, that is, by using projection operators such as local Lagrangian interpolation, mortar projection, and local average contact methods to ensure that the non-penetration condition on the mortar boundary is established under integration. But the approximation method in this discrete space leads to  $\mathcal{K}_h\not\subset \mathcal{K}$. Since in the layer decomposition algorithm, the original problem is decomposed into multiple unilateral contact sub-problems, mortar technology can be used to solve the calculation problems caused by incompatibility of contact boundary mesh nodes on these sub-problems. 
However, it is worth emphasizing that the overall convergence properties of the algorithm have not been explored after using the mortar method in multi-layer elastic contact systems and the Neumann-Neumann domain decomposition method. Therefore, even if it is feasible in terms of technical demonstration, theoretical analysis and numerical experimental verification are necessary, and this is also our subsequent research direction.
\end{remark}

\begin{figure}[!t]
\centering
\begin{minipage}{0.3\linewidth}
\centering
\includegraphics[width=0.9\linewidth]{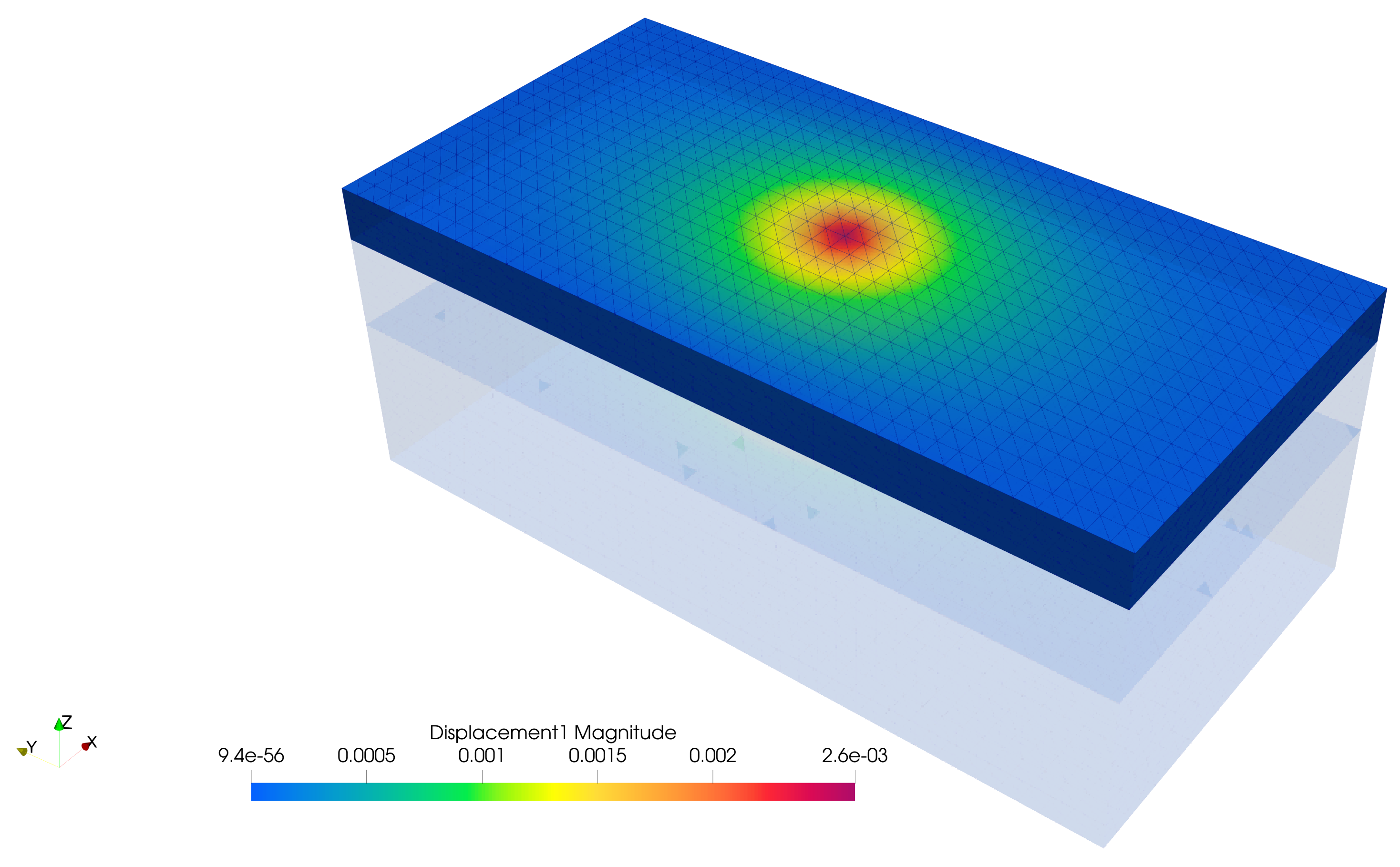}
\caption{Displacement field of the $\Omega^1$ under $\boldsymbol{f}_{2}^{1}$.}
\label{3:fig3}
\end{minipage}
\begin{minipage}{0.3\linewidth}
\centering
\includegraphics[width=0.9\linewidth]{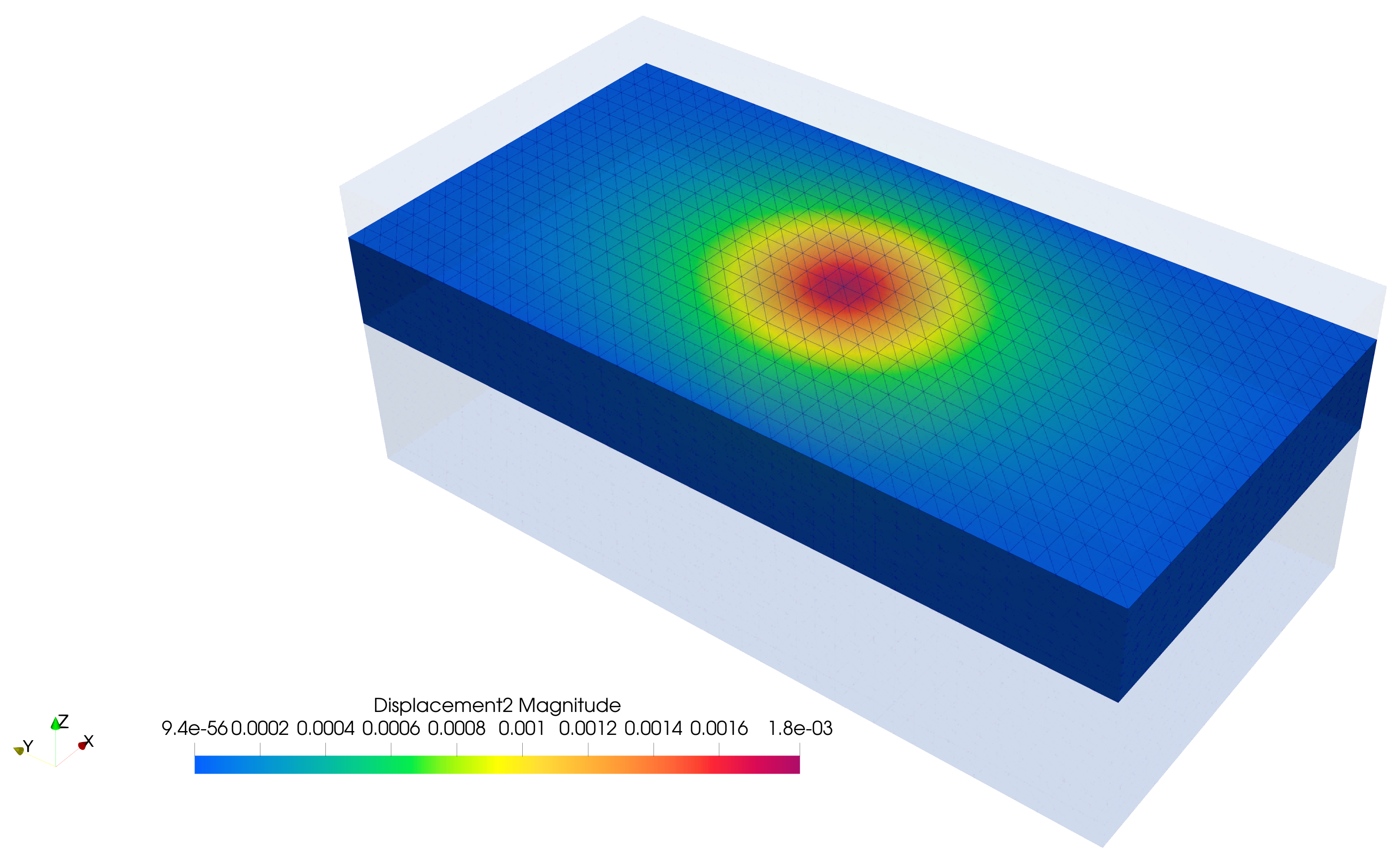}
\caption{Displacement field of the $\Omega^2$ under $\boldsymbol{f}_{2}^{1}$.}
\label{3:fig4}
\end{minipage}
\begin{minipage}{0.3\linewidth}
\centering
\includegraphics[width=0.9\linewidth]{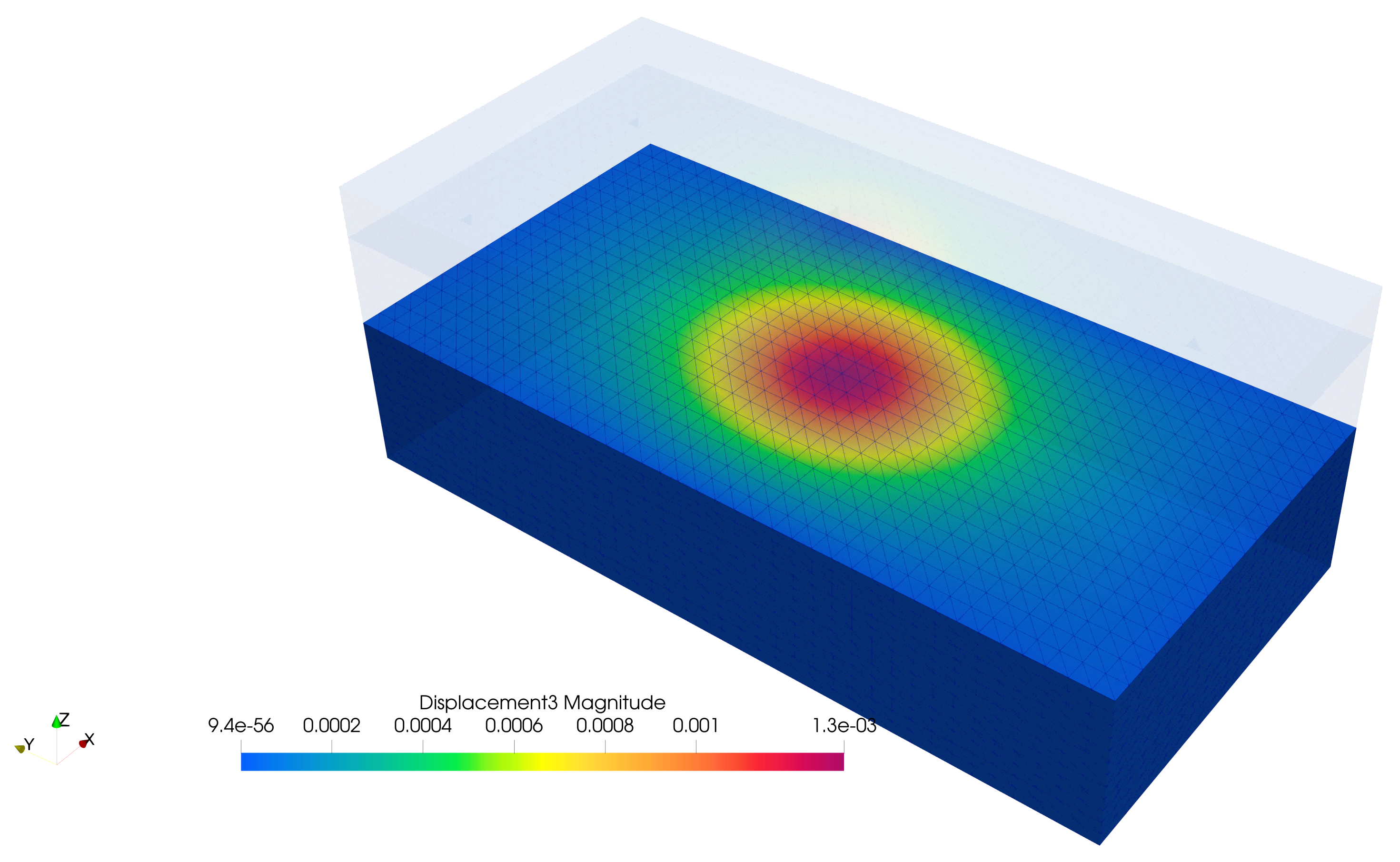}
\caption{Displacement field of the $\Omega^3$ under $\boldsymbol{f}_{2}^{1}$.}
\label{3:fig5}
\end{minipage}

\begin{minipage}{0.3\linewidth}
\centering
\includegraphics[width=0.9\linewidth]{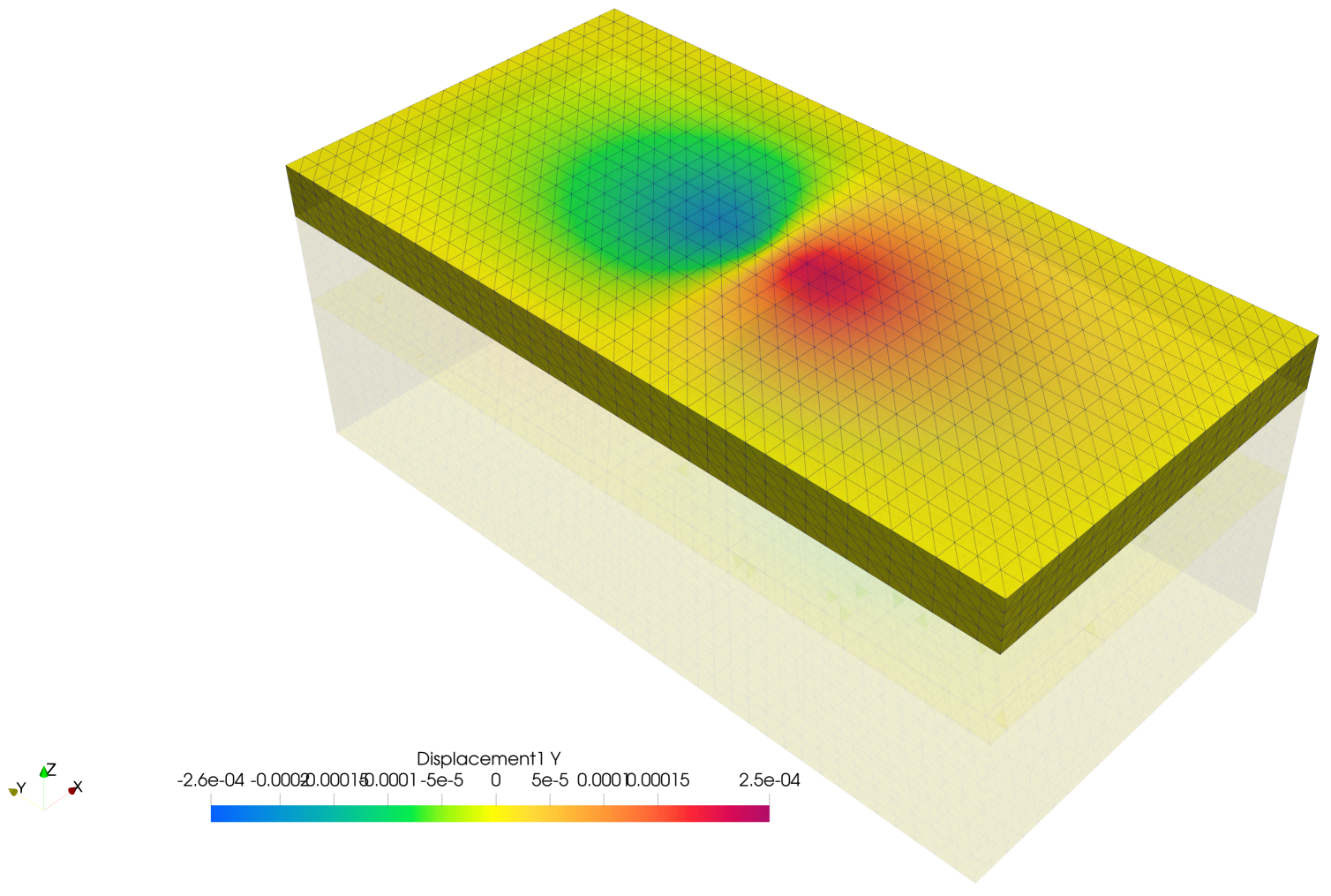}
\caption{Displacement field along the Y-axis of the $\Omega^1$ under $\boldsymbol{f}_{2}^{1}$.}
\label{3:fig6}
\end{minipage}
\begin{minipage}{0.3\linewidth}
\centering
\includegraphics[width=0.9\linewidth]{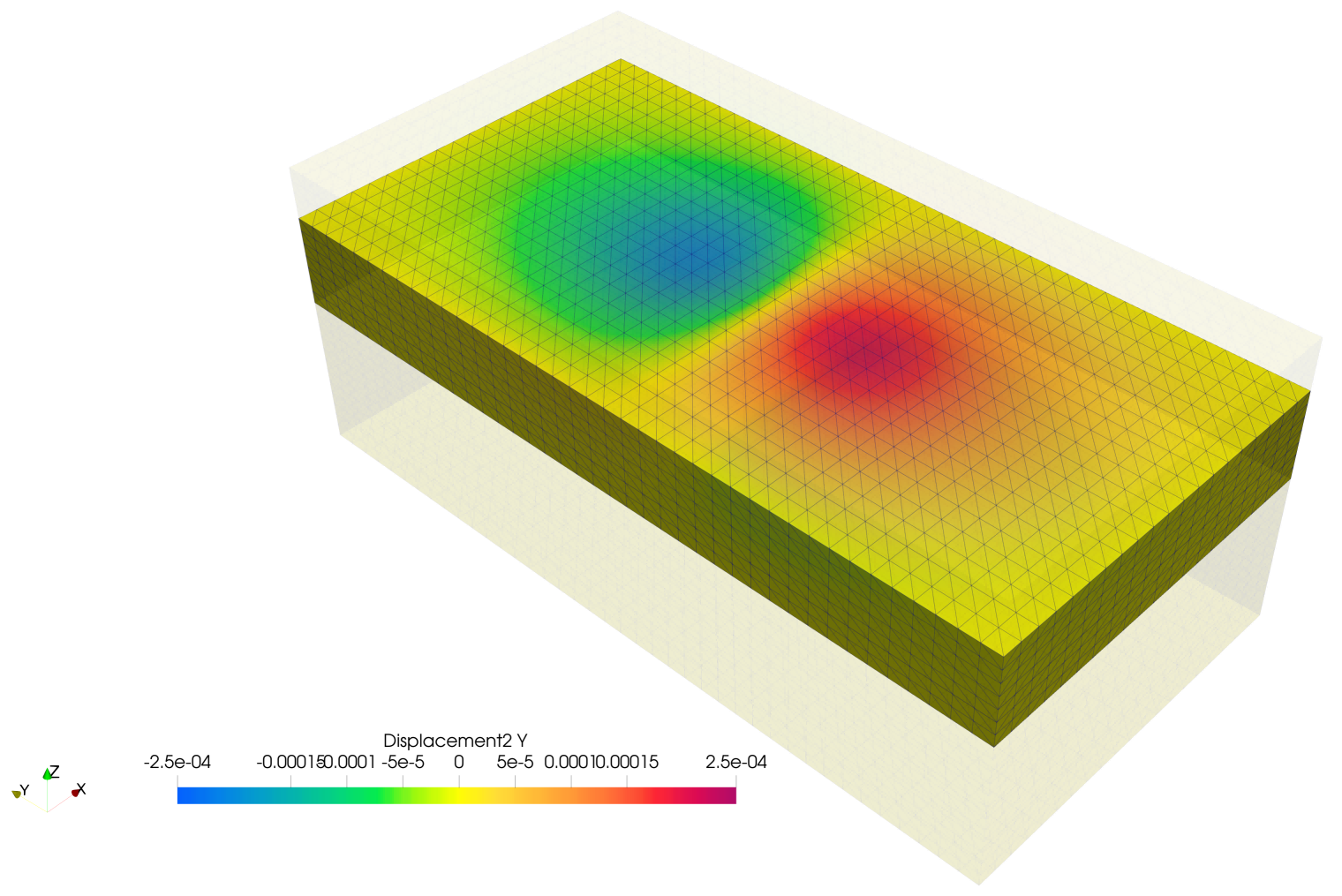}
\caption{Displacement field along the Y-axis of the $\Omega^2$ under $\boldsymbol{f}_{2}^{1}$.}
\label{3:fig7}
\end{minipage}
\begin{minipage}{0.3\linewidth}
\centering
\includegraphics[width=0.9\linewidth]{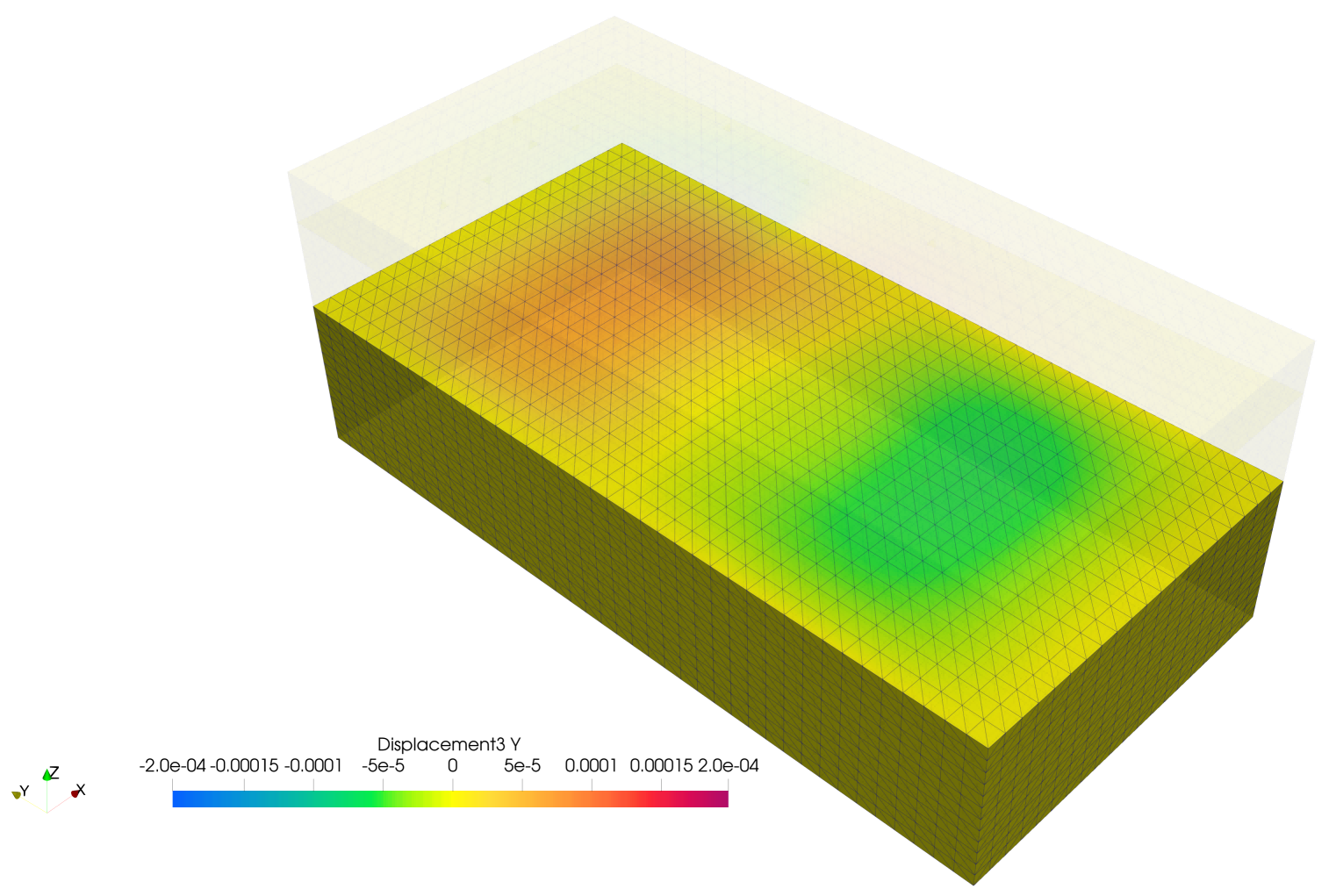}
\caption{Displacement field along the Y-axis of the $\Omega^3$ under $\boldsymbol{f}_{2}^{1}$.}
\label{3:fig8}
\end{minipage}
\end{figure}

\begin{figure}[!t]
\centering
\begin{minipage}{0.3\linewidth}
\centering
\includegraphics[width=0.9\linewidth]{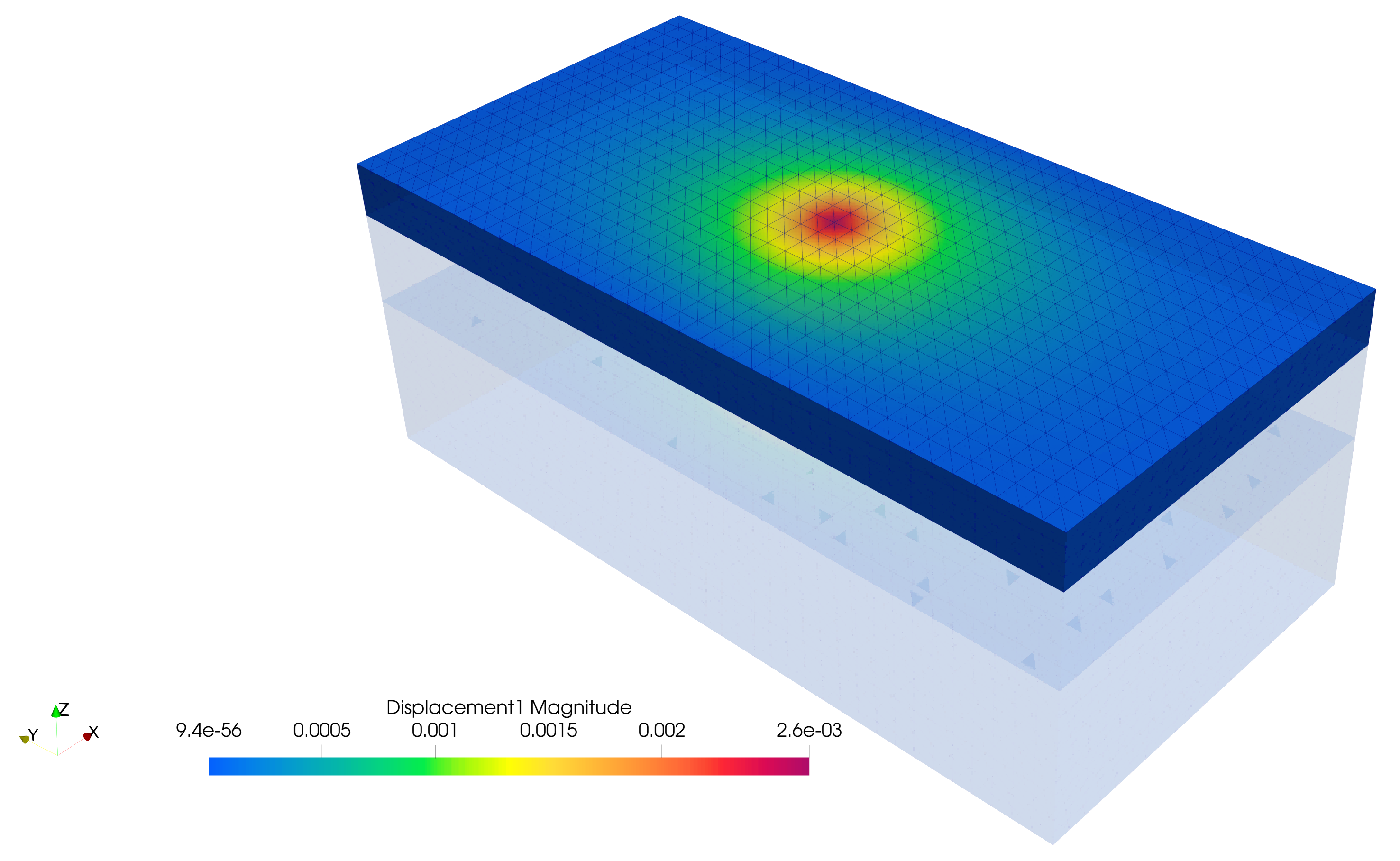}
\caption{Displacement field of the $\Omega^1$ under $\boldsymbol{f}_{2}^{2}$.}
\label{3:fig11}
\end{minipage}
\begin{minipage}{0.3\linewidth}
\centering
\includegraphics[width=0.9\linewidth]{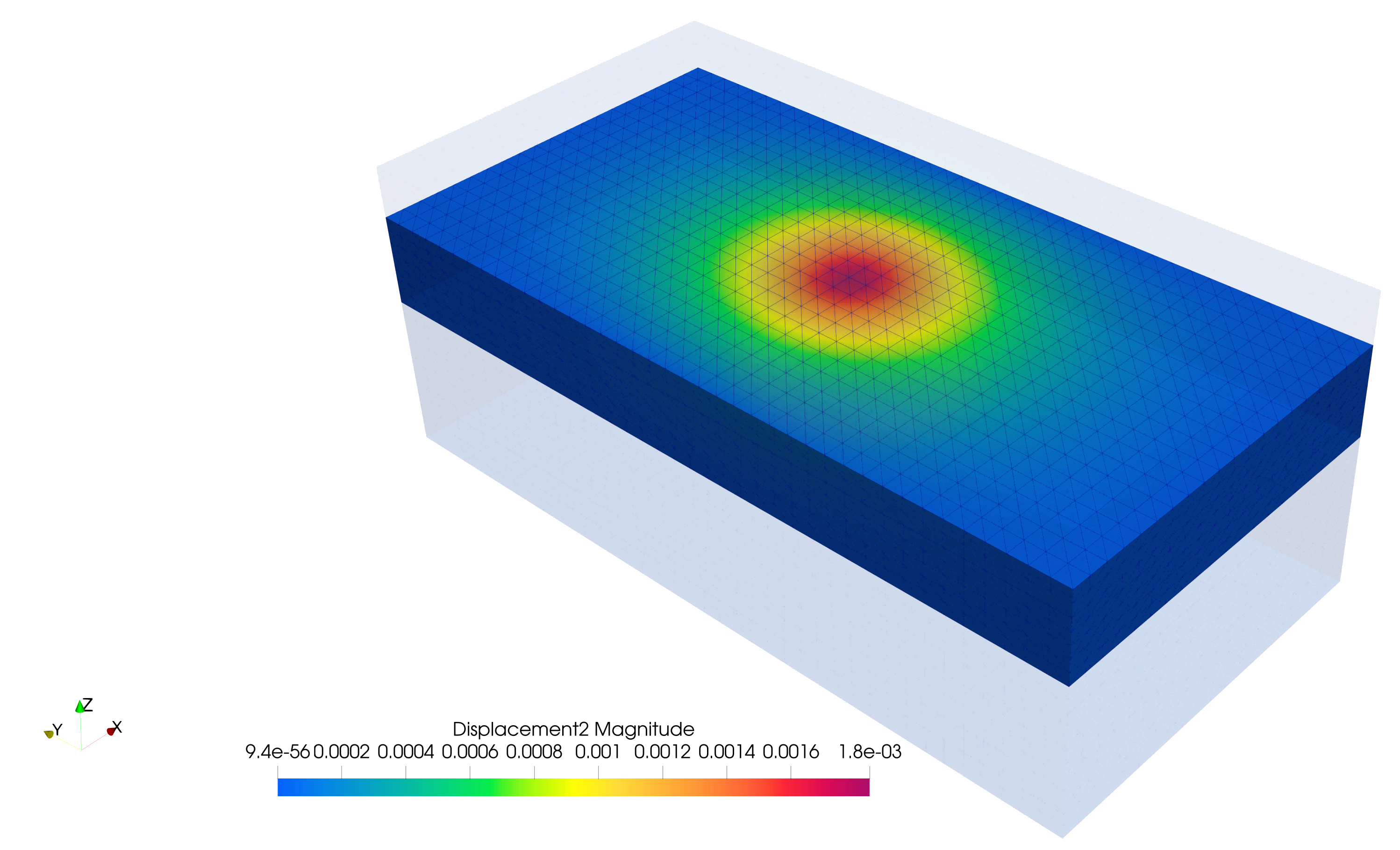}
\caption{Displacement field of the $\Omega^2$ under $\boldsymbol{f}_{2}^{2}$.}
\label{3:fig12}
\end{minipage}
\begin{minipage}{0.3\linewidth}
\centering
\includegraphics[width=0.9\linewidth]{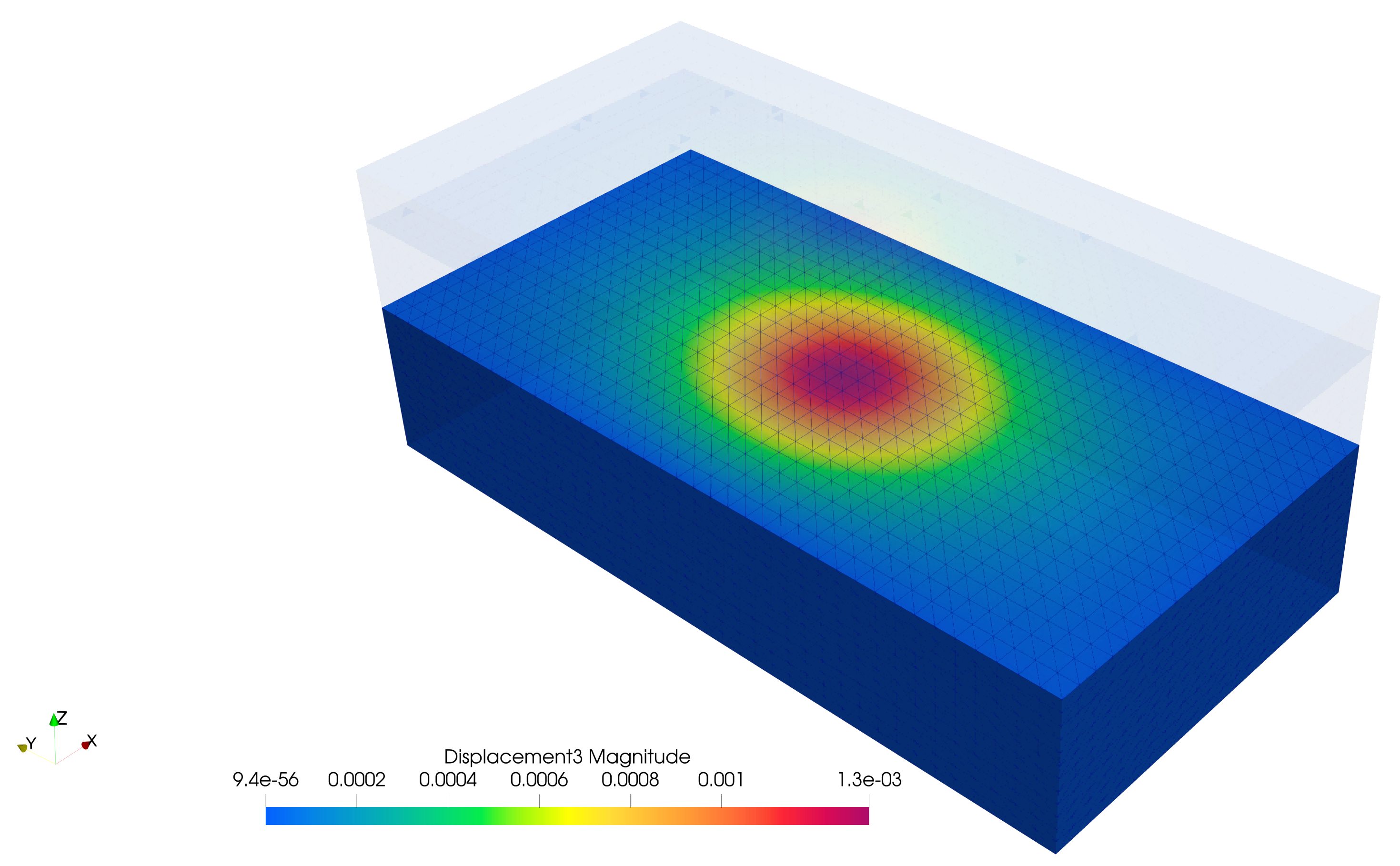}
\caption{Displacement field of the $\Omega^3$ under $\boldsymbol{f}_{2}^{2}$.}
\label{3:fig13}
\end{minipage}

\begin{minipage}{0.3\linewidth}
\centering
\includegraphics[width=0.9\linewidth]{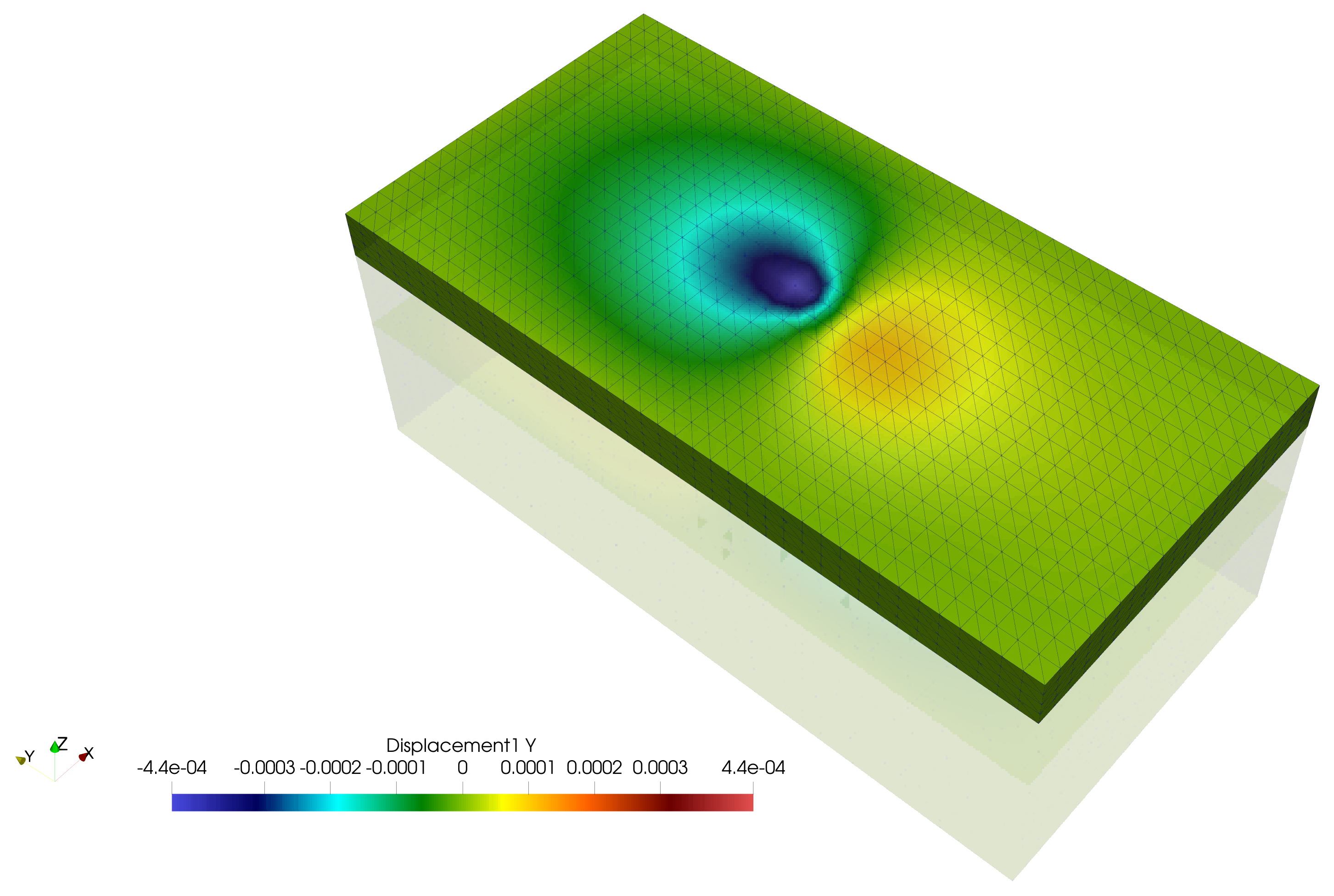}
\caption{Displacement field along the Y-axis of the $\Omega^1$ under $\boldsymbol{f}_{2}^{2}$.}
\label{3:fig14}
\end{minipage}
\begin{minipage}{0.3\linewidth}
\centering
\includegraphics[width=0.9\linewidth]{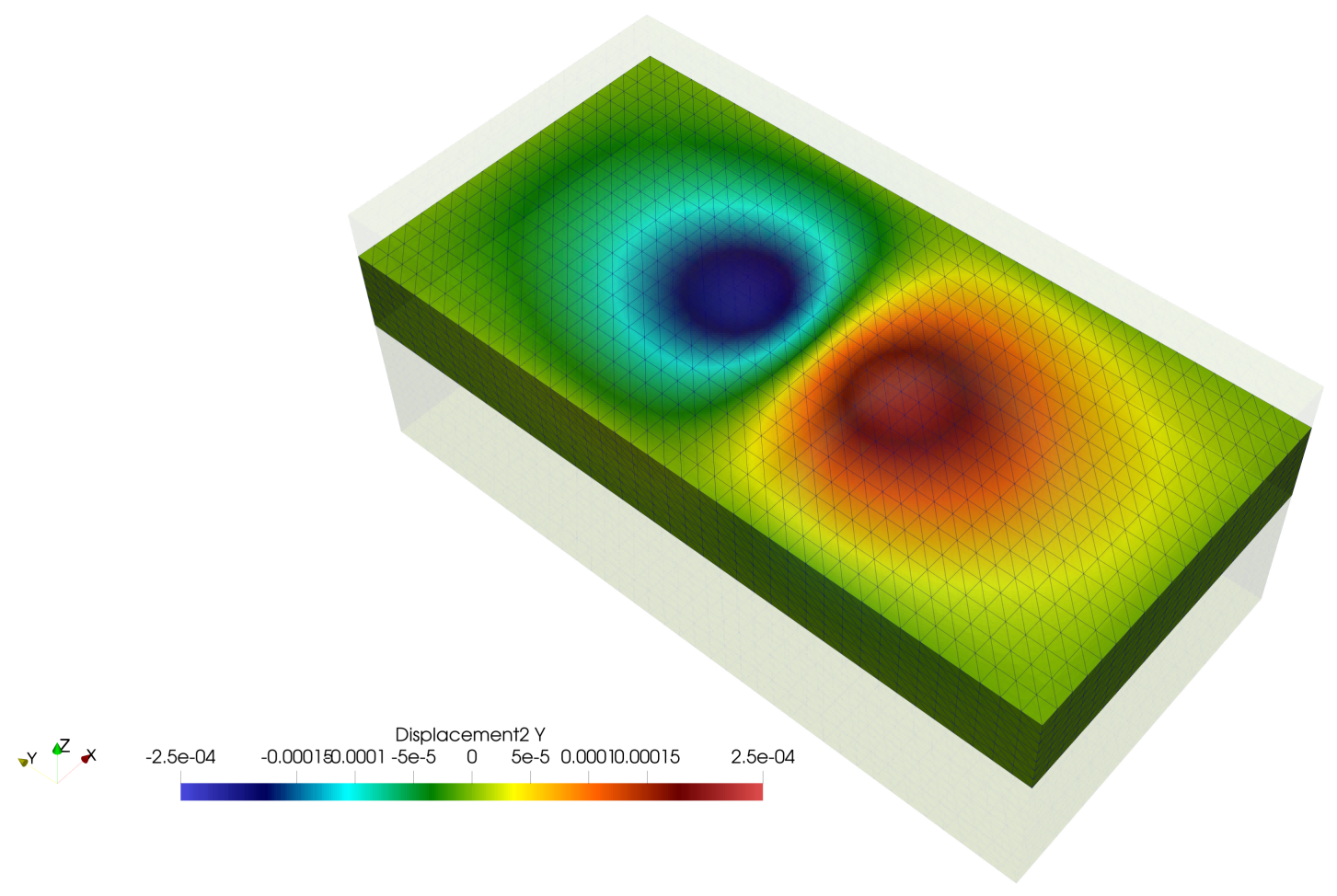}
\caption{Displacement field along the Y-axis of the $\Omega^2$ under $\boldsymbol{f}_{2}^{2}$.}
\label{3:fig15}
\end{minipage}
\begin{minipage}{0.3\linewidth}
\centering
\includegraphics[width=0.9\linewidth]{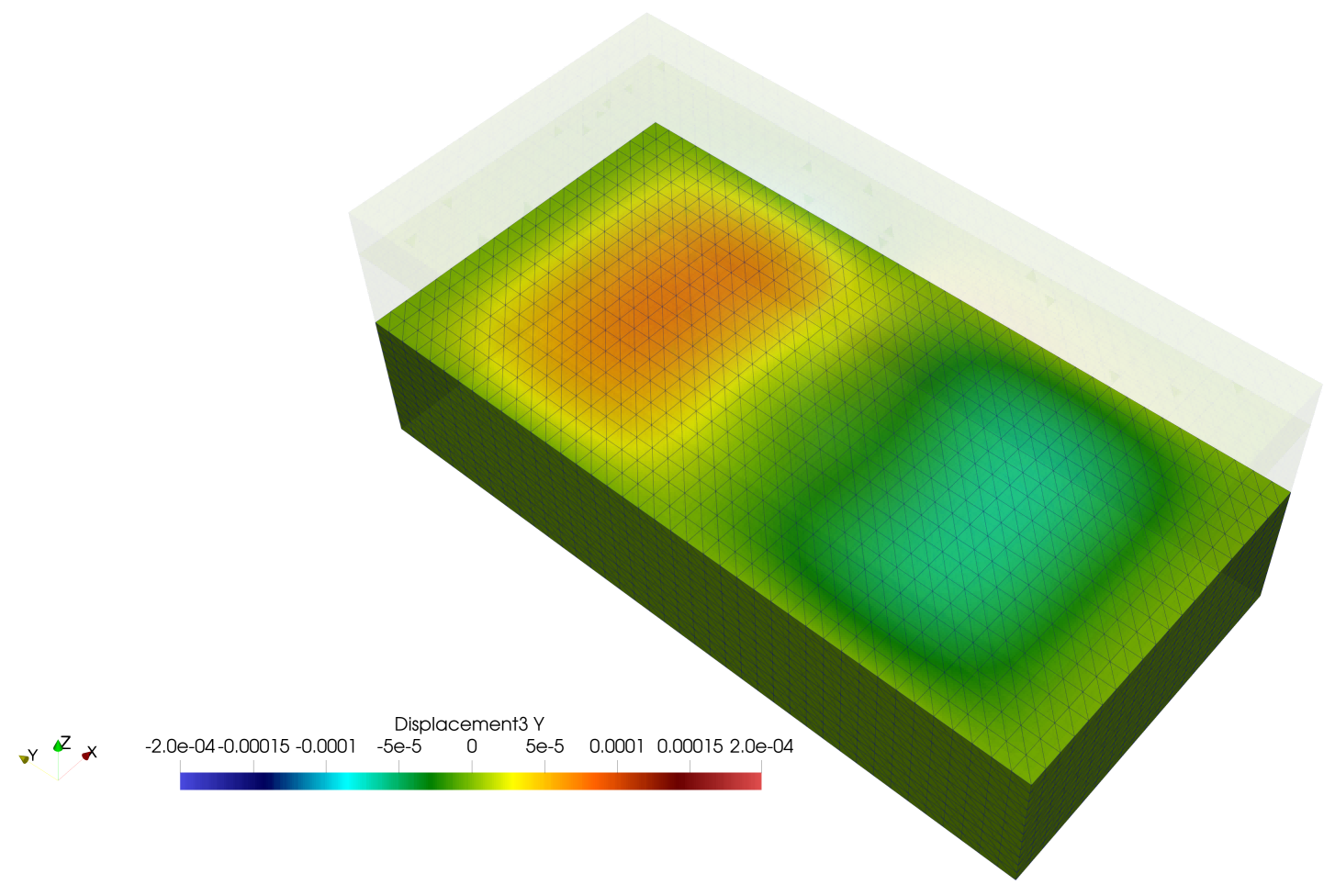}
\caption{Displacement field along the Y-axis of the $\Omega^3$ under $\boldsymbol{f}_{2}^{2}$.}
\label{3:fig16}
\end{minipage}
\end{figure}

\section{Numerical experiments}

In this section, a three-layer elastic system will be constructed in which a pavement mechanical model is simulated by setting the Young modulus and Poisson ratio for each layer. And two kinds of stress boundary conditions are applied to the model respectively to analyze the displacement response of the elastic system under different pavement stress states. The FreeFEM++ finite element software platform was used in this numerical experiment, and the Interior Point OPTimizer algorithm based on the internal point method was used to solve the inner-loop optimization problems \cite{3hecht2012new,3wachter2006implementation}. The following numerical experiments are performed on $\operatorname{Intel}(\mathrm{R})$ Core(TM) i7-11800H CPU, $2.30 \mathrm{GHz}$ with $32.0 \mathrm{~GB}$ RAM.

First, the three-dimension model of the three-layer elastic system is shown in Fig.\ref{3:fig2}. The system consist of three layers, which are denoted from top to bottom as $\Omega_1$, $\Omega_2$ and $\Omega_3$, respectively. Their definition in three-dimensional Euclidean space are as follow:
$\Omega^1 = (0,3) \times (0,6) \times (1.9, 2.3)$,
$\Omega^2 = (0,3) \times (0,6) \times (1.2, 1.9)$,
$\Omega^3 = (0,3) \times (0,6) \times (0, 1.2)$,
where the unit of length here is set to meters. Then, the material of each layer of elastic bodies is assumed to be homogeneous and isotropic, but with different modulus. Therefore, according to the theory of elasticity, the linear elasticity operator is defined as follows:
$$
\sigma_{k l}^i=\frac{E^i P_o^i}{\left(1+P_o^i\right)\left(1-2 P_o^i\right)}\left(\nabla \cdot \mathbf{u}^i\right) \delta_{k l}+\frac{E^i}{1+P_o^i} \varepsilon_{k l}\left(\mathbf{u}^i\right), i=1,2,3,
$$
According to the existing experimental research on pavement materials \cite{3ma2021analytical,3kim2011numerical}, the Young's modulus $E^{i}$ and Poisson ratio $P_o^{i}$ of each material in the three-layer elastic system are respectively set to:
$E^1 = 5\cdot 10^3$, $P_o^{1} = 0.25$,
$E^2 = 2\cdot 10^3$, $P_o^{2} = 0.25$,
$E^3 = 2\cdot 10^2$, $P_o^{3} = 0.4$.

The specific mesh generation is shown in Fig.\ref{3:fig2}, where $h=0.15$. According to the mesh on the multi-layer system, it can be judged that there are a total of $22386$ mesh nodes in the numerical experiment and the degree of freedom of the finite element space is $67158$.
The contact functions $g^1(x)$ and $g^2(x)$ at the contact surfaces $\Gamma_{c}^{1}$ and $\Gamma_{c}^{2}$ are defined as $g^1(x) = 0.2$ on $\Gamma_{c}^{1}$ and
$g^2(x) = 0.05$ on $\Gamma_{c}^{2}$.
In addition, in this numerical experiment, the parameters in Algorithm \ref{algorithm:3.2} are set as follows: the tolerance error $tol=10^{-4}$, the parameter $\theta = 0.04$ and the initial vector $\boldsymbol{\lambda}_{0h} = \boldsymbol{0}$.
Then, two stress conditions will be introduced and analyzed respectively.

\begin{figure}[!t]
\centering
\begin{minipage}{0.24\linewidth}
\centering
\includegraphics[width=0.9\linewidth]{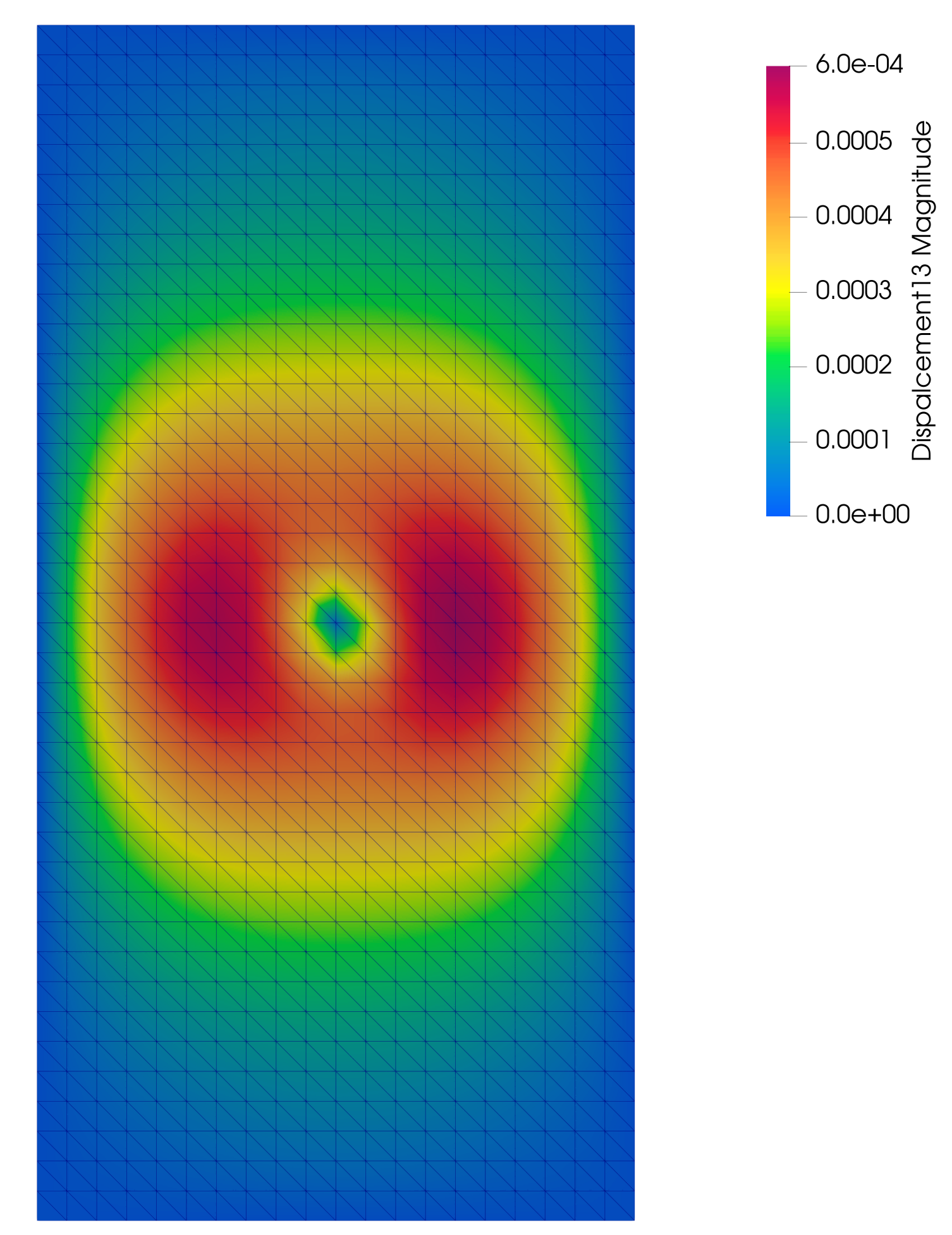}
\caption{The difference between the displacement of $\Omega^1$ and $\Omega^2$ in the contact zone $\Gamma^1_{c}$ under $\boldsymbol{f}_{2}^{1}$.}
\label{3:fig9}
\end{minipage}
\begin{minipage}{0.24\linewidth}
\centering
\includegraphics[width=0.9\linewidth]{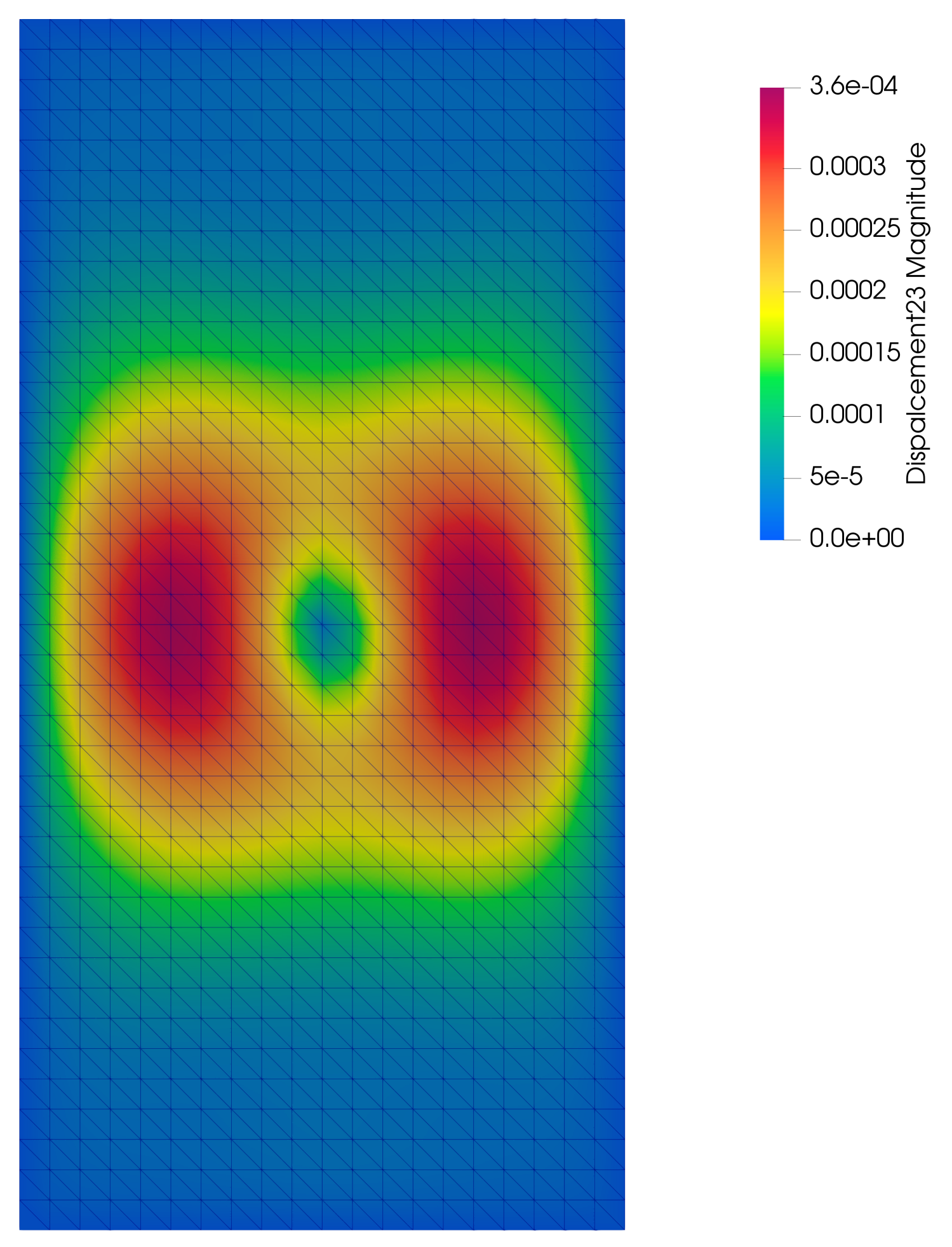}
\caption{The difference between the displacement of $\Omega^2$ and $\Omega^3$ in the contact zone $\Gamma^2_{c}$ under $\boldsymbol{f}_{2}^{1}$.}
\label{3:fig10}
\end{minipage}
\begin{minipage}{0.24\linewidth}
\centering
\includegraphics[width=0.9\linewidth]{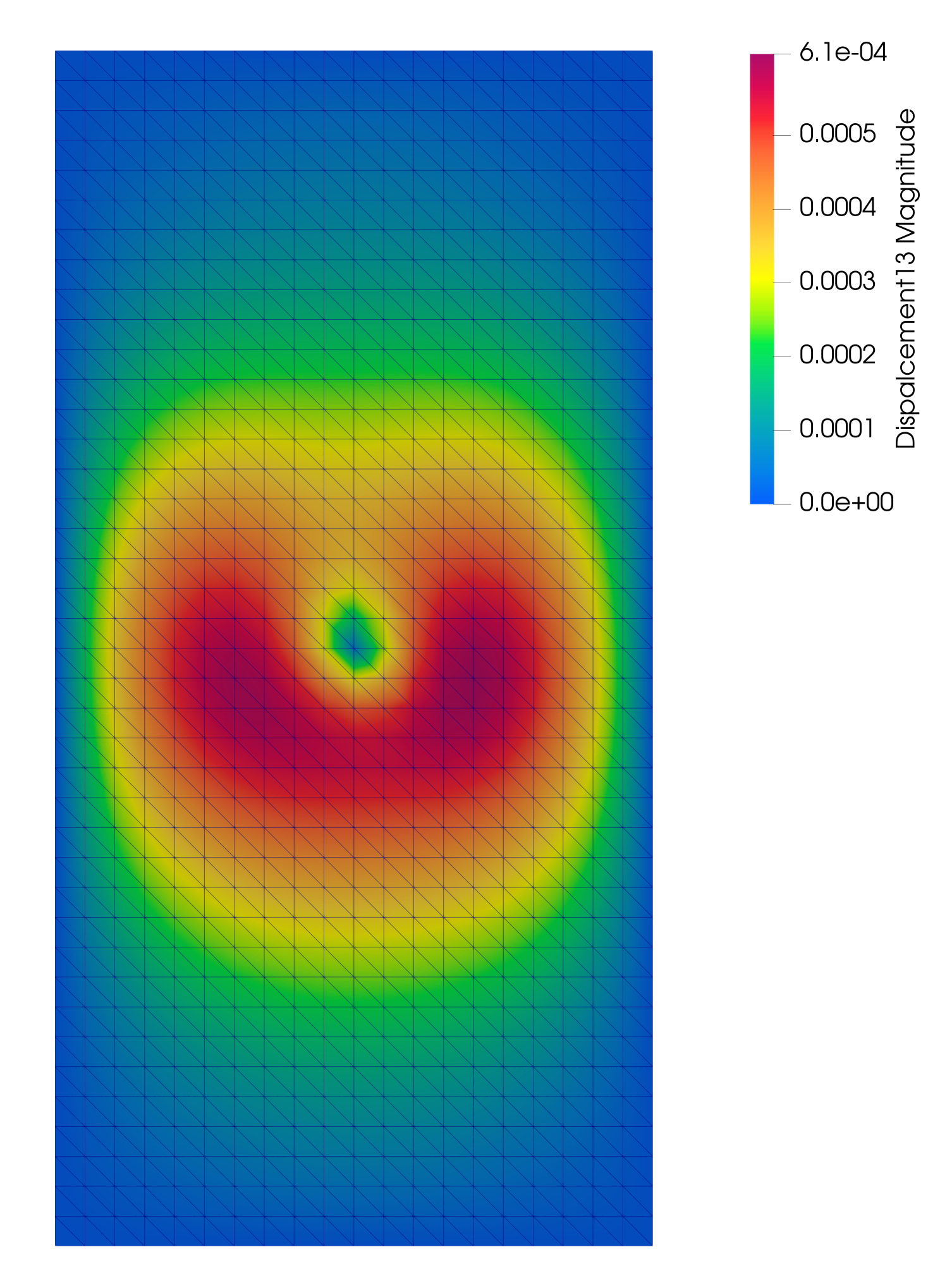}
\caption{The difference between the displacement of $\Omega^1$ and $\Omega^2$ in the contact zone $\Gamma^1_{c}$ under $\boldsymbol{f}_{2}^{2}$.}
\label{3:fig17}
\end{minipage}
\begin{minipage}{0.24\linewidth}
\centering
\includegraphics[width=0.9\linewidth]{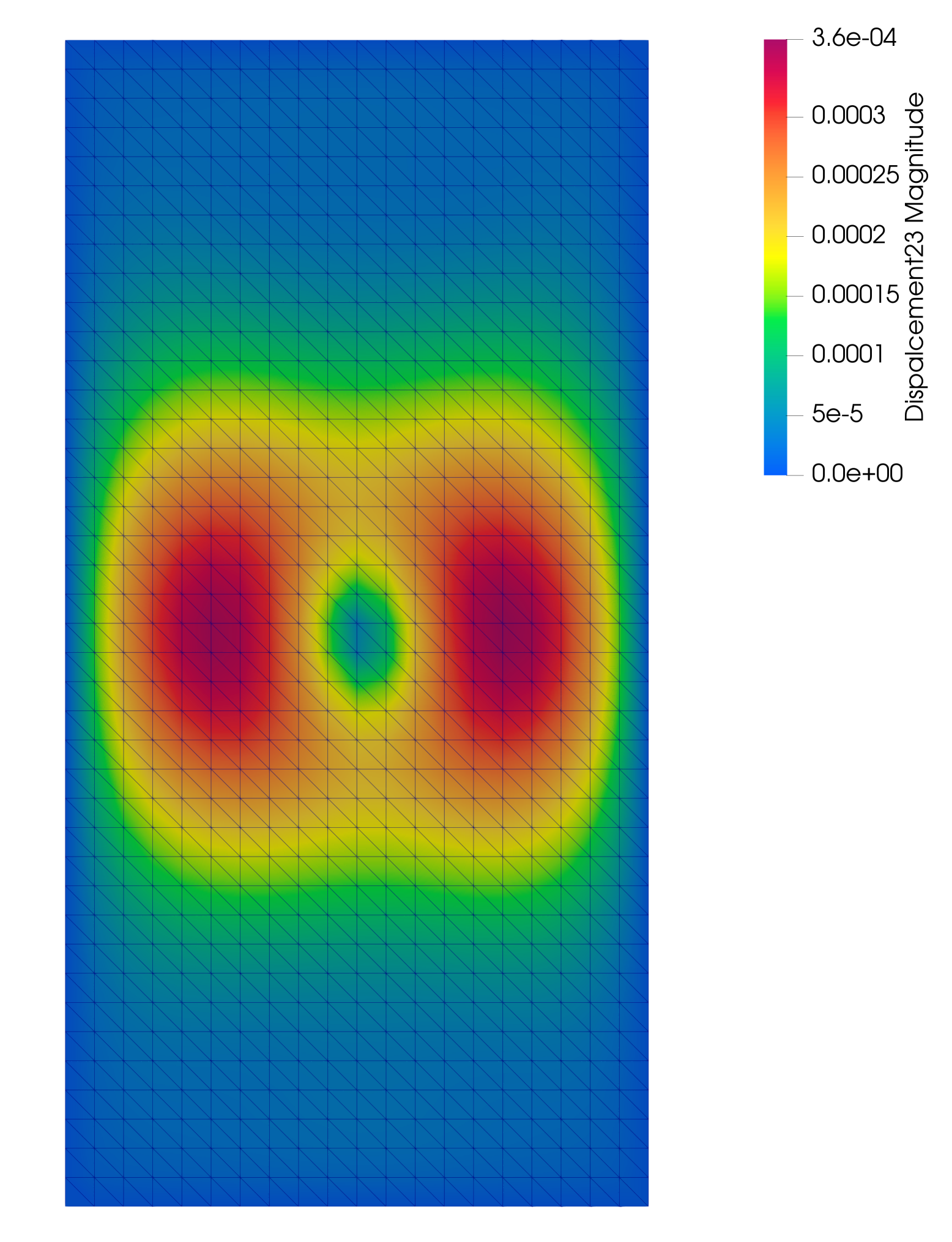}
\caption{The difference between the displacement of $\Omega^1$ and $\Omega^2$ in the contact zone $\Gamma^1_{c}$ under $\boldsymbol{f}_{2}^{2}$.}
\label{3:fig18}
\end{minipage}
\end{figure}

\subsection{Vertical stress boundary condition}

In the first example, the three-layer elastic system is subjected to gravity $\boldsymbol{f}_{0}= [0,0,-0.05]^{\top}$, and the surface $\Gamma_2^{1}$ of $\Omega^1$ will be subjected to vertical downward pressure $\boldsymbol{f}_{2}^{1}$, which is defined as
\begin{small}
$$
\boldsymbol{f}_{2}^{1} = \left\{ 
\begin{array}{ll}
[0,0,-22.5]^{\top}, & x\in[1.34, 1.66],~ y\in[2.84,3.16],~ z=2.3,\\
{[0,0,0]}^{\top} , & \text{other}.
\end{array}
\right.
$$
\end{small}
The mechanical units of force $\boldsymbol{f}_{0}$ and surface force $\boldsymbol{f}_{2}$ here are $N/m^{3}$ and $N/m^2$, respectively. The stress boundary condition simulates the stress state of the road surface when the vehicle is stationary or traveling at a constant speed. 

Under this stress boundary condition, the displacement function calculation results of the three-layer elastic system are shown in Fig.\ref{3:fig3}-\ref{3:fig8}. Among them, Fig.\ref{3:fig3}-\ref{3:fig5} shows that under vertical stress, the displacements of the elastic layers are symmetric along the X-axis or Y-axis, and Fig.\ref{3:fig6}-\ref{3:fig8} can also prove this. Fig.\ref{3:fig9} and Fig.\ref{3:fig10} show the difference in displacement between the two layers of elastic system at the contact zone $\Gamma_c^{1}$ and $\Gamma_c^{2}$, respectively. Obviously, if the displacement difference is $0$, then the contact surface is in a sticking state, if the displacement difference is not $0$, then it is in a sliding state. It can be confirmed from the images that compared with $\Gamma_c^{2}$, the area in the sliding state on $\Gamma_c^{1}$ is wider and the slip amplitude is larger.

\subsection{Non-vertical stress boundary conditions}

In order to simulate the mechanical model of the pavement in the case of vehicle acceleration or deceleration, the force $\boldsymbol{f}_{2}^{2}$ on the boundary $\Gamma_{2}^{1}$ is defined as follows:
\begin{small}
$$
\boldsymbol{f}_{2}^{2} = \left\{ 
\begin{array}{ll}
[0,-4.5,-22.5]^{\top}, & x\in[1.34, 1.66],~ y\in[2.84,3.16],~ z=2.3,\\
{[0,0,0]}^{\top} , & \text{other}.
\end{array}
\right.
$$
\end{small}
The volume force on the three-layer elastic system is still set to $\boldsymbol{f}_0=[0,0,-0.05]^{\top}$. 

Under the action of stress boundary condition $\boldsymbol{f}_{2}^{2}$, the three-layer elastic system will displace, and the displacement results are shown in Fig.\ref{3:fig11}-\ref{3:fig16}. It can be seen from Fig.\ref{3:fig14}-\ref{3:fig16} that under the action of non-vertical surface forces $\boldsymbol{f}_{2}^{2}$, the displacement field along the Y-axis is not symmetric with respect to the X-axis.
Similarly, under the action of external force $\boldsymbol{f}_{2}^{2}$, the difference between the interlayer displacements on the contact surface $\Gamma_c^1$ and $\Gamma_c^1$ shifts significantly to the negative direction of the Y-axis, as can be proved in Fig.\ref{3:fig17} and Fig.\ref{3:fig18}.

\begin{figure}[t]
\centering
\begin{minipage}{0.495\linewidth}
\centering
\includegraphics[width=1\linewidth]{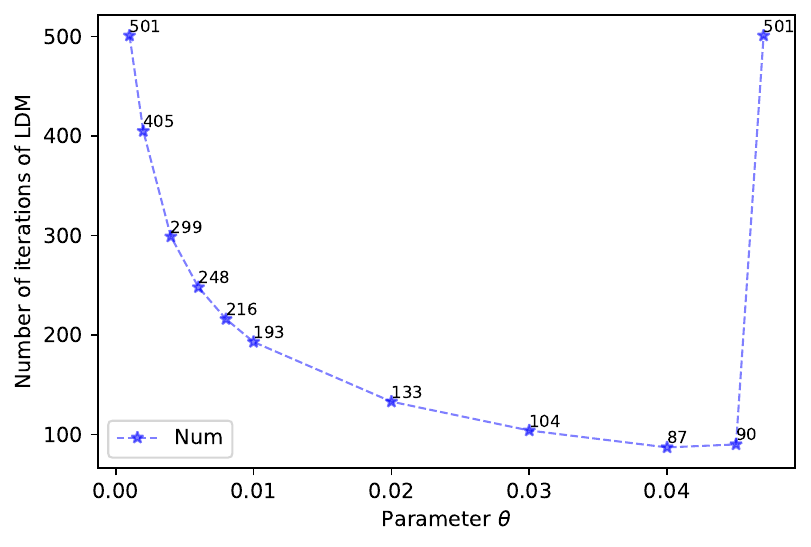}
\caption{The relationship between the number of iterations of the LD algorithm and the parameter $\theta$ with $h=0.15$.}
\label{3:fig:Num_iterations_LDM}
\end{minipage}
\begin{minipage}{0.495\linewidth}
\centering
\includegraphics[width=1\linewidth]{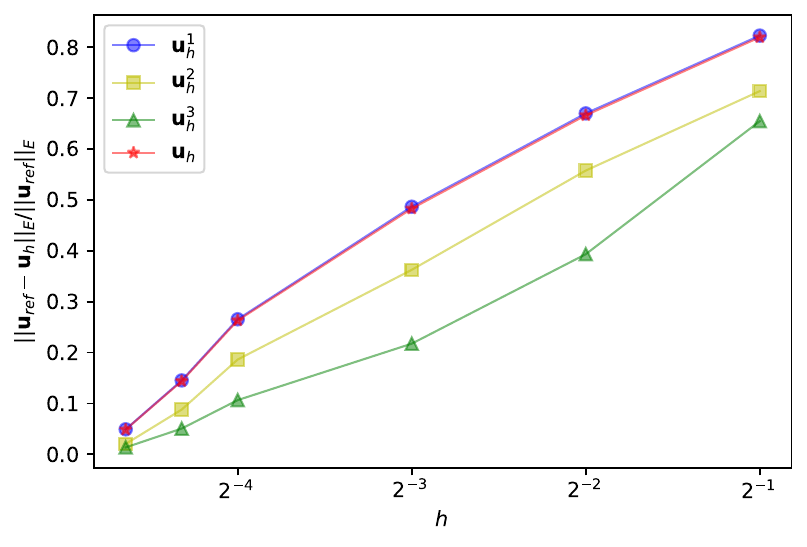}
\caption{Relative Numerical errors: $\|\boldsymbol{u}_{h}^{i}-\boldsymbol{u}_{r}^{i}\|_{E}/\|\boldsymbol{u}_{r}^{i}\|_{E}$ and $\|\boldsymbol{u}_{h}-\boldsymbol{u}_{r}\|_{E}/\|\boldsymbol{u}_{r}\|_{E}$ based on LD method.}
\label{3:fig:Num_err_LDM}
\end{minipage}
\end{figure}

\subsection{Verification of convergence}

According to Theorem \ref{thm:3.2}, it can be found that two independent convergence properties of the discrete LD algorithm \ref{algorithm:3.2} will affect the numerical results.
\begin{itemize}
\item The parameter $\theta$ will determine whether the algorithm converges and affects the number of iterations required for the algorithm to converge to a given tolerance error;
\item The mesh diameter $h$ determines the error between the finite element numerical solution and the exact solution. When the number of iterations is sufficiently large, the error is determined only by $h$.
\end{itemize}

First, in order to verify the influence of parameter $\theta$ on the convergence of the algorithm, multiple parameters were selected, and the number of iterations required for the algorithm to converge to the tolerance error $tol$ when using this parameter was recorded. The results are shown in the Fig.\ref{3:fig:Num_iterations_LDM}. 
It can be found that by increasing the parameter $\theta$, the convergence rate of Algorithm \ref{algorithm:3.2} is accelerated, but when $\theta$ is greater than the upper bound, the algorithm will no longer converge, which is consistent with the conclusions of Theorem \ref{thm:3.1} and Theorem \ref{thm:3.2}.

Then, the convergence properties of the numerical solution with respect to diameter $h$ under this LD algorithm need to be verified through numerical experiments. For this purpose, we define the energy norm as follows:
$$
\left\|\boldsymbol{v}_h^i\right\|_E=\frac{1}{\sqrt{2}}\left(a^i\left(\boldsymbol{v}_h^i, \boldsymbol{v}_h^i\right)\right)^{1 / 2}, \quad\left\|\boldsymbol{v}_h\right\|_E=\frac{1}{\sqrt{2}}\left(a\left(\boldsymbol{v}_h, \boldsymbol{v}_h\right)\right)^{1 / 2}
$$
Under this fixed tolerance error $tol$, the numerical solution with $h=1/32$ is used as the reference solution $\boldsymbol{u}_{ref}$, and the relative errors $\|\boldsymbol{u}_h - \boldsymbol{u}_{ref}\|_{E}/ \|\boldsymbol{u}_{ref}\|_{E}$ and $\|\boldsymbol{u}^i_h - \boldsymbol{u}^i_{ref}\|_{E}/ \|\boldsymbol{u}^i_{ref}\|_{E}$ ($i=1,2,3$) of the numerical solution under different mesh diameters $h$ is calculated, and the results are shown in Fig.\ref{3:fig:Num_err_LDM}.
This numerical experiment also proves that when $tol$ is sufficiently small, the finite element numerical solution calculated by Algorithm \ref{algorithm:3.2} can converge to the exact solution as $h\to0$. 

\section{Conclusion}

In this paper, a layer decomposition algorithm for solving multi-layer elastic contact systems with interlayer Tresca friction contact conditions is studied. 
In this algorithm, based on the Neumann-Neumann method, the original coupled multi-layer contact problem is decomposed into multiple unilateral contact sub-problems. Different from the LATIN method \cite{oumaziz2017non, giacoma2015toward}, no additional assumptions are given, so the algorithm has wide applicability in theory. In addition, based on the independence of these sub-problems, the LD algorithm has the possibility to be combined with the FETI method \cite{dostal2010scalable, dostal2012theoretically}, so that the original problem can be decomposed into more sub-problems and solved. This combination algorithm needs further exploration and research.

It is also worth emphasizing that the application background of layered elastic systems with interlayer friction conditions is not only in pavement mechanics modeling. There are many potential applications in composite layered structures, such as spliced wood, rail and so on. Finally, there are important derived models and algorithms based on the multi-layer elastic contact system worth exploring and studying. For example, dynamic models involving the propagation of forces, viscoelastic constitutive equation models involving viscoelastic materials, semi-variational inequality models involving complex interlayer contact conditions, and parallel computational methods for layer decomposition algorithms, etc. The expansion of such models and algorithms will enrich the application scenarios and improve accuracy of multi-layer elastic contact systems, and also provide a new direction for the study of variational inequality theory.

%% file: LDM_for_MLECS_Elsarticle.bbl
\begin{thebibliography}{10}

\bibitem{3abide2022unified}
St{\'e}phane Abide, Mika{\"e}l Barboteu, Soufiane Cherkaoui, and Serge Dumont.
\newblock Unified primal-dual active set method for dynamic frictional contact problems.
\newblock {\em Fixed Point Theory and Algorithms for Sciences and Engineering}, 2022(1):1--22, 2022.

\bibitem{3adams2003sobolev}
Robert~A Adams and John~JF Fournier.
\newblock {\em Sobolev spaces}.
\newblock Elsevier, 2003.

\bibitem{3bayada2008convergence}
Guy Bayada, Jalila Sabil, and Taoufik Sassi.
\newblock Convergence of a neumann-dirichlet algorithm for two-body contact problems with non local coulomb's friction law.
\newblock {\em ESAIM: Mathematical Modelling and Numerical Analysis}, 42(2):243--262, 2008.

\bibitem{3BEAUDE2023116124}
L.~Beaude, F.~Chouly, M.~Laaziri, and R.~Masson.
\newblock Mixed and nitsche’s discretizations of coulomb frictional contact-mechanics for mixed dimensional poromechanical models.
\newblock {\em Computer Methods in Applied Mechanics and Engineering}, 413:116124, 2023.

\bibitem{belgacem1998mortar}
F~Ben Belgacem, P~Hild, and P~Laborde.
\newblock The mortar finite element method for contact problems.
\newblock {\em Mathematical and Computer Modelling}, 28(4-8):263--271, 1998.

\bibitem{bjorstad1986iterative}
Petter~E Bj{\o}rstad and Olof~B Widlund.
\newblock Iterative methods for the solution of elliptic problems on regions partitioned into substructures.
\newblock {\em SIAM Journal on Numerical Analysis}, 23(6):1097--1120, 1986.

\bibitem{3brown1996soil}
SF~Brown.
\newblock Soil mechanics in pavement engineering.
\newblock {\em G{\'e}otechnique}, 46(3):383--426, 1996.

\bibitem{3brownjohn2007structural}
James~MW Brownjohn.
\newblock Structural health monitoring of civil infrastructure.
\newblock {\em Philosophical Transactions of the Royal Society A: Mathematical, Physical and Engineering Sciences}, 365(1851):589--622, 2007.

\bibitem{3burmister1945general}
Donald~M Burmister.
\newblock The general theory of stresses and displacements in layered systems. i.
\newblock {\em Journal of Applied Physics}, 16(2):89--94, 1945.

\bibitem{chouly2023finite}
Franz Chouly, Patrick Hild, and Yves Renard.
\newblock {\em Finite element approximation of contact and friction in elasticity}.
\newblock Springer, 2023.

\bibitem{dostal2010scalable}
Zdenek Dost{\'a}l, Tom{\'a}s Kozubek, Petr Horyl, Tom{\'a}s Brzobohat{\`y}, and Alexandros Markopoulos.
\newblock A scalable tfeti algorithm for two-dimensional multibody contact problems with friction.
\newblock {\em Journal of computational and applied mathematics}, 235(2):403--418, 2010.

\bibitem{dostal2012theoretically}
Zden{\v{e}}k Dost{\'a}l, Tom{\'a}{\v{s}} Kozubek, Alexandros Markopoulos, Tom{\'a}{\v{s}} Brzobohat{\`y}, V{\'\i}t Vondr{\'a}k, and Petr Horyl.
\newblock A theoretically supported scalable tfeti algorithm for the solution of multibody 3d contact problems with friction.
\newblock {\em Computer methods in applied mechanics and engineering}, 205:110--120, 2012.

\bibitem{drouet2017accurate}
Guillaume Drouet and Patrick Hild.
\newblock An accurate local average contact method for nonmatching meshes.
\newblock {\em Numerische Mathematik}, 136:467--502, 2017.

\bibitem{3eck1998existence}
Christof Eck and Jii'i Jarusek.
\newblock Existence results for the static contact problem with coulomb friction.
\newblock {\em Mathematical Models and Methods in Applied Sciences}, 8(03):445--468, 1998.

\bibitem{fernandez2003numerical}
Jos{\'e}~R Fern{\'a}ndez, Patrick Hild, and Juan~M Viano.
\newblock Numerical approximation of the elastic-viscoplastic contact problem with non-matching meshes.
\newblock {\em Numerische Mathematik}, 94:501--522, 2003.

\bibitem{3FRANCESCHINI2022114632}
Andrea Franceschini, Massimiliano Ferronato, Matteo Frigo, and Carlo Janna.
\newblock A reverse augmented constraint preconditioner for lagrange multiplier methods in contact mechanics.
\newblock {\em Computer Methods in Applied Mechanics and Engineering}, 392:114632, 2022.

\bibitem{giacoma2015toward}
Anthony Giacoma, David Dureisseix, Anthony Gravouil, and Michel Rochette.
\newblock Toward an optimal a priori reduced basis strategy for frictional contact problems with latin solver.
\newblock {\em Computer Methods in Applied Mechanics and Engineering}, 283:1357--1381, 2015.

\bibitem{3han2019numerical}
Weimin Han and Mircea Sofonea.
\newblock Numerical analysis of hemivariational inequalities in contact mechanics.
\newblock {\em Acta Numerica}, 28:175--286, 2019.

\bibitem{3haslinger2014domain}
Jaroslav Haslinger, Radek Ku{\v{c}}era, Julien Riton, and Taoufik Sassi.
\newblock A domain decomposition method for two-body contact problems with tresca friction.
\newblock {\em Advances in Computational Mathematics}, 40(1):65--90, 2014.

\bibitem{3hecht2012new}
Fr{\'e}d{\'e}ric Hecht.
\newblock New development in freefem++.
\newblock {\em Journal of Numerical Mathematics}, 20(3-4):251--266, 2012.

\bibitem{hild2000numerical}
Patrick Hild.
\newblock Numerical implementation of two nonconforming finite element methods for unilateral contact.
\newblock {\em Computer Methods in Applied Mechanics and Engineering}, 184(1):99--123, 2000.

\bibitem{3hild2010stabilized}
Patrick Hild and Yves Renard.
\newblock A stabilized lagrange multiplier method for the finite element approximation of contact problems in elastostatics.
\newblock {\em Numerische Mathematik}, 115:101--129, 2010.

\bibitem{3huang2004pavement}
Yang~Hsien Huang et~al.
\newblock {\em Pavement analysis and design}, volume~2.
\newblock Pearson Prentice Hall Upper Saddle River, NJ, 2004.

\bibitem{3hueber2008primal}
Stefan H{\"u}eber, Georg Stadler, and Barbara~I Wohlmuth.
\newblock A primal-dual active set algorithm for three-dimensional contact problems with coulomb friction.
\newblock {\em SIAM Journal on Scientific Computing}, 30(2):572--596, 2008.

\bibitem{3hung2001elastic}
H-H Hung and Y-B Yang.
\newblock Elastic waves in visco-elastic half-space generated by various vehicle loads.
\newblock {\em Soil Dynamics and Earthquake Engineering}, 21(1):1--17, 2001.

\bibitem{3ito2003semi}
Kazufumi Ito and Karl Kunisch.
\newblock Semi--smooth newton methods for variational inequalities of the first kind.
\newblock {\em ESAIM: Mathematical Modelling and Numerical Analysis}, 37(1):41--62, 2003.

\bibitem{3kikuchi1988contact}
Noboru Kikuchi and John~Tinsley Oden.
\newblock {\em Contact problems in elasticity: a study of variational inequalities and finite element methods}.
\newblock SIAM, 1988.

\bibitem{3kim2011numerical}
Hyunwook Kim, Martin Arraigada, Christiane Raab, and Manfred~N Partl.
\newblock Numerical and experimental analysis for the interlayer behavior of double-layered asphalt pavement specimens.
\newblock {\em Journal of Materials in Civil Engineering}, 23(1):12--20, 2011.

\bibitem{3krause2009nonsmooth}
Rolf Krause.
\newblock A nonsmooth multiscale method for solving frictional two-body contact problems in 2d and 3d with multigrid efficiency.
\newblock {\em SIAM Journal on Scientific Computing}, 31(2):1399--1423, 2009.

\bibitem{3kuvcera2013interior}
Radek Ku{\v{c}}era, Jitka Machalov{\'a}, Horym{\'\i}r Netuka, and Pavel {\v{Z}}en{\v{c}}{\'a}k.
\newblock An interior-point algorithm for the minimization arising from 3d contact problems with friction.
\newblock {\em Optimization Methods and Software}, 28(6):1195--1217, 2013.

\bibitem{3laborde2008fixed}
Patrick Laborde and Yves Renard.
\newblock Fixed point strategies for elastostatic frictional contact problems.
\newblock {\em Mathematical Methods in the Applied Sciences}, 31(4):415--441, 2008.

\bibitem{3liu2011coupled}
Zhen Liu and Xiong Yu.
\newblock Coupled thermo-hydro-mechanical model for porous materials under frost action: theory and implementation.
\newblock {\em Acta Geotechnica}, 6:51--65, 2011.

\bibitem{3ma2021analytical}
Xianyong Ma, Weiwen Quan, Chundi Si, Zejiao Dong, and Yongkang Dong.
\newblock Analytical solution for the mechanical responses of transversely isotropic viscoelastic multi-layered asphalt pavement subjected to moving harmonic load.
\newblock {\em Soil Dynamics and Earthquake Engineering}, 147:106822, 2021.

\bibitem{3nawaz2013soil}
Muhammad~Farrakh Nawaz, Guilhem Bourrie, and Fabienne Trolard.
\newblock Soil compaction impact and modelling. a review.
\newblock {\em Agronomy for Sustainable Development}, 33:291--309, 2013.

\bibitem{oumaziz2017non}
Paul Oumaziz, Pierre Gosselet, Pierre-Alain Boucard, and St{\'e}phane Guinard.
\newblock A non-invasive implementation of a mixed domain decomposition method for frictional contact problems.
\newblock {\em Computational Mechanics}, 60:797--812, 2017.

\bibitem{3raposeiras2012influence}
Aitor~C Raposeiras, {\'A}ngel Vega-Zamanillo, Miguel~{\'A}ngel Calzada-P{\'e}rez, and Daniel Castro-Fresno.
\newblock Influence of surface macro-texture and binder dosage on the adhesion between bituminous pavement layers.
\newblock {\em Construction and Building Materials}, 28(1):187--192, 2012.

\bibitem{3simo1992augmented}
J~Ci Simo and TA1143885 Laursen.
\newblock An augmented lagrangian treatment of contact problems involving friction.
\newblock {\em Computers \& Structures}, 42(1):97--116, 1992.

\bibitem{3wachter2006implementation}
Andreas W{\"a}chter and Lorenz~T Biegler.
\newblock On the implementation of an interior-point filter line-search algorithm for large-scale nonlinear programming.
\newblock {\em Mathematical Programming}, 106:25--57, 2006.

\bibitem{wohlmuth2012abstract}
BI~Wohlmuth, A~Popp, MW~Gee, and WA29351821312 Wall.
\newblock An abstract framework for a priori estimates for contact problems in 3d with quadratic finite elements.
\newblock {\em Computational Mechanics}, 49(6):735--747, 2012.

\bibitem{3yusoff2011modelling}
Nur Izzi~Md Yusoff, Montgomery~T Shaw, and Gordon~D Airey.
\newblock Modelling the linear viscoelastic rheological properties of bituminous binders.
\newblock {\em Construction and Building Materials}, 25(5):2171--2189, 2011.

\bibitem{3zang2011contact}
Mengyan Zang, Wei Gao, and Zhou Lei.
\newblock A contact algorithm for 3d discrete and finite element contact problems based on penalty function method.
\newblock {\em Computational Mechanics}, 48:541--550, 2011.

\bibitem{3zhang2016weak}
Yuqing Zhang, Bjorn Birgisson, and Robert~L Lytton.
\newblock Weak form equation--based finite-element modeling of viscoelastic asphalt mixtures.
\newblock {\em Journal of Materials in Civil Engineering}, 28(2):04015115, 2016.

\bibitem{3zhang2022variational}
Zhizhuo Zhang, Xiaobing Nie, and Jinde Cao.
\newblock Variational inequalities of multilayer elastic systems with interlayer friction: existence and uniqueness of solution and convergence of numerical solution.
\newblock {\em arXiv preprint arXiv:2207.13260}, 2022.

\bibitem{3zhangvariational}
Zhizhuo Zhang, Xiaobing Nie, and Jinde Cao.
\newblock Variational inequalities of multilayer viscoelastic systems with interlayer tresca friction: Existence and uniqueness of solution and convergence of numerical solution.
\newblock {\em Mathematical Methods in the Applied Sciences}, 47(2):1170--1194, 2024.

\bibitem{3zokaei2014finite}
Ali Zokaei-Ashtiani, Cesar Carrasco, and Soheil Nazarian.
\newblock Finite element modeling of slab--foundation interaction on rigid pavement applications.
\newblock {\em Computers and Geotechnics}, 62:118--127, 2014.

\end{thebibliography}
